\documentclass[11pt]{article} 
\usepackage{amsfonts,latexsym,amstext} 
\usepackage{amsmath} 
\usepackage{amssymb} 
\usepackage{comment} 
\usepackage{amsthm} 
\usepackage[latin1]{inputenc} 
\usepackage{color} 
 
\theoremstyle{plain} 
\newtheorem{theorem}{Theorem}[section] 
\newtheorem{lemma}[theorem]{Lemma} 
\newtheorem{proposition}[theorem]{Proposition} 
\newtheorem{corollary}[theorem]{Corollary} 

\newtheorem{definition}[theorem]{Definition} 
\newtheorem{hypothesis}[theorem]{Hypothesis} 
\theoremstyle{remark} 
\newtheorem{remark}[theorem]{Remark} 
\newtheorem{example}[theorem]{Example}

\allowdisplaybreaks 
 
\makeatletter 
\@addtoreset{equation}{section} 
\def\theequation{\thesection.\arabic{equation}} 
\makeatother 
 
\def\qed{{\hfill\hbox{\enspace${ \square}$}} \smallskip} 
\def\sqr#1#2{{\vcenter{\vbox{\hrule height .#2pt \hbox{\vrule 
 width .#2pt height#1pt \kern#1pt \vrule 
width .#2pt} \hrule height .#2pt}}}} 
\def\square{\mathchoice\sqr54\sqr54\sqr{4.1}3\sqr{3.5}3}

\def\ds{\begin{displaystyle}} 
\def\eds{\end{displaystyle}} 
 
\def\<{\langle } 
\def\>{\rangle }

\def\R{\mathbb R} 
\def\N{\mathbb N}

\def\P{\mathbb P} 
\def\Q{\mathbb Q}

\def\calf{{\cal F}}

\newcommand{\sper}[1]{\mathbb{E} \left[ #1 \right]}                               
 
\DeclareMathAlphabet{\mathonebb}{U}{bbold}{m}{n}                           %
\newcommand{\one}{\ensuremath{\mathonebb{1}}}                               
\newcommand{\limit}{\operatornamewithlimits{\longrightarrow}}   

\title{Weak Dirichlet processes with jumps 
} 

\usepackage{authblk}

\author[1]{Elena Bandini\thanks{ebandini@luiss.it}
} 
\author[2]{Francesco Russo\thanks{francesco.russo@ensta-paristech.fr}} 
\affil[1]{LUISS Roma, Dipartimento di Economia e Finanza, via Romania 32, I-00197 Roma, Italy} 
\affil[2]{ENSTA ParisTech, Universit\'e Paris-Saclay, 
 Unit\'e de Math\'ematiques appliqu\'ees, 828, Boulevard des Mar\'echaux, F-91120 Palaiseau, France}

\date{} 
 
\begin{document} 
 
\allowdisplaybreaks 
\maketitle

\begin{abstract} 
This paper develops systematically the stochastic calculus via regularization 
in the case of jump processes. In particular one continues  the analysis  
of real-valued c\`adl\`ag weak Dirichlet processes with respect to a 
given filtration. Such a process is the sum of a local martingale 
and an adapted process $A$ such that $[N,A] = 0$, for any  
continuous local martingale $N$. 
Given a  
function $u:[0,T] \times \R \rightarrow \R$, which is of  
class $C^{0,1}$ (or sometimes less), we provide a chain 
rule type expansion for $u(t,X_t)$  which stands in applications 
for a chain It\^o type rule.

\end{abstract} 
{\bf Key words:} Weak Dirichlet processes; Calculus via regularizations; 
Random measure; Stochastic integrals for jump processes; Orthogonality. 

{\small\textbf{MSC 2010:}  60J75; 60G57; 60H05} 
 
\section{Introduction} 
The present paper  extends  the stochastic calculus via regularizations 
to the case of jump processes, and  carries  on the investigations on 
the so called weak Dirichlet processes in the discontinuous case. 
Applications of that calculus appear in the companion paper \cite{BandiniRusso2}, where  we provide  
the identification of 
the solution of a forward 
backward stochastic  differential 
equation driven by a random measure.
 
 
Stochastic calculus via regularization was essentially known 
in the case of continuous integrators $X$, see e.g. 
\cite{RVCRAS, rv1}, with a survey in \cite{RVSem}. 
In this case a fairly complete 
theory was developed, 
see for instance It\^o formulae for processes with finite quadratic (and more general) variation, 
stochastic differential equations, It\^o-Wentzell type formulae  \cite{flru1}, and generalizations to the 
case of Banach space type integrators, see e.g. \cite{DGR1}. 
The notion of covariation $[X,Y]$ (resp. quadratic variation $[X,X]$) for two processes $X,Y$ (resp. a process $X$) has 
been introduced in the framework of regularizations (see \cite{rv95}) and of discretizations as well (see \cite{fo}). 
Even if there is no direct theorem relating the two approaches, those
 coincide in all the  examples considered in the literature. 
If $X$ is a finite quadratic variation continuous  process, 
an It\^o formula  was proved  for the expansion of $F(X_t)$, when $F \in C^2$, see \cite{rv95}; 
this constitutes the counterpart of the related result for discretizations, 
see \cite{fo}.  Moreover, for  $F$  of class 
$C^1$  and $X$  a reversible continuous semimartingale, an It\^o expansion has been established in  \cite{rv96}. 
 
A natural extension of the notion of continuous semimartingale 
 is the one of   Dirichlet (possibly c\`adl\`ag) process (with respect to a given filtration $(\calf_t)$): it   was introduced by  \cite{fo} and \cite{ber} in the discretizations framework. 
A Dirichlet process $X$ is the sum of a local martingale $M$ and an adapted  process $A$ with zero quadratic variation: $A$ is the generalization of a bounded variation process.
If $A_0 = 0$, then the decomposition $M + A$ constitutes a the analogous
of the Doob-Meyer decomposition for a semimartingale $X$.
 However,  requiring $A$ to have 
zero quadratic variation imposes that $A$ is continuous, 
see Lemma \ref{L_bracket_jumps}; since a bounded variation process with jumps has a non zero finite quadratic variation,  the generalization of
 the semimartingale for jump processes is not 
completely represented by the notion of Dirichlet process. 
A natural generalization should then at least include the possibility that  
$A$   is a bounded variation process with jumps. 

The concept of $(\calf_t)$-weak Dirichlet process indeed  constitutes a
true generalization of the concept of (c\`adl\`ag) semimartingale.
Such a process is defined as the sum of a local  martingale $M$ and an adapted  process $A$ such that $[A,N]=0$ for every continuous 
local martingale $N$. 
For a continuous process $X$, that notion 
 was  introduced in \cite{er2}.
 In \cite{gr}
a chain rule  was established for $F(t,X_t)$ when  $F$ belongs to class $C^{0,1}$ and $X$ is a  weak Dirichlet process with finite quadratic variation. 
Such a process is indeed again a weak Dirichlet process (with possibly no finite quadratic variation). Relevant applications  to stochastic control were considered in \cite{gr1}. 
Contrarily to the continuous case, the decomposition $X = M +A$ 
is generally not unique. $X$ is denominated  \emph{special weak Dirichlet process} 
if it admits a decomposition of the same type
but where $A$ is predictable.
This concept is compatible with the one introduced 
in   \cite{cjms} 
using the discretization language.
The authors of \cite{cjms} were the first to introduce a notion of weak Dirichlet process
in the framework of jump processes. 
The  decomposition 
 of a special weak Dirichlet process is now unique, see Proposition
 \ref{P_unique_decomp}, at least fixing $A_0 = 0$.
 We remark that a continuous 
weak Dirichlet process is a special weak Dirichlet one. 
If the concept of (non necessarily special) weak Dirichlet process 
extends the notion of semimartingale,
the one of special weak Dirichlet process appears 
to be a generalization of the one of special semimartingale.

Towards calculus via regularization  in the jump case only a few steps were done in \cite{rv95}, \cite{rv93}, and several other authors, 
see Chapter 15 of \cite{nunno-book} and references therein. 
 For instance no  It\^o type formulae have been established 
in the framework of regularization and in the discretization framework  only very  few chain rule results 
are available for $F(X)$, when $F(X)$ is not a semimartingale.
 In that direction two peculiar results are available: 
the expansion of $F(X_t)$ when $X$ is a reversible semimartingale and 
$F$ is of class $C^1$ with some H\"older conditions on  the derivatives (see \cite{erv}),  and a chain rule for $F(X_t)$ 
when $X$ is a (c\`adl\`ag   special) weak Dirichlet process with finite quadratic variation and finite energy  and $F$ is of class $C^1_b$, see Corollary 3.2 in \cite{cjms}. 
The work in \cite{erv} has been continued by several authors, see 
e.g. \cite{eisen2007}   and references therein, expanding the remainder making use of local time type processes. 
A systematic study of  calculus via regularization was missing and this paper 
fills out this gap. 
  
Let us now go through the description of the main results of the paper. 
As we have already mentioned, our first basic objective consists in developing 
  calculus via regularization in the case of finite quadratic variation 
c\`adl\`ag processes. To this end, we revisit the definitions given by \cite{rv95} concerning 
forward integrals (resp. covariations). 
Those objects are introduced as u.c.p. (uniform convergence in probability) limit of 
 the  expressions of the type 
\eqref{Appr-ucp-int} (resp. \eqref{Appr-ucp-cov}). 
That convergence ensures that the limiting objects are c\`adl\`ag, since the approximating expressions have the same property. 
For instance a c\`adl\`ag process $X$ will be called  \emph{finite quadratic variation process} 
 whenever the limit (which will be denoted by $[X,X]$)     of 
\begin{equation} 
\label{Appr_cov_ucpI}	 
[X,X]^{ucp}_{\varepsilon}(t):= \,\int_{]0,\,t]}\,\frac{(X((s+\varepsilon)\wedge t)-X(s))^2}{\varepsilon}\,ds, 
\end{equation} 
exists u.c.p.   In \cite{rv95}, 
the authors  introduced a slightly different approximation 
of $[X,X]$ when $X$ is continuous, namely
\begin{equation} 
\label{Appr_cov_I} 
C_{\varepsilon}(X,X)(t):= \,\int_{]0,\,t]}\,\frac{(X(s+\varepsilon)-X(s))^2}{\varepsilon}\,ds. 
\end{equation} 
When the  u.c.p. limit of $C_{\varepsilon}(X,X)$ 
 exists, it is automatically a continuous process, since the approximating processes 
 are continuous. 
For this reason,  when $X$ is  a jump process,  the choice of approximation  \eqref{Appr_cov_I} 
would not be suitable, since its quadratic variation is expected to be a 
jump process. 
In that case,  however, the u.c.p. convergence of  \eqref{Appr_cov_ucpI} can  be shown to be equivalent to the a.s. pointwise convergence  (up to subsequences) of $ C_{\varepsilon}(X,X)$,
see Remark  \ref{R_equiv_mutual_brackets}. 
Both formulations will be used in the development of the calculus.
 
For a c\`adl\`ag finite quadratic variation process $X$, we establish, 
via  regularization techniques, 
 a quasi pathwise It\^o formula for $C^{1,2}$ functions of $X$. This is the object of 
 Theorem \ref{P_Ito_C2_cadlag}, whose proof is based on an 
accurate separation between the  neighborhood of "big" and "small" jumps, 
where  specific tools  are used, see for instance the preliminary  results Lemma 
 \ref{L_ucp_big_jumps} and  Lemma \ref{Lem_Billingsley}. 
Another significant instrument is  Lemma \ref{L_Fn_F_G}, which is
of  Dini type  for c\`adl\`ag functions. 
Finally, from Theorem \ref{P_Ito_C2_cadlag}  easily follows 
an It\^o formula under weaker regularity conditions on $F$, see Proposition \ref{P_RV_Itoformula_C1lambda}. 
We remark that a similar formula was stated in  \cite{erv}, 
using a discretization definition of the covariation, when $F$ 
is time-homogeneous. 

The second main task of the paper consists in  investigating 
 weak Dirichlet jump processes. 
Beyond some basic properties, two  significant results are Theorem \ref{T_C1_dec_weak_Dir} 
and Theorem  \ref{T_C1_dec_specialweak_Dir}. 
They both concern expansions of $F(t,X_t)$  where $F$ is of class $C^{0,1}$ and 
$X$ is a weak Dirichlet process of finite quadratic variation. 
Theorem \ref{T_C1_dec_weak_Dir} states that $F(t,X_t)$ 
will be again a weak Dirichlet process, however not necessarily of finite quadratic variation. 
 Theorem  \ref{T_C1_dec_specialweak_Dir} concerns the cases when $X$ and $F(t,X_t)$ are special weak Dirichlet processes. 
A first significant step in this sense was done in \cite{cjms}, where $X$ belongs to a 
 bit different  class of special weak 
Dirichlet jump processes (of finite energy) and $F$ does not depend on time 
and has bounded derivative. 
They show that $F(X)$ is again a special weak Dirichlet process. 
In  \cite{cjms} the underlying process has finite energy, which requires 
a control of the expectation of the approximating sequences of the quadratic variation. 
On the other hand, our techniques do not require that type of control. 
Moreover, the integrability condition \eqref{E_C} 
that we ask on $F(t,X_{t})$ in order to get the chain rule in Theorem \ref{T_C1_dec_specialweak_Dir} 
is automatically verified under the  hypothesis on the first order derivative considered in \cite{cjms}, see Remark  
\ref{R_bounded_derivative}.
 In particular, this allows us to recover the result in \cite{cjms} by means of our techniques, see Corollary \ref{C_Coquet}.
Finally, in some cases a chain rule may hold even when $F$ is not necessarily differentiable in the second variable, if we know a priori some information on the process $F(t,X_t)$. A result in this direction is provided by Proposition \ref{P_C00_chain_rule}, and does not require even that $X$ is a weak Dirichlet process.
 
In the present paper we also introduce  a  subclass of weak Dirichlet processes, called \emph{particular}, see 
Definition \ref{D_PWD}. 
 Those processes  inherit some of the 
semimartingales features: as in the semimartingale case, the particular weak Dirichlet processes  admit an integral representation 
(see Proposition \ref{P_char_part_weak_D})   and a 
 (unique) \emph{canonical decomposition} holds  when 
$ 
x\,\one_{\{|x| >1\}}\ast \mu \in \mathcal{A}_{\rm loc}. 
$ 
Under that condition, those particular processes are indeed   special weak Dirichlet processes, see  Proposition \ref{P_int_dec_par_WD} 
and \ref{P_char_part_weak_D}.

The paper is organized as follows. 
In Section 
\ref{Sec_tecnicalities_JumpCalculus} 
we introduce the  notations and 
we give some  preliminary results to the   development of 
the calculus via regularization with jumps. 
Section \ref{Section_Ito_C12} is devoted to the proof of the $C^{1,2}$ It\^o's formula for c\`adl\`ag processes with finite quadratic variation. 
In Section 
\ref{S32}  we recall some basic results on  the stochastic integration with respect to integer-valued random measures, and  we use them to reexpress the $C^{1,2}$ It\^o's formula of Section \ref{Section_Ito_C12} in terms of the (compensated) random measure associated to the c\`adl\`ag process.
Section \ref{Section_WDprocesses} concerns the study of weak Dirichlet processes, and presents the   expansion of $F(t,X_t)$ for $X$ weak Dirichlet, 
when $F$ is of class $C^{0,1}$. 
Finally,  we report in the Appendix \ref{Section_Appendix} some additional comments and technical results on calculus via regularizations.

\section{Preliminaries, calculus via regularization with jumps 
and related technicalities}\label{Sec_tecnicalities_JumpCalculus} 

Let $T>0$ be a finite horizon. We will indicate by $C^{1,2}$ (resp. $C^{0,1}$) 
the space of all functions 
$$ 
u: [0,\,T]\times \R \rightarrow \R, \quad (t,x)\mapsto u(t,x),
$$ 
that are continuous  together their derivatives $\partial_t u$, $\partial_x u$, $\partial_{xx} u$ (resp. $\partial_x u$). 
$C^{1,2}$ is equipped with the topology of uniform convergence on each compact of $u$, $\partial_x u$, $\partial_{xx} u$, $\partial_t u$, $C^{0,1}$ is equipped with the same topology on each compact of $u$ and $\partial_x u$. 
Given a topological space $E$, in the sequel $\mathcal{B}(E)$ will denote 
the Borel $\sigma$-field associated with $E$.

In the whole article, we are given a probability space $(\Omega,\mathcal{F},\P)$
and a filtration $\mathcal{F}=(\mathcal{F}_t)_{t \in [0,\,T]}$, fulfilling the usual conditions.
The symbols $\mathbb{D}^{ucp}$ and  $\mathbb{L}^{ucp}$ 
will denote  the  space of all adapted c\`adl\`ag and  c\`agl\`ad 
 processes endowed with the  u.c.p. (uniform convergence in probability) 
topology. By convention, any c\`adl\`ag process defined on $[0,\,T]$ is  extended on $\R_+$ by continuity.

Let $f$ and $g$ be two  functions defined on $\R$, and set 
\begin{eqnarray} 
	I^{-ucp}(\varepsilon, t,f, dg)&=&\int_{]0,\,t]} f(s)\,\frac{g((s+\varepsilon)\wedge t)-g(s)}{\varepsilon}\,ds, \label{Appr-ucp-int} \\ 
	\left[f,g\right]^{ucp}_{\varepsilon}(t)&=& \,\int_{]0,\,t]}\,\frac{(f((s+\varepsilon)\wedge t)-f(s)) 
(g((s+\varepsilon)\wedge t)-g(s))}{\varepsilon}\,ds. \label{Appr-ucp-cov} 
\end{eqnarray} 
Notice that the function $I^{-ucp}(\varepsilon, t,f, dg)$ is c\`adl\`ag  and 
  admits the   decomposition 
\begin{equation}\label{int_ucp_dec} 
I^{-ucp}(\varepsilon, t,f, dg) = \int_{0}^{(t-\varepsilon)_+}f(s)\,\frac{g(s + \varepsilon)-g(s)}{\varepsilon}\,ds 
+\int_{(t-\varepsilon)_+}^{t}f(s)\,\frac{g(t)- g(s)}{\varepsilon}\,ds. 
\end{equation} 
\begin{definition}\label{D_ucp_integral} 
Let $X$ be a c\`adl\`ag process and $Y$ be a  process belonging to $L^1([0,\,T])$ a.s. 
Suppose the existence of a process  $(I(t))_{t \in [0,\,T]}$ such that 
$(I^{-ucp}(\varepsilon, t,Y, dX))_{t \in [0,\,T]}$ converges u.c.p. to  $(I(t))_{t \in [0,\,T]}$, namely 
	\begin{displaymath} 
	\lim_{\varepsilon \rightarrow 0}\,\P\left(\sup_{0 \leq s \leq t}|I^{-ucp}(\varepsilon, t,Y, dX) - I(t)| > \alpha \right)= 0\qquad \textup{for every}\,\,\,\alpha >0. 
	\end{displaymath} 
Then we will set 
 $\int_{]0,\,t]}\,Y_s \, d^{-} X_s := I(t)$. 
That process will be called 
 \textbf{the forward integral of $Y$ with respect to $X$}.  
\end{definition} 
\begin{remark}\label{R_RV95} 
	In \cite{rv95}   a very similar notion of forward integral 
	 is considered: 
	\[ 
	I^{-RV}(\varepsilon, t,f, dg) =\int_{\R} f_{t]}(s)\,\frac{g_{t]}(s+\varepsilon)-g_{t]}(s)}{\varepsilon}\,ds, 
	\] 
	with 
	\begin{eqnarray*} 
		f_{t]}= 
		\left\{ 
		\begin{array}{ll} 
			f(0_+) &\textup{if}\,\, x \leq 0,\\ 
			f(x) &\textup{if}\,\, 0 < x \leq t,\\ 
			f(t_+) &\textup{if}\,\, x > t. 
		\end{array} 
		\right. 
	\end{eqnarray*} 
	The u.c.p. limit of  $I^{-RV}(\varepsilon, t,f, dg)$, when it exists, coincides with  that of $I^{-ucp}(\varepsilon, t,f, dg)$. 
	As a matter of fact, the process $I^{-RV}(\varepsilon, t,f, dg)$ is c\`adl\`ag  and can be rewritten as 
	\begin{eqnarray}\label{int_RV_dec} 
		I^{-RV}(\varepsilon, t,f, dg) &=& I^{-ucp}(\varepsilon, t,f, dg) -f(0_+)\,\frac{1}{\varepsilon}\int_0^\varepsilon [g(s)-g(0_+)]\,ds. 
	\end{eqnarray} 
	In particular 
	\[ 
	\sup_{t \in [0,\,T]} 
	[ 
	I^{-ucp}(\varepsilon, t,f, dg) - I^{-RV}(\varepsilon, t,f, dg) 
	] 
	= 
	f(0_+)\,\frac{1}{\varepsilon}\int_0^\varepsilon [g(s)-g(0_+)]\,ds, 
	\] 
	and therefore 
	\[ 
	\limsup_{\varepsilon \rightarrow 0} \sup_{t \in [0,\,T]} 
	[ 
	I^{-RV}(\varepsilon, t,f, dg) - I^{-ucp}(\varepsilon, t,f, dg) 
	] = 0. 
	\] 
\end{remark} 
\begin{proposition}\label{ucp_predic} 
	Let  $A$ be a c\`adl\`ag predictable process and $Y$ be a  process belonging to $L^1([0,\,T])$ a.s. Then the 
forward integral 
	\begin{displaymath} 
		\int_{]0,\,\cdot]}Y_s\,d^{-}A_s, 
	\end{displaymath} 
when it exists,	is a predictable process. 
\end{proposition} 
\proof 
From decomposition \eqref{int_ucp_dec} it follows that  
$
I^{-ucp}(\varepsilon, \cdot,Y, dA) = C^\varepsilon + D^\varepsilon,
$
where $C^\varepsilon$ is a continuous adapted process, therefore predictable, and  
$$
D^\varepsilon_t = A_t \,\frac{1}{\varepsilon} \int_{(t-\varepsilon)+}^t Y_s\,ds, \quad t \in [0,\,T].
$$
$D^\varepsilon$ is the product of the predictable process $A$ and another continuous adapted process, so it is itself predictable.
By definition, the u.c.p. stochastic integral, when it exists,  is the u.c.p. limit  of $I^{-ucp}(\varepsilon, \cdot,Y, dA)$: 
since the u.c.p. convergence 
preserves the predictability, the claim follows. 
\endproof 
\begin{definition}\label{D_ucp_cov} 
Let  $X, Y$ be two c\`adl\`ag processes. 
Suppose the existence of a process  $(\Gamma(t))_{t\geq 0}$ such that 
$\left[X,Y\right]^{ucp}_{\varepsilon}(t)$ converges u.c.p. to $(\Gamma(t))_{t\geq 0}$, 
namely 
	\begin{displaymath} 
	\lim_{\varepsilon \rightarrow 0}\,\P\left(\sup_{0 \leq s \leq t}|[X,Y]^{ucp}_{\varepsilon}(t) - \Gamma(t)| > \alpha \right)= 0\qquad \textup{for every}\,\,\,\alpha >0, 
	\end{displaymath} 
Then we will set 
$[X,Y]_t := \Gamma(t)$. 
That process will be called 
\textbf{the covariation between $X$ and $Y$}.	 
In that case we say that 	 
\textbf{the covariation between $X$ and $Y$  exists}. 
\end{definition} 
\begin{definition}\label{D_quadr_var} 
	We say that a c\`adl\`ag process $X$ \textbf{is a finite quadratic variation}  if $[X,X]$  exists. In that case $[X,X]$ is called the quadratic variation of $X$.
\end{definition} 
\begin{definition}\label{D_mutual_brackets} 
	We say that a pair of c\`adl\`ag processes $(X,Y)$ \textbf{admits all its mutual brackets} if $[X,X]$, $[X,Y]$, $[Y,Y]$  exist. 
\end{definition} 
\begin{remark}\label{Rem_quadr_var_increasing} 
	Let $X$, $Y$ be two c\`adl\`ag processes.  
\begin{enumerate} 
\item By definition	$[X,Y]$ is necessarily a c\`adl\`ag process. 
\item $[X,X]$ is an increasing process.  
\end{enumerate} 
\end{remark} 
$[X,X]^c$ denotes the continuous part of  $[X,X]$.

Forward integrals and covariations generalize  It\^o integrals and the classical square brackets of semimartingales. 
\begin{proposition}\label{P_Ito_cov} 
Let  $X, Y$ be two c\`adl\`ag semimartingales, $M^1, M^2$ two c\`adl\`ag  local martingales, 
$H, K$ two c\`adl\`ag adapted  process. Then the following properties hold.
\begin{itemize} 
\item[(i)] $[X,Y]$ exists and it is the usual bracket. 
\item[(ii)] $\int_{]0,\,\cdot]} H\, d^{-}X$  is the usual stochastic integral 
 $\int_{]0,\,\cdot]} H_{s-} dX_{s}$. 
\item[(iii)]$\left[\int_0^{\cdot} H_{s-}\,dM^1_{s},\,\int_0^{\cdot} 
K_{s-}\,dM^2_{s}\right]$  is the usual bracket and it equals 
$\int_0^{\cdot} H_{s-}\, K_{s-}\,d[M^1,M^2]_s$. 
\end{itemize} 
\end{proposition} 
\proof 
Items (i) and (ii) are consequence of 
 Proposition 1.1 in \cite{rv95} and  Remark \ref{R_RV95}. 
Item (iii) follows from (i) and the corresponding properties 
for classical brackets of local martingales, see Theorem 29, chapter 2 
of \cite{protter}. 
\endproof 
\begin{remark}\label{R_StandAss} 
In the sequel we will make often use of the following assumption on a  c\`adl\`ag  process $X$: 
\begin{equation}\label{StandAss} 
\sum_{s \in ]0,\, T]} |\Delta X_s|^2 < \infty, 	\,\, \textup{a.s.} 
\end{equation} 
	Condition \eqref{StandAss} holds for instance in the case of processes $X$ of finite quadratic variation. 
\end{remark} 
\begin{lemma}\label{L_bracket_jumps} 
	Suppose that $X$ is a  c\`adl\`ag, finite quadratic variation  process. Then 
	\begin{itemize} 
	\item[(i)] $\forall s \in [0,\,T]$,\,\, 
	$ 
		\Delta [X,X]_s = (\Delta X_s)^2; 
	$ 
	\item[(ii)] 
		$ 
		[X,X]_s= [X,X]_s^c + \sum_{t \leq s}(\Delta X_t)^2,\quad \forall s \in [0,\,T],\,\,\textup{a.s.} 
		$ 
		 
		In particular, \eqref{StandAss} holds. 
	\end{itemize} 
\end{lemma} 
\proof 
(i) Since  $X$ is a  finite quadratic variation process, $[X,X]^{ucp}_\varepsilon$ converges u.c.p. to $[X,X]$. 
This implies the existence of a sequence $(\varepsilon_n)$ such that 
 $[X,X]^{ucp}_{\varepsilon_n}$ converges uniformly a.s. 
to $[X,X]$. 
We fix a realization $\omega$ outside a suitable null set, which will be omitted in the sequel. Let $\gamma >0$. 
 There is $\varepsilon_0$ such that 
\begin{equation} \label{EF1} 
	\varepsilon_n < \varepsilon_0 \Rightarrow |[X,X]_s - 
[X,X]^{ucp}_{\varepsilon_n}(s)| \leq \gamma, \quad \forall s \in [0,\,T]. 
\end{equation} 
We fix $s \in ]0,\,T]$. 
Let $\varepsilon_n < \varepsilon_0$. 
For every $\delta \in [0,\,s[$, we have  
\begin{equation} \label{EF2} 
	|[X,X]_{s} - [X,X]^{ucp}_{\varepsilon_n}(s- \delta)| \leq \gamma. 
\end{equation} 
Since $[X,X]$ is  c\`adl\`ag, we need to show that the quantity 
\begin{equation}\label{def_quantity} 
|[X,X]_{s} - [X,X]_{s- \delta}- (\Delta X_s)^2| 
\end{equation} 
goes to zero, when $\delta \rightarrow 0$. 
 For $\varepsilon:= \varepsilon_n < \varepsilon_0$, \eqref{def_quantity}  is smaller or equal 
than 
\begin{eqnarray*} 
	&&2\gamma + |[X,X]^{ucp}_{\varepsilon}(s)- [X,X]^{ucp}_{\varepsilon}(s- \delta)- (\Delta X_s)^2|\\ 
	&&= 2 \gamma 
	+ \bigg|\frac{1}{\varepsilon}\int_{s-\varepsilon - \delta}^{s}(X_{(t+ \varepsilon )\wedge s}- X_t)^2 \,dt 
	-\frac{1}{\varepsilon}\int_{s-\varepsilon- \delta}^{s- \delta}(X_{ s- \delta}- X_t)^2 \,dt 
	-(\Delta X_s)^2 
	\bigg|\\ 
	&&\leq  2 \gamma + \frac{1}{\varepsilon}\int_{s-\varepsilon- \delta}^{s- \delta}(X_{ s- \delta}- X_t)^2 \,dt + |I(\varepsilon,\,\delta,\,s)|,\quad \forall \delta \in [0,\,s[, 
\end{eqnarray*} 
where 
\begin{displaymath} 
	I(\varepsilon,\,\delta,\,s) = 
	\frac{1}{\varepsilon}\int_{s-\varepsilon - \delta}^{s- \varepsilon}(X_{t+ \varepsilon}- X_t)^2 \,dt 
	+ \frac{1}{\varepsilon}\int_{s-\varepsilon}^{s}[(X_{s}- X_t)^2- (\Delta X_s)^2] \,dt. 
\end{displaymath}

At this point, 
we have 
\begin{displaymath} 
	|[X,X]_{s} - [X,X]_{s- \delta}- (\Delta X_s)^2|\leq 2 \gamma + \frac{1}{\varepsilon}\int_{s-\varepsilon- \delta}^{s- \delta}(X_{ s- \delta}- X_t)^2 \,dt + |I(\varepsilon,\,\delta,\,s)|, \quad \forall s \in [0,\,T]. 
\end{displaymath} 
We take the $\limsup_{\delta \rightarrow 0}$ on both sides to get, since $X$ has a left limit at $s$, 
\begin{displaymath} 
	|\Delta [X,X]_{s} - (\Delta X_s)^2|\leq 2 \gamma + \frac{1}{\varepsilon}\int_{s-\varepsilon}^{s}(X_{s-}- X_t)^2 \,dt + \frac{1}{\varepsilon}\int_{s-\varepsilon}^{s}|(X_{s}- X_t)^2- (\Delta X_s)^2| \,dt, \quad \textup{for } 
\varepsilon:= \varepsilon_n < \varepsilon_0. 
\end{displaymath} 
We take the limit when $n \rightarrow \infty$ 
and we  get 
$ |\Delta [X,X]_{s} - (\Delta X_s)^2|\leq 2 \gamma$,
and this concludes the proof of (i). 
 
(ii) We still work fixing a priori a realization $\omega$. 
 Set $Y_s = [X,X]_s$,\,\, $s \in [0,\,T]$. Since $Y$ is an increasing c\`adl\`ag process, it can be decomposed as 
$Y_s= Y_s^c + \sum_{t \leq s}\Delta Y_t$, for all $s \in [0,\,T]$, a.s.,
and  the result follows 
from point (i). 
In particular, setting $s=T$, we get 
$$ 
\sum_{s \leq T}(\Delta X_s)^2 \leq \sum_{s \leq T}(\Delta X_s)^2  + [X,X]_T^c = [X,X]_T < \infty,\,\,\textup{a.s.} 
$$
\endproof

We now state and prove some fundamental preliminary results, that we will deeply use in the sequel. 
\begin{lemma}\label{L_ucp_big_jumps} 
	Let  $Y_t$  be a c\`adl\`ag function with values in $\R^n$. 
	Let $\phi: \R^n \times \R^n \rightarrow \R$ be an equicontinuous function on each compact, such that $\phi(y,y)=0$ for every $y \in \R^n$. 
	Let $0 \leq t_1 \leq t_2 \leq ... \leq t_N \leq T$. 
	We have 
	\begin{align} 
		\sum_{i=1}^N \frac{1}{\varepsilon} \int_{t_i - \varepsilon}^{t_i} \one_{]0,\,s]}(t)\,\phi(Y_{(t+ \varepsilon)\wedge s},Y_t)\,dt \overset{\varepsilon \rightarrow 0}{\longrightarrow}\,\, \sum_{i=1}^N  \one_{]0,\,s]}(t_i)\, \phi(Y_{t_i},Y_{t_{i}-}), \label{u_conv_prod_big_j} 
	\end{align} 
	uniformly in $s \in [0,\,T]$. 
\end{lemma} 
\proof 
Without restriction of generality, we consider the case  $n=1$. 
Let us fix $\gamma >0$. Taking into account that $\phi$ is equicontinuous on compacts, by definition of left and right limits, there exists $\delta >0$ such that, for every $i \in \{1,...,N\}$, 
\begin{align} 
	\ell < t_{i},\, u > t_{i},\, |\ell-{t_{i}}| \leq \delta, \,|u-{t_{i}}| \leq \delta &\Rightarrow |\phi(Y_u, Y_\ell)-\phi(Y_{t_i}, Y_{{t_i}-})| < \gamma,\label{DS}\\ 
	\ell_2 < \ell_1 < t_{i},\, |\ell_1-{t_{i}}| \leq \delta, \,|\ell_2-{t_{i}}| \leq \delta &\Rightarrow |\phi(Y_{\ell_1}, Y_{\ell_2})| = |\phi(Y_{\ell_1}, Y_{\ell_2})-\phi(Y_{t_{i}-}, Y_{{t_i}-})| < \gamma.\label{SS} 
\end{align} 
Since the  sum in \eqref{u_conv_prod_big_j} is finite, it is enough to show 
the uniform convergence in $s$ of  the integrals on  $]t_{i}-\varepsilon,\,t_i]$, for a fixed $t_i \in [0,\,T]$, namely that 
\begin{align} 
	I(\varepsilon, s):=\frac{1}{\varepsilon} \int_{t_i - \varepsilon}^{t_i} \one_{]0,\,s]}(t)\,\phi(Y_{(t+ \varepsilon)\wedge s},Y_t)\,dt -\one_{]0,\,s]}(t_i)\, \phi(Y_{t_i},Y_{t_{i}-})\label{E:TermToEvaluate1} 
\end{align} 
converges to zero uniformly  in $s$, when $\varepsilon$ goes to zero. 
Let thus fix $t_i \in [0,\,T]$, and choose $\varepsilon < \delta$. 
We distinguish the  cases (i), (ii), (iii), (iv) concerning 
the position of $s$ with respect to $t_i$. 
\begin{itemize} 
	\item[(i)] $s < t_i - \varepsilon$. \eqref{E:TermToEvaluate1} 
	vanishes. 
	\item [(ii)]  $s \in [t_i-\varepsilon,\,t_i[$. 
	By \eqref{SS} we get 
	$$ 
	|I(\varepsilon,s)| \leq \frac{1}{\varepsilon} \int_{t_i - \varepsilon}^{t_i}|\phi(Y_{s},Y_t)|\,dt \leq \gamma. 
	$$ 
	\item[(iii)]  $s \in [t_i,\,t_i+\varepsilon[$. 	By \eqref{DS} we get 
	$$ 
	|I(\varepsilon,s)| \leq \frac{1}{\varepsilon} \int_{t_i - \varepsilon}^{t_i}|\phi(Y_{(t+ \varepsilon)\wedge s },Y_t)-\phi(Y_{t_i},Y_{t_{i}-})|\,dt \leq \gamma. 
	$$ 
	\item[(iv)]  $s \geq t_i+\varepsilon$. 
	By \eqref{DS} we get 
	$$ 
	|I(\varepsilon,s)| \leq \frac{1}{\varepsilon} \int_{t_i - \varepsilon}^{t_i}|\phi(Y_{t+ \varepsilon },Y_t)-\phi(Y_{t_i},Y_{t_{i}-})|\,dt \leq \gamma. 
	$$ 
\end{itemize} 
Collecting all the cases above, we see that 
$ \limsup_{\varepsilon \rightarrow 0}\sup_{s \in [0,\,T]}|I(\varepsilon,s)| 
	 \leq \gamma$,
and letting   $\gamma$ go to zero we get 
 the uniform convergence. 
\endproof 
 
\begin{lemma}\label{Lem_Billingsley} 
Let $X$ be a c\`adl\`ag (c\`agl\`ad) real process. 
	Let $\gamma>0$, $t_0,\,t_1 \in \R$ and  $I = [t_0,\,t_1]$ be a subinterval of $[0,\,T]$ such that 
	\begin{equation} \label{EJump} 
		|\Delta X_t|^2 \leq \gamma^2,\quad \forall t \in I. 
	\end{equation} 
	Then there is $\varepsilon_0 >0$ such that 
	\begin{displaymath} 
		\sup_{\underset{|a-t| \leq \varepsilon_0}{a,\,t \in I}}|X_a - X_t|\leq 3 \gamma. 
	\end{displaymath} 
\end{lemma} 
\proof 
We only treat the c\`adl\`ag case, the c\`agl\`ad one is a consequence of an obvious time reversal argument. 
 
Also in this proof a realization $\omega$ will be fixed, but omitted. 
According to Lemma 1, Chapter 3,  in \cite{billingsley}, applied to $[t_0,\,t_1]$ replacing $[0,\,1]$, there exist points 
\begin{displaymath} 
	t_0 = s_0 < s_1 <...< s_{l-1} < s_l = t_1 
\end{displaymath} 
such that for every $j \in \{1,..., l\}$ 
\begin{equation}\label{sup} 
	\sup_{d,\,u\, \in [s_{j-1},\,s_{j}[}|X_d -X_u| < \gamma. 
\end{equation} 
Since $X$ is c\`adl\`ag, we can  choose $\varepsilon_0$ such that, $\forall j \in \{0,\,...,\,l-1\}$, 
\begin{eqnarray} 
	|d-s_j| \leq \varepsilon_0 &\Rightarrow& |X_d - X_{s_{j}-}| \leq \gamma,\label{remind_1}\\ 
	|u-s_j| \leq \varepsilon_0 &\Rightarrow& |X_u - X_{s_{j}}| \leq \gamma.\label{remind_2} 
\end{eqnarray} 
Let $t \in [s_{j-1},\,s_j[$ for some $j$ and $a$ such that $|t-a| \leq \varepsilon$ for $\varepsilon < \varepsilon_0$. Without restriction of generality we can take 
$ t < a$. 
There are two cases. 
\begin{itemize} 
	\item[(i)] $a,\,t \in [s_{j-1},\,s_j[$. 
		In this case, \eqref{sup} gives 
	$	|X_{a}- X_t| < \gamma$.
	\item[(ii)] $s_{j-1} \leq t <  s_{j} \leq a$. 
	 Then, 
	\begin{eqnarray*} 
		|X_a- X_t| \leq |X_{a}- X_{s_j}| 
		+ |X_{s_j}- X_{s_j-}|+ |X_{s_j-}- X_t|\leq  3\gamma, 
	\end{eqnarray*} 
	where the first absolute value is bounded by \eqref{remind_2}, the second by 
\eqref{EJump} 
and the third  by 
 \eqref{remind_1}. 
\end{itemize} 
\endproof 
 
\begin{remark}\label{R_Billingsley} 
Let $I = [t_{0},\,t_1] \subset [0,\,T]$, 
 let $\varepsilon > 0$. 
Let $t \in ]t_{0},\,t_1 - \varepsilon]$ and $s>t$.  
We will apply Lemma \ref{Lem_Billingsley} to the couple $(a,t)$, where $a= (t + \varepsilon)\wedge s$. 
Indeed $a \in I$ because $ a \le t + \varepsilon \le t_1.$ 
\end{remark} 
\begin{proposition}\label{P_quad_var} 
	Let  $(Z_t)$ be a c\`adl\`ag process, $(V_t)$ be a bounded variation process. Then $[Z,V]_s$ 
 exists and equals 
	\[ 
	\sum_{t \leq s} \Delta Z_t\,\Delta V_t,\quad \forall s \in [0,\,T]. 
	\] 
	In particular,  $V$ is a finite quadratic variation process. 
\end{proposition} 
\proof 
We need to prove the u.c.p.  convergence  to zero of 
\begin{equation}\label{approx_cov} 
	\frac{1}{\varepsilon}\int_{]0,\,s]}(Z_{(t + \varepsilon)\wedge s}-Z_{t})(V_{(t + \varepsilon)\wedge s}-V_{t})\,dt -\sum_{t \leq s} \Delta Z_t\,\Delta V_t. 
\end{equation} 
As usual the realization   $\omega \in \Omega$ will be fixed, but often omitted. 
 Let $(t_i)$ be the enumeration of 
  all 
 the jumps of $Z(\omega)$ in $[0,\,T]$. 
 We have 
\[ 
\lim_{i \rightarrow \infty} |\Delta Z_{t_i}(\omega)| = 0. 
\] 
Indeed, if it were not the case, it would exists $a>0$ and a 
subsequence $(t_{i_l})$ of $(t_{i})$ such that $|\Delta Z_{t_{i_l}}| \geq a$. This is not 
possible   since a c\`adl\`ag function admits at most  
a finite number of jumps exceeding any $a > 0$, 
see considerations below Lemma 1, Chapter 2 of \cite{billingsley}. 
 
At this point, let $\gamma >0$ and $N= N(\gamma)$ such that 
\begin{equation} \label{E_N} 
n \geq N, \quad |\Delta Z_{t_n}| \leq \gamma. 
\end{equation} 
We introduce 
\begin{align}
A(\varepsilon,N) &=\bigcup_{i=1}^{N} \, ]t_{i}-\varepsilon,t_{i}], \label{Aepsilon}\\
	B(\varepsilon,N) &= \bigcup_{i=1}^{N} \, ]t_{i-1},t_i - \varepsilon]=[0,\,T] \setminus A(\varepsilon,N), \label{Bepsilon}
\end{align}
and we decompose \eqref{approx_cov} into 
\begin{equation}\label{approx_cov_1} 
	I_A(\varepsilon,N,s)+ I_{B1}(\varepsilon,N,s)+I_{B2}(\varepsilon,N,s) 
\end{equation} 
where 
\begin{align*} 
	I_A(\varepsilon,N,s)&=\frac{1}{\varepsilon}\int_{]0,\,s] \cap A(\varepsilon,N) }(Z_{(t + \varepsilon)\wedge s}-Z_{t})(V_{(t + \varepsilon)\wedge s}-V_{t})\,dt - \sum_{i=1}^{N}\one_{]0,\,s[}(t_i)\,\Delta Z_{t_i}\,\Delta V_{t_i},\\ 
	I_{B1}(\varepsilon,N,s)&=\frac{1}{\varepsilon}\int_{]0,\,s]\cap B(\varepsilon,N) }(Z_{(t + \varepsilon)\wedge s}-Z_{t})(V_{(t + \varepsilon)\wedge s}-V_{t})\,dt,\\ 
	I_{B2}(N,s)&=- \sum_{i=N+1}^{\infty}\one_{]0,\,s[}(t_i)\,\Delta Z_{t_i}\,\Delta V_{t_i}. 
\end{align*} 
Applying Lemma \ref{L_ucp_big_jumps} to  $Y=(Y^1,Y^2)=(Z,V)$ and $\phi(y_1,y_2)=(y_1^1-y^1_2)(y_1^2-y^2_2)$ 
 we get 
\[ 
I_{A}(\varepsilon,N,s) \underset{\varepsilon \rightarrow 0}{\rightarrow} 0, \,\, \textup{uniformly in } s.
\]  
On the other hand,  for  $t \in ]t_{i-1},\,t_i - \varepsilon[$ and $s > t$, 
by Remark \ref{R_Billingsley} we know  that 
$(t + \varepsilon)\wedge s \in [t_{i-1},\,t_i]$. Therefore Lemma \ref{Lem_Billingsley} 
with $X = Z$, 
 applied successively  to the intervals $I = [t_{i-1},\,t_i]$ implies that 
\begin{align*} 
	\vert I_{B1}(\varepsilon,N,s)\vert & 
=\frac{1}{\varepsilon}\int_{]0,\,s]\cap B(\varepsilon, N) }\vert Z_{(t + \varepsilon)\wedge s}-Z_{t}\vert 
\vert V_{(t + \varepsilon)\wedge s}-V_{t}\vert\,dt\\ 
	&\leq 3 \,\gamma \, \frac{1}{\varepsilon}\int_{]0,\,s]\cap B(\varepsilon,N) } 
\vert V_{(t + \varepsilon)\wedge s}-V_{t}\vert \,dt\\ 
	&\leq 3 \,\gamma \, \int_{]0,\,s]}|V_{(t + \varepsilon)\wedge s}-V_{t}|\,\frac{dt}{\varepsilon}\\ 
	&=3 \,\gamma \, \int_{]0,\,s]}\frac{dt}{\varepsilon}\int_{]t,\,(t+ \varepsilon)\wedge s]} 
d\Vert V \Vert_r\\ 
	&=3 \,\gamma \, \int_{]0,\,s]}d\Vert V \Vert_r\int_{[(r-\varepsilon)^+,\,r[} 
\frac{dt}{\varepsilon}\\ 
	&\leq 3 \,\gamma\,||V||_T, 
\end{align*} 
where $r \mapsto \Vert V \Vert_r$ denotes the total variation function of $V$. 
Finally, concerning $I_{B2}(N,s)$, by \eqref{E_N} we have 
\begin{align*} 
	\vert I_{B2}(N,s) \vert 
\leq \gamma\,\sum_{i=N+1}^{\infty}\one_{]0,\,s[}(t_i)\,\,|\Delta V_{{t_i}}| \leq \gamma\, ||V||_T. 
\end{align*} 
Therefore, collecting the previous estimations we get 
\[ 
\limsup_{\varepsilon \rightarrow 0} \sup_{s \in [0,\,T]} 
\vert I_{A}(\varepsilon,N,s) + I_{B1}(\varepsilon,N,s) + I_{B2}(N,s)\vert 
 \leq 4\,\gamma\, ||V||_T, 
\] 
and we conclude by the arbitrariness of $\gamma > 0$. 
\endproof 
\begin{proposition}\label{P_2.15}
Let $X$ be a c\`adl\`ag  process with finite quadratic variation and $A$ be a process such that $[A,A]=0$. Then $[X,A]$ exists and equals zero. 
\end{proposition}
\proof
By Remark \ref{R_equiv_mutual_brackets}, $[X,X]$  and $[A,A]$ exist  in the pathwise sense.
By Lemma \ref{L_bracket_jumps}, $A$ is continuous. 
Consequently, by Proposition 1-6) of \cite{RVSem}, $[X,A]=0$ in the pathwise sense. Finally, by Remark \ref{R_equiv_mutual_brackets},  $(X, A)$ admits all its mutual brackets and $[X,A]=0$.
 \endproof
Finally we  give a generalization of Dini type lemma in the c\`adl\`ag case. 
\begin{lemma}\label{L_Fn_F_G} 
	Let $(G_n, \, n \in \N)$ be a sequence of continuous increasing functions, let $G$ (resp. $F$) from $[0,\,T]$ to $\R$ be a c\`adl\`ag (resp. continuous) function. We set $F_n = G_n + G$ and suppose that $F_n \rightarrow F$ pointwise. Then 
	\begin{displaymath} 
		\limsup_{n \rightarrow \infty}\sup_{s \in [0,\,T]}|F_n(s)-F(s)| \leq 2 \sup_{s \in [0,\,T]}|G(s)|. 
	\end{displaymath} 
\end{lemma} 
\proof 
Let us fix $m\in \N^\ast$. Let $0=t_0 < t_1<...<t_m = T$ such that $t_i = \frac{i}{m}$, $i=0,...,m$. 
If $s \in [t_i,\,t_{i+1}]$, $0 \leq i \leq m-1$, we have 
\begin{equation}\label{EA1} 
	F_n(s)-F(s) \leq F_n(t_{i+1})-F(s)+G(s)-G(t_{i+1}). 
\end{equation} 
Now 
\begin{align} \label{EA1bis} 
	F_n(t_{i+1})-F(s) &\leq F_n(t_{i+1})-F(t_{i+1})+ F(t_{i+1})-F(s)\nonumber\\ 
	&\leq \delta\left(F,\,\frac{1}{m}\right) +F_n(t_{i+1})-F(t_{i+1}),
\end{align} 
where $\delta(F,\cdot)$ denotes the modulus of continuity of $F$.
From \eqref{EA1} and \eqref{EA1bis} it follows 
\begin{align} 
	F_n(s)-F(s) &\leq F_n(t_{i+1})-F(t_{i+1})+ G(s)-G(t_{i+1})+ 
\delta\left(F,\,\frac{1}{m}\right) \nonumber\\ 
	&\leq 2 ||G||_{\infty}+\delta\left(F,\,\frac{1}{m}\right) +|F_n(t_{i+1})-F(t_{i+1})|, \label{E_first_ineq} 
\end{align} 
where $||G||_{\infty}= \sup_{s \in [0,\,T]}|G(s)|$. 
Similarly, 
\begin{equation}\label{E_second_ineq} 
	F(s)-F_n(s) \geq - 2||G||_{\infty} - \delta\left(F,\,\frac{1}{m}\right) - 
 |F_n(t_{i})-F(t_i)|. 
\end{equation} 
So, collecting \eqref{E_first_ineq} and \eqref{E_second_ineq} we have,  $\forall s\in [t_i,\,t_{i+1}]$,
\begin{equation*} 
	|F_n(s)-F(s)| \leq 2 ||G||_{\infty}+ \delta\left(F,\,\frac{1}{m}\right) + 
|F_n(t_{i})-F(t_i)|+|F_n(t_{i+1})-F(t_{i+1})|. 
\end{equation*} 
Consequently, 
\begin{eqnarray}\label{EA4} 
	\sup_{s\in [0,\,T]}|F_n(s)-F(s)| \leq 2 ||G||_{\infty}+ \delta\left(F,\,\frac{1}{m}\right) 
 + \sum_{i=1}^{m}|F_n(t_{i})-F(t_i)|. 
\end{eqnarray} 
Recalling that $F_n \rightarrow F$ pointwise, taking the $\limsup$ in \eqref{EA4} we get 
\begin{displaymath} 
	\limsup_{n \rightarrow \infty}\sup_{s \in [0,\,T]}|F_n(s)-F(s)| \leq 2 ||G||_{\infty}+ 
\delta\left(F,\,\frac{1}{m}\right). 
\end{displaymath} 
Since $F$ is uniformly continuous, we let $m$ go to infinity,  
so that the result follows. 
\endproof

\section{It\^{o}'s formula for $C^{1,2}$ functions: the basic formulae}\label{Section_Ito_C12} 
 
 
We start with the It\^o formula for finite quadratic  variation processes 
in the sense of calculus via regularizations. 
 
\begin{theorem}\label{P_Ito_C2_cadlag} 
Let $X$ be a finite quadratic variation c\`adl\`ag process and $F: [0,\,T]\times \R \rightarrow \R $ a function of class $C^{1,2}$. Then we have 
\begin{align}\label{Ito_formula_C2} 
F(t,X_t)&=F(0,X_0)+\int_0^t\partial_s F(s,X_s)\,ds +\int_0^t\partial_x F(s,X_{s})\,d^{-}X_s+\frac{1}{2}\int_0^t\partial_{xx}^2 F(s,X_{s-})\,d[X,X]_s^c \nonumber\\ 
&+ \sum_{s\leq t}\,[F(s,X_s)-F(s,X_{s-})- \partial_x F(s,X_{s-})\,\Delta \,X_{s}]. 
\end{align} 
\end{theorem}
\proof 
Since $X$ is a finite quadratic variation process, by Lemma \ref{L_ucp_conv}, 
taking into account Definition \ref{D_cov} and Corollary 
\ref{C_id_RV_fw}-2), for a given c\`adl\`ag process $(g_t)$ 
we have 
$$ \int_{0}^{s}g_t\,(X_{(t + \varepsilon)\wedge s}-X_t)^2\,\frac{dt}{\varepsilon} \overset{\varepsilon \rightarrow 0}{\longrightarrow}\int_{0}^{s}g_{t-}\,d[X,X]_t\quad \textup{u.c.p.} 
$$ 
Setting $g_t = 1$ and $g_t = \frac{\partial_{xx}^2 F(t,\,X_{t})}{2}$, 
there exists  a positive sequence $\varepsilon_n$ such that 
\begin{align} 
\lim_{n \rightarrow \infty}\int_0^s(X_{(t + \varepsilon_n)\wedge s}-X_t)^2\,\frac{dt}{\varepsilon_n} &= [X,X]_s, 
\label{A2} \\ 
\lim_{n \rightarrow \infty}\int_0^s \frac{\partial_{xx}^2 F(t,\,X_{t})}{2}\,(X_{(t + \varepsilon_n)\wedge s}-X_t)^2\,\frac{dt}{\varepsilon_n} &= \int_{]0,\,s]} \frac{\partial_{xx}^2 F(t,\,X_{t-})}{2}\,d[X,X]_t, 
\label{A1} 
\end{align} 
uniformly in $s$, a.s. 
Let then $\mathcal N$ be a null set such that \eqref{A2}, \eqref{A1} hold for every $\omega \notin \mathcal N$. 
 
In the sequel we fix $\gamma >0$, $\varepsilon >0$, and $\omega \notin \mathcal N$, 
and we enumerate the jumps of $X(\omega)$  on $[0,\,T]$ by $(t_i)_{i \geq 0}$. 
Let $N = N(\omega)$ such that 
\begin{equation}\label{gamma_N+1_condition} 
\sum_{i = N+1}^{\infty} |\Delta X_{t_i}(\omega)|^2 \leq \gamma^2. 
\end{equation} 
From now on the dependence on $\omega$ will be often neglected. 
The quantity 
\begin{equation}\label{J0} 
	J_0(\varepsilon,\,s) =	\frac{1}{\varepsilon} \int_0^s [F((t + \varepsilon)\wedge s,\,X_{(t + \varepsilon)\wedge s})- F(t,\,X_t)]\,dt,\quad s \in[0,\,T], 
\end{equation} 
converges  to $F(s,\,X_s)-F(0,\,X_0)$ uniformly in $s$. 
As a matter of fact, setting $Y_t = (t,\,X_t)$, we have 
\begin{align}\label{J0_conv} 
J_0(\varepsilon,\,s) 
&=	\frac{1}{\varepsilon}  \int_{[0,\,s[} F(Y_{(t + \varepsilon)\wedge s})\,dt- \frac{1}{\varepsilon}  \int_{[0,\,s[} F(Y_t)\,dt\nonumber\\ 
&=	\frac{1}{\varepsilon}  \int_{[\varepsilon,\,s+ \varepsilon[} F(Y_{t  \wedge s})\,dt- \frac{1}{\varepsilon}  \int_{[0,\,s[} F(Y_t)\,dt\nonumber\\ 
&=	\frac{1}{\varepsilon}  \int_{[s,\,s+ \varepsilon[} F(Y_{t  \wedge s})\,dt- \frac{1}{\varepsilon}  \int_{[0,\,\varepsilon[} F(Y_{t})\,dt\nonumber\\ 
&=	 F(Y_s)- \frac{1}{\varepsilon}  \int_{[0,\,\varepsilon[} F(Y_t)\,dt\nonumber\\ 
&\underset{\varepsilon \rightarrow 0}{\longrightarrow} F(Y_s)- F(Y_0),\,\, \textup{uniformly in $s$.} 
\end{align} 
We define $A(\varepsilon,N)$ and $B(\varepsilon,N)$ as in \eqref{Aepsilon}-\eqref{Bepsilon}.
 $J_0(\varepsilon,\,s)$ can be also rewritten as 
 \begin{equation} 
 \label{J0_passtolimit} 
 J_0(\varepsilon,\,s) = J_A(\varepsilon,\,N,\,s) + J_B(\varepsilon,\,N,\,s), 
 \end{equation} 
 where 
\begin{align} 
J_A(\varepsilon,\,N,\,s) &=  \frac{1}{\varepsilon} \int_0^s [F((t + \varepsilon)\wedge s,\,X_{(t + \varepsilon)\wedge s})- F(t,\,X_t)]\,\one_{A(\varepsilon,N)}(t)\,dt,\label{JA}\\ 
J_B(\varepsilon,\,N,\,s) &= \frac{1}{\varepsilon} \int_0^s [F((t + \varepsilon)\wedge s,\,X_{(t + \varepsilon)\wedge s})- F(t,\,X_t)]\,\one_{B(\varepsilon,N)}(t)\,dt.\label{JB} 
\end{align} 
Applying  Lemma \ref{L_ucp_big_jumps} with $n =2$ to  $Y=(Y^1,Y^2)=(t,X)$ and $\phi(y_1,y_2)=F(y_1^1,y_1^2)-F(y_2^1,y_2^2)$, we have 
\begin{align}\label{J1_conv} 
J_A(\varepsilon,\,N,\,s) &= \sum_{i=1}^N \frac{1}{\varepsilon} \int_{t_i - \varepsilon}^{t_i} [F((t + \varepsilon)\wedge s,\,X_{(t + \varepsilon)\wedge s})- F(t,\,X_t)]\,dt\nonumber\\ 
&\quad  \limit^{\varepsilon \rightarrow 0} \,\,\sum_{i=1}^N  \one_{]0,\,s]}(t_i)\,[F(t_i,\,X_{t_i})-F(t_i,\,X_{t_{i-}})],\,\, 
\textup{uniformly in $s$}. 
\end{align} 
Concerning  $J_B(\varepsilon,\,N,\,s)$, it can be decomposed into the sum of the two terms 
\begin{align*} 
J_{B1}(\varepsilon,\,N,\,s)&=\frac{1}{\varepsilon} \int_0^s [F((t + \varepsilon)\wedge s,\,X_{(t + \varepsilon)\wedge s})- F(t,\,X_{(t + \varepsilon)\wedge s})]\,\one_{B(\varepsilon,N)}(t)\,dt,\\ 
J_{B2}(\varepsilon,\,N,\,s)&= \frac{1}{\varepsilon} \int_0^s [F(t,\,X_{(t + \varepsilon)\wedge s})- F(t,\,X_t)]\,\one_{B(\varepsilon,N)}(t)\,dt. 
\end{align*} 
By finite increments theorem   in time we get  
\begin{equation}\label{JB1_decomp} 
	J_{B1}(\varepsilon,\,N,\,s) = J_{B10}(\varepsilon,\,s) + J_{B11}(\varepsilon,\,N,\,s)+ J_{B12}(\varepsilon,\,N,\,s)+J_{B13}(\varepsilon,\,N,\,s), 
\end{equation} 
where 
\begin{align*} 
	J_{B10}(\varepsilon,\,s) &=  \int_{0}^{s} \partial_t F(t,\,X_t)\,\frac{(t+ \varepsilon)\wedge s-t}{\varepsilon}\,dt,\\ 
	J_{B11}(\varepsilon,\,N,\,s) &= - 
	\sum_{i=1}^N \int_{t_i - \varepsilon}^{t_i} 
	\partial_t F(t,\,X_t)\,\frac{(t+ \varepsilon)\wedge s-t}{\varepsilon}\,dt,\\ 
	J_{B12}(\varepsilon,\,N,\,s) &= \int_{0}^{s}R_1(\varepsilon,\,t,\,s)\,\one_{B(\varepsilon,N)}(t)\, \frac{(t+ \varepsilon)\wedge s-t}{\varepsilon}\,dt,\\ 
	J_{B13}(\varepsilon,\,N,\,s) &= \int_{0}^{s} R_2(\varepsilon,\,t,\,s)\,\one_{B(\varepsilon,N)}(t)\,\frac{(t+ \varepsilon)\wedge s-t}{\varepsilon}\,dt, 
\end{align*} 
and 
\begin{align} 
R_1(\varepsilon,\,t,\,s) &= \int_0^1 [\partial_{t} F(t + a \,((t +\varepsilon)\wedge s -t),\,X_{(t + \varepsilon)\wedge s})- \partial_{t} F(t,\,X_{(t + \varepsilon)\wedge s})]\,da,\label{Rest_1}\\ 
R_2(\varepsilon,\,t,\,s) &= \partial_{t} F(t,\,X_{(t + \varepsilon)\wedge s})- \partial_{t} F(t,\,X_t).\label{Rest_2} 
\end{align} 
A Taylor expansion  in space up to second order  gives 
\begin{equation}\label{JB2_decomp} 
	J_{B2}(\varepsilon,\,N,\,s) = J_{B20}(\varepsilon,\,s) + J_{B21}(\varepsilon,\,s)+ J_{B22}(\varepsilon,\,N,\,s)+ J_{B23}(\varepsilon,\,N,\,s), 
\end{equation} 
where 
\begin{align} 
	J_{B20}(\varepsilon,\,s) &= \frac{1}{\varepsilon} \int_{0}^{s} \partial_x F(t,\,X_t)\,(X_{(t + \varepsilon)\wedge s}-X_t)\,dt,\label{appr_forward_int}\\ 
	J_{B21}(\varepsilon,\,s) &= \frac{1}{\varepsilon} \int_{0}^{s} \frac{\partial_{xx}^2 F(t,\,X_t)}{2}\,(X_{(t + \varepsilon)\wedge s}-X_t)^2\,dt,\nonumber\\ 
	J_{B22}(\varepsilon,\,N,\,s) &= -\frac{1}{\varepsilon} 
	\sum_{i=1}^N \int_{t_i - \varepsilon}^{t_i} 
	\left[\partial_x F(t,\,X_t)\,(X_{(t + \varepsilon)\wedge s}-X_t)+ \frac{\partial_{xx}^2 F(t,\,X_t)}{2}\,(X_{(t + \varepsilon)\wedge s}-X_t)^2\right]\,dt,\nonumber\\ 
	J_{B23}(\varepsilon,\,N,\,s) &=  \int_{0}^{s} R_3(\varepsilon,\,t,\,s)\,\one_{B(\varepsilon,N)}(t)\,\frac{(X_{(t + \varepsilon)\wedge s}-X_t)^2}{\varepsilon}\,dt,\nonumber 
\end{align} 
and 
\begin{equation} 
R_3(\varepsilon,\,t,\,s) = \int_0^1  [\partial_{xx}^2 F(t,\,X_t + a (X_{(t + \varepsilon)\wedge s}- X_t))- \partial_{xx}^2 F(t,\,X_t)]\,da.\label{Rest_3} 
\end{equation} 
Let us consider the term  $J_{B22}(\varepsilon,\,N,\,s)$. 
Applying   Lemma \ref{L_ucp_big_jumps} with $n =2$ to $Y=(Y^1,Y^2)=(t,X)$ and $\phi(y_1,y_2)= \partial_x F(y_2^1,y_2^2)(y_1^2-y_2^2)+ \partial_{xx}^2 F(y_2^1,y_2^2)(y_1^2-y_2^2)^2$, 
 we get 
\begin{equation}\label{J22_conv} 
\lim_{\varepsilon \rightarrow 0} J_{B22}(\varepsilon,\,N,\,s) = 
-\sum_{i=1}^N \one_{]0,\,s]}(t_i)\,\left[\partial_x F(t_i,\,X_{t_{i-}})\,(X_{t_{i}}-X_{t_{i}-})+ \frac{\partial_{xx}^2 F(t_i,\,X_{t_{i-}})}{2}(X_{t_{i}}-X_{t_{i}-})^2\right] 
\end{equation} 
uniformly in $s$. 
Moreover, the term  $J_{B10}(\varepsilon,\,N,\,s)$ can be decomposed into
\begin{equation}\label{J_B10} 
J_{B10}(\varepsilon,\,s) = \int_0^s \partial_t F(t,\,X_t)\, dt+ J_{B10'}(\varepsilon,\,s) +J_{B10''}(\varepsilon,\,s), 
\end{equation} 
with 
\begin{align} 
J_{B10'}(\varepsilon,\,s) &= \int_{s-\varepsilon}^{s}\partial_{t} F(t,X_t)\,\frac{s-t}{\varepsilon}\, dt,\label{J_B10'}\\ 
J_{B10''}(\varepsilon,\,s)&= -\int_{s-\varepsilon}^{s}\partial_{t} F(t,X_t)\, dt.\label{J_B10''} 
\end{align} 
 
At this point we remark that identity \eqref{J0_passtolimit} can be rewritten as 
\begin{align}\label{Jovarepsilo_BIS} 
J_0(\varepsilon,s)&=J_A(\varepsilon,N,s) +\int_0^s \partial_t F(t,\,X_t)\, dt+ J_{B10'}(\varepsilon,s)+J_{B10''}(\varepsilon,s)+J_{B11}(\varepsilon,N,s) 
+J_{B12}(\varepsilon,N,s)\nonumber\\&+J_{B13}(\varepsilon,N,s) 
+J_{B20}(\varepsilon,s)+J_{B21}(\varepsilon,s)+J_{B22}(\varepsilon,N,s)+J_{B23}(\varepsilon,N,s). 
\end{align} 
Passing 
to the limit in \eqref{Jovarepsilo_BIS} on both the left-hand and right-hand sides, uniformly in $s$,  as $\varepsilon$ goes to zero,  
taking into account   convergences \eqref{J0_conv}, \eqref{J1_conv}, \eqref{J22_conv}. 
we get 
\begin{align} 
F(s,\,X_s)- F(0,\,X_0) &= \int_0^s \partial_t F(t,\,X_t)\, dt + \sum_{i = 1}^N \one_{]0,\,s]}(t_i)\,\left[ F(t_i,\,X_{t_i})-F(t_i,\,X_{t_{i-}}) \right]\nonumber \\ 
&-\sum_{i = 1}^N \one_{]0,\,s]}(t_i)\,\left[ \partial_x F(t_i,\,X_{t_{i}-})\,(X_{t_{i}}-X_{t_{i}-})- \frac{\partial_{xx}^2 F(t_i,\,X_{t_{i}-})}{2}(X_{t_{i}}-X_{t_{i}-})^2\right]\nonumber\\ 
& + \lim_{\varepsilon \rightarrow 0}\, ( 
J_{B20}(\varepsilon,\,N, \,s)  
+J_{B21}(\varepsilon, \,s)+L(\varepsilon, N,s)) 
\label{collect} 
\end{align} 
where the previous limit is intended uniformly in $s$, and we have set 
\begin{align*} 
L(\varepsilon, N,s) &:=J_{B10'}(\varepsilon,\,s)+ J_{B10''}(\varepsilon,\,s)+J_{B11}(\varepsilon,\,N,\,s) +J_{B12}(\varepsilon,\,N, \,s) \\ 
&+ J_{B13}(\varepsilon, \,N,\,s)  + J_{B23}(\varepsilon,\,N, \,s). 
\end{align*} 
We evaluate previous limit uniformly in $s$, for every 
 $\omega \notin 
 \mathcal N$. Without restriction of generality it is enough to show the   uniform convergence in $s$ 
for the subsequence $\varepsilon_n$ introduced in \eqref{A2}-\eqref{A1}, when $n \rightarrow \infty$. 
 
According to \eqref{A1}, 
we get 
\begin{equation} \label{J21_conv} 
\lim_{n \rightarrow \infty} J_{B21}(\varepsilon_n, \,s) = \int_{]0,\,s]}\frac{\partial_{xx}^2 F(t,\,X_{t-})}{2}\,d[X,X]_t, 
\end{equation} 
 uniformly in $s$.

We now should discuss $J_{B12}(\varepsilon_n,\,N,\,s)$, $J_{B13}(\varepsilon_n,\,N,\,s)$ and $J_{B23}(\varepsilon_n,\,N,\,s)$. 
In the sequel,  $\delta(f,\cdot)$ will denote  the modulus of continuity of a function $f$, and by  $I_l$ the interval $[t_{l-1},\,t_l]$, $l \geq 0$. 
Since $\frac{(t+\varepsilon)\wedge s -t}{\varepsilon} \leq 1$ for every $t, s$, 
by Remark \ref{R_Billingsley} we get 
\begin{align*} 
\one_{B(\varepsilon,N)}(t)\,|R_1(\varepsilon,\,t,\,s)| 
\leq& \delta\left( \partial_{t} F,\, \varepsilon \right),\\ 
\one_{B(\varepsilon,N)}(t)\,|R_2(\varepsilon,\,t,\,s)| 
\leq& \delta\Big( \partial_{t} F,\, \sup_l\sup_{\underset{\vert t - a\vert \le \varepsilon}{t, a \in I_l}} 
|X_{a}- X_t|) \Big),\\ 
\one_{B(\varepsilon,N)}(t)\,|R_3(\varepsilon,\,t,\,s)| 
\leq& \delta\Big( \partial_{xx}^2 F,\, \sup_l\sup_{\underset{\vert t - a\vert \le \varepsilon}{t, a \in I_l}} 
|X_{a}- X_t|) \Big). 
\end{align*} 
Considering the two latter inequalities, Lemma \ref{Lem_Billingsley} applied successively to the intervals $I_l$ 
implies 
\begin{align*} 
\one_{B(\varepsilon,N)}(t)\,|R_2(\varepsilon,\,t,\,s)| &\leq \delta(\partial_{t} F,\,3 \gamma),\\ 
\one_{B(\varepsilon,N)}(t)\,|R_3(\varepsilon,\,t,\,s)| &\leq \delta(\partial_{xx}^2 F,\,3 \gamma). 
\end{align*} 
Then, using again  $\frac{(t+\varepsilon_n)\wedge s-t}{\varepsilon} \leq1$, we get 
\begin{align} 
\sup_{s \in [0,\,T]}|J_{B12}(\varepsilon_{n},\,N,\,s)| &\leq \delta(\partial_{t} F,\,\varepsilon_n)\cdot T, 
\nonumber \\ 
\sup_{s \in [0,\,T]}|J_{B13}(\varepsilon_{n},\,N,\,s)| &\leq \delta(\partial_{t} F,\,3 \gamma)\cdot T, 
\nonumber\\ 
\sup_{s \in [0,\,T]}|J_{B23}(\varepsilon_{n},\,N,\,s)| &\leq \delta(\partial_{xx}^2 F,\,3 \gamma)\cdot 
 \sup_{n \in \N, s \in [0,\,T]} [X,X]_{\varepsilon_{n}}^{ucp}(s), \label{ESupremumBracket} 
\end{align} 
where we remark that the supremum in the right-hand side of \eqref{ESupremumBracket} is finite taking into 
account \eqref{A2}. 
Therefore 
\begin{align} 
\limsup_{n \rightarrow \infty}\sup_{s \in [0,\,T]}|J_{B23}(\varepsilon_{n},\,N,\,s)| 
	&=\delta(\partial_{xx}^2 F,\,3 \gamma)\cdot \!\!\!\!\!\!\sup_{n \in \N, s \in [0,\,T]} [X,X]_{\varepsilon_{n}}^{ucp}(s),\label{J23_limsup}\\ 
\limsup_{n \rightarrow \infty}\sup_{s \in [0,\,T]}|J_{B13}(\varepsilon_{n},\,N,\,s)| 
&=\delta(\partial_{t} F,\,3 \gamma)\cdot T \label{J13_limsup}, 
\end{align} 
while 
\begin{equation}\label{J12_conv} 
\lim_{n \rightarrow \infty}\sup_{s \in [0,\,T]}|J_{B12}(\varepsilon_{n},\,N,\,s)|= 0. 
\end{equation} 
Let now  consider the terms  $J_{B10'}(\varepsilon_{n},\,s)$, $J_{B10''}(\varepsilon_{n},\,s)$ and  $J_{B11}(\varepsilon_{n},\,N,\,s)$. 
\begin{align*} 
\sup_{s \in [0,\,T]}|J_{B10'}(\varepsilon_{n},\,s)| &\leq \sup_{y \, \in \,\mathbb{K}^X(\omega) \times [0,\,T]}|\partial_t F(y)|\cdot \varepsilon_n,\\ 
\sup_{s \in [0,\,T]}|J_{B10''}(\varepsilon_{n},\,s)| &\leq \sup_{y \, \in \, \mathbb{K}^X(\omega)\times [0,\,T]}|\partial_t F(y)|\cdot \varepsilon_n,\\ 
\sup_{s \in [0,\,T]}|J_{B11}(\varepsilon_{n},\,N,\,s)| &\leq \sup_{y \, \in \, \mathbb{K}^X(\omega) \times [0,\,T]}|\partial_t F(y)|\,N\cdot\varepsilon_n, 
\end{align*} 
where $\mathbb{K}^X(\omega)$ is the (compact)  set $\{X_t(\omega),\, t \in [0,\,T]\}$. 
So, it follows 
\begin{equation}\label{J10primes_conv} 
\lim_{n \rightarrow \infty}\sup_{s \in [0,\,T]}|J_{B10'}(\varepsilon_{n},\,s)|=\lim_{n \rightarrow \infty}\sup_{s \in [0,\,T]}|J_{B10''}(\varepsilon_{n},\,s)|=\lim_{n \rightarrow \infty}\sup_{s \in [0,\,T]}|J_{B11}(\varepsilon_{n},\,N,\,s)|=0. 
\end{equation} 
Taking into account \eqref{J10primes_conv}, \eqref{J13_limsup}, \eqref{J23_limsup}, and \eqref{J21_conv}, we see that  
\begin{equation}\label{L_conv} 
\limsup_{n \rightarrow \infty}\sup_{s \in [0,\,T]} |L(\varepsilon_n, N,s)| = \delta(\partial_{xx}^2 F,\,3 \gamma)\cdot \!\!\!\!\!\!\sup_{n \in \N, s \in [0,\,T]} [X,X]_{\varepsilon_{n}}^{ucp}(s) + \delta(\partial_{t} F,\,3 \gamma)\cdot T. 
\end{equation} 
 
 
Recalling that $J_{B20}(\varepsilon,s)$ in \eqref{appr_forward_int} is the $\varepsilon$-approximation of the forward integral $\int_0^t\partial_x F(s,X_{s})\,d^{-}X_s$, to conclude  
it remains to show that 
\begin{eqnarray}\label{final} 
&&\sup_{s \in [0,\,T]}\big| J_{B20}(\varepsilon_{n},\,s) 
- J(s)\big| 
\limit_{n \rightarrow \infty} 0\quad \textup{a.s.}, 
\end{eqnarray} 
where 
\begin{align}\label{J0_final} 
J(s)&=F(s,\,X_s)- F(0,\,X_0)- \int_{]0,\,s]}\partial_t F(t,\,X_t)\,dt -\sum_{t \leq s} \left[F(t,\,X_t)-F(t,\,X_{t-})\right]\nonumber\\ 
& +\sum_{0 <t \leq s} \left[ \partial_x F(t,\,X_{t-})\,(X_{t}-X_{t-}) +\frac{\partial_{xx}^2 F(t,\,X_{t-})}{2}\,\,(X_{t}-X_{t-})^2 \right] 
\nonumber\\ 
&- \frac{1}{2} \int_{]0,\,s]} \partial_{xx}^2 F(t,\,X_{t-})\,d[X,X]_{t}. 
\end{align} 
In particular this would  imply that 
$\int_{]0,\,s]}\partial_x F(t,\,X_t)\,d^{-}X_t$ 
 exists and equals $J(s)$. 
Taking into account \eqref{Jovarepsilo_BIS}, we have  
\begin{align}\label{JB2_final} 
J_{B20}(\varepsilon_n,\,s) 
 	=J_0(\varepsilon_n, s) -J_A(\varepsilon_{n},\,N,\,s) 
 	-\int_0^s \partial_t F(t,\,X_t)\, dt 
 	-L(\varepsilon_n,N,s) 
 	- J_{B21}(\varepsilon_n,\,s) 
 	 -J_{B22}(\varepsilon_n,\,N,\,s). 
\end{align} 
Taking into account \eqref{J0_final} and \eqref{JB2_final}, 
we see that the term inside the absolute value in \eqref{final} equals 
\begin{align*} 
& J_{0}(\varepsilon_{n},\,s) - (F(s,\,X_s)- F(0,\,X_0))\\ 
&-J_{A}(\varepsilon_{n},\,N,\,s)+ \sum_{i=1}^N \one_{]0,\,s]}(t_i) 
[ F(t_i,\,X_{t_i})-F(t_i,\,X_{t_{i-}})]\\ 
&  -J_{B22}(\varepsilon_{n},\,N,\,s) 
- \sum_{i=1}^N \one_{]0,\,s]}(t_i) 
\left[\partial_x F(t_i,\,X_{{t_i}-})\,(X_{t_i}-X_{t_{i-}}) + \frac{\partial_{xx}^2 F(t_i,\,X_{{t_i}-})}{2}\,(X_{t_i}-X_{t_{i}-})^2 \right]\\ 
& -J_{B21}(\varepsilon_{n},\,s) +\frac{1}{2} \int_{]0,\,s]} \partial_{xx}^2 F(t,\,X_{t-})\,d[X,X]_t \\ 
&- L(\varepsilon_n,N,s)\\ 
& + \sum_{i=N+1}^{\infty} \one_{]0,\,s]}(t_i)\left[F(t_i,\,X_{t_i})-F(t_i,\,X_{t_{i-}}) - \partial_x F(t_i,\,X_{{t_i}-})\, (X_{t_i}-X_{t_{i-}}) - \frac{\partial_{xx}^2 F(t_i,\,X_{{t_i}-})}{2}\,(X_{t_i}-X_{t_{i-}})^2 \right]. 
\end{align*} 
Taking into account   \eqref{J0_conv}, \eqref{J1_conv}, \eqref{J22_conv}, \eqref{J12_conv}, \eqref{L_conv},%
we have 
\begin{align}\label{last_ineq} 
	&\limsup_{n \rightarrow \infty}\sup_{s\in [0,\,T]}\bigg |J_{B20}(\varepsilon_{n},\,s)- 
	J(s) \bigg|\nonumber\\ 
	&\leq 
	\limsup_{n \rightarrow \infty}\sup_{s\in [0,\,T]} |L(\varepsilon_{n},\,N,\,s)|\nonumber\\ 
	 &
	 +\sup_{s\in [0,\,T]} 
	 \sum_{i=N+1}^{\infty} \one_{]0,\,s]}(t_i)\,\bigg|F(t_i,\,X_{t_i})-F(t_i,\,X_{t_{i-}}) - \partial_x F(t_i,\,X_{{t_i}-})\,\Delta \,X_{t_i} - \frac{\partial_{xx}^2 F(t_i,\,X_{{t_i}-})}{2}\,(\Delta \, X_{t_{i}})^2 \bigg|\nonumber\\ 
	&= 
	\limsup_{n \rightarrow \infty}\sup_{s\in [0,\,T]} |L(\varepsilon_{n},\,N,\,s)| \nonumber \\ 
	 &+ 
	 \sup_{s\in [0,\,T]}\sum_{i=N+1}^\infty (\Delta X_s)^2\, 
	\one_{]0,\,s]}(t_i)\,\frac{1}{2}\,\Big|\int_0^1\partial_{xx}^2 F(t_i,\,X_{{t_i}-} + a(\Delta X_{t_i}))\,da-\partial_{xx}^2 F(t_i,\,X_{{t_i}-})\Big|\nonumber\\ 
	&\leq 
	\delta(\partial_{t} F,\,3 \gamma)\cdot T+\delta(\partial_{xx}^2 F,\,3 \gamma)\cdot \!\!\!\!\!\!\sup_{n \in \N, s \in [0,\,T]} [X,X]_{\varepsilon_{n}}^{ucp}(s) +\gamma^2 \sup_{y \in \mathbb{K}^X(\omega) \times [0,\,T]}|\partial_{xx}^2 F(y)|, 
\end{align} 
where  the last term on the right-hand side of \eqref{last_ineq}  is obtained  using \eqref{gamma_N+1_condition}. 
Since $\gamma$ is arbitrarily small, 
we conclude that 
\begin{displaymath} 
\lim_{n \rightarrow \infty}\sup_{s\in [0,\,T]}\bigg |J_{B20}(\varepsilon_{n},\,s)- 
J(s)\bigg|=0, \quad \forall \omega \notin \mathcal{N}. 
\end{displaymath} 
This concludes the proof of the It\^o formula. 
\endproof 
From Theorem \ref{P_Ito_C2_cadlag}, Proposition \ref{P_Ito_cov}-ii), and by classical Banach-Steinhaus theory (see, e.g.,  \cite{ds}, Theorem 1.18 p. 55) for $F$-type spaces, we have the following. 
\begin{proposition}\label{P_RV_Itoformula_C1lambda} 
	Let $F:[0,\,T]\times \R \rightarrow \R$ be a function of class $C^1$ such that $\partial_xF$ is H\"older continuous with respect to the second variable for some $\lambda \in [0,\,1[$. 
	Let $(X_t)_{t \in [0,\,T]} $ be a reversible semimartingale, satisfying moreover 
	\[ 
	\sum_{0 <s \leq t}|\Delta X_s|^{1+ \lambda} < \infty\quad \textup{a.s.} 
	\] 
	Then 
	\begin{equation*}\label{Ito_formula_C1lambda} 
	F(t,X_t)=F(0,X_0)+\int_0^t\partial_s F(s,X_s)\,ds +\int_0^t\partial_x F(s,X_{s-})\,d X_s+\frac{1}{2}\,[\partial_{x}F(\cdot,X),\,X]_t + J(F,X)(t), 
	\end{equation*} 
	where 
	\begin{equation*}\label{JFX} 
	 J(F,X)(t) = \sum_{0 <s \leq t}\left[F(s,X_s) - F(s,X_{s-})  - \frac{\partial_x F(s,X_s)+\partial_x 
F(s,X_{s-})}{2}\,\Delta X_s\right] 
	\end{equation*} 
\end{proposition} 
\begin{remark}\label{R_Itoformula_C1lambda} 
	\begin{itemize} 
		\item[(i)]Previous result can be   easily extended to the case when $X$ is multidimensional. 
		\item[(ii)] When $F$ does not depend on time,  previous statement was the object of 
 \cite{erv}, Theorem 3.8, example 3.3.1. In that case however, stochastic integrals and covariations 
were defined  by discretizations means. 
		\item[(iii)] The proof of Proposition \ref{P_RV_Itoformula_C1lambda} follows the same lines as the one of  Theorem 3.8. in \cite{erv}. 
	\end{itemize} 
\end{remark}

 \section{It\^{o}'s formula for $C^{1,2}$ functions  related to random measures}\label{S32}
 
The concept of random measure allows  a very tractable description of the jumps of a c\`adl\`ag process. 
The object of the present section is to reexpress the statement 
of Theorem \ref{P_Ito_C2_cadlag} making use of the 
jump measure 
associated with a c\`adl\`ag process $X$.

We start by recalling 
 the main definitions and some  properties that we will extensively use in the sequel; for a complete discussion on this topic  and the unexplained  notations we refer to   Chapter II, Section 1, in \cite{JacodBook},  Chapter XI, Section 1, in \cite{chineseBook}, and also to the Appendix in \cite{BandiniRusso2}. 
We set  $\tilde \Omega= \Omega \times [0,T] \times \R$, and $\tilde{\mathcal{P}} =\mathcal{P} \otimes \mathcal{B}(\R)$, 
where $\mathcal{P}$
is the predictable 
$\sigma$-field 
on $\Omega \times [0,T]$. 
A function $W$ defined on $\tilde \Omega$ which is $\tilde{\mathcal{P}}$-measurable will be 
called predictable.
We will  also indicate by    $\mathcal{A}$ (resp $\mathcal{A}_{\textup{loc}}$)  the collection of all adapted processes with   integrable variation (resp.  with locally integrable variation), and by    $\mathcal{A}^+$ (resp $\mathcal{A}_{\textup{loc}}^+$)  the collection of all adapted integrable increasing (resp. adapted locally integrable)  processes. 
The significance of {\rm locally} is the usual one,   which refers 
to  localization by stopping times, see e.g. (0.39) of \cite{jacod_book}. We only remark that our stopping times will take values in $[0,\,T] \cup \{+ \infty\}$ instead of 
$[0,\, \infty]$. For instance, adapting the definition of locally bounded process stated before Theorem 15, Chapter IV, in \cite{protter}, to the processes indexed by $[0,\,T]$, we can state  the following. 

\begin{definition}\label{D_locally_bounded} 
A process $(X_t)_{t \in [0,\,T]}$ is locally bounded if there exists a sequence of stopping times $(\tau_n)_{n \geq 1}$ in $[0,\,T] \cup \{+\infty\}$  increasing to $\infty$ a.s., such that 
		$(X_{\tau_n \wedge t}\,\one_{\{\tau_n>0\}})_{t \in [0,\,T]}$ is bounded. 
\end{definition} 
\begin{remark}\label{R_localisation} 
	\begin{itemize} 
		\item[(i)] 
		Any c\`agl\`ad
		 process is locally bounded, see the lines above Theorem 15, Chapter IV, in \cite{protter}. 
		\item[(ii)]Let $X$ be a c\`adl\`ag process satisfying condition \eqref{StandAss}. 
			Set $(Y_t)_{t \in [0,\,T]}=$ $(X_{t-},\sum_{s <t} |\Delta X_s|^2)_{t \in [0,\,T]}$. The process $Y$ is c\`agl\`ad, therefore  locally bounded by item (i). In particular,  we can fix a  sequence of stopping times $(\tau_n)_{n \geq 1}$ in $[0,\,T] \cup \{+\infty\}$  increasing to $\infty$ a.s., such that 	$(Y_{\tau_n \wedge t}\,\one_{\{\tau_n>0\}})_{t \in [0,\,T]}$ is bounded.  
	\end{itemize} 
\end{remark}

For any $X= (X_t)$  adapted  real valued c\`adl\`ag process 
on $[0,\,T]$, we call  \textbf{jump measure} of $X$ 
the integer-valued random measure on $\R_+ \times \R$ defined as 
\begin{equation}\label{X_jump_measure} 
	\mu^X(\omega; dt\,dx):= \sum_{s\in ]0,\,T]} \,\one_{\{\Delta X_s(\omega)\neq 0\}}\,\delta_{(s,\,\Delta X_s(\omega))}(dt\,dx). 
\end{equation} 
The compensator of 	$\mu^X(ds\,dy)$ 
is called 
the L\'evy system of $X$, and will be denoted by $\nu^X(ds\,dy)$;
we also set  
\begin{equation}\label{hatnu} 
\hat{\nu}^X_t = \nu^X(\{t\},\,dy)\quad  \textup{for every} \,\,t \in [0,\,T]. 
\end{equation}

\begin{remark} \label{R_Jumps_measure} 
	The jump measure $\mu^X$ 
	acts in the following way: for any positive  function $W: \Omega \times [0,\,T] \times \R\rightarrow \R$,  measurable  with respect to the product sigma-field,
	we have 
	\begin{equation}\label{muX_def}
		\sum_{s \in ]0,\,T]}\one_{\{\Delta X_s \neq 0\}} W_s(\cdot,\Delta X_s) = \int_{]0,T] \times \R} W_s(\cdot, x) \, \mu^X(\cdot, ds \, dx). 
	\end{equation}
	Clearly, the same result holds if the left or the right hand-side of \eqref{muX_def} is finite a.s., replacing the process $W$ by the process $|W|$.
\end{remark}

\begin{proposition}\label{P_X_jump_comment} 
	Let $X$ be a  c\`adl\`ag process on $[0,\,T]$ 
	satisfying condition \eqref{StandAss}, 
	and let  $F$ be a function of class $C^{1,2}$. 
	Then 
	$$ 
	|(F(s,X_{s-} + x)-F(s,X_{s-})- x\,\partial_x F(s,X_{s-})|\,\one_{\{|x|\leq 1\}}\ast \mu^X \in \mathcal{A}_{\textup{loc}}. 
	$$ 
\end{proposition} 
\proof 
Let $(\tau_n)_{n \geq 1}$ be the sequence of stopping times  
introduced in Remark \ref{R_localisation}-(ii) for the 
 process $Y_t = (X_{t-},  \sum_{s < t}|\Delta X_{s}|^2)$. Fix $\tau = \tau_n$, and  let $M$ such that $\sup_{t \in [0,\,T]}|Y_{t \wedge \tau}\,\one_{\{\tau >0 \}}|\leq M$. 
So, by an obvious Taylor expansion, taking into account Remark 
\ref{R_Jumps_measure}, we have 
\begin{align*}  
	&\sper{\int_{]0,\,t \wedge \tau]\times \R}	|(F(s,X_{s-} + x)-F(s,X_{s-})- x\,\partial_x F(s,X_{s-})|\,\one_{\{\vert x \vert \leq 1\}}\,\mu^X(ds,\,dx)}\nonumber\\ 
	&=\sper{\sum_{0 <s \leq t \wedge \tau}\left[F(s,X_s)-F(s,X_{s-})- \partial_x F(s,X_{s-})\,\Delta X_s\right]}\\ 
	&=\sper{\sum_{0<s \leq t\wedge \tau}(\Delta X_s)^2\,\one_{\{\tau>0\}}\,\frac{1}{2}\,\int_0^1
	\partial_{xx}^2 F(s,X_{s-}+ a\,\Delta X_s)
	\,da} \\ 
	&\leq \frac{1}{2} \sup_{ \underset{t \in [0,\,T]}{y \in  [-M,\,M]}}|\partial^2_{xx} F|(t,\,y)\, \sper{\sum_{0<s < {t \wedge \tau}}|\Delta X_s|^2 \,\one_{\{|\Delta X_s| \leq 1\}}\,\one_{\{\tau>0\}} +|\Delta X_{\tau}|^2 \,\one_{\{|\Delta X_{\tau}| \leq 1\}}\,\one_{\{\tau>0\}}}\nonumber\\ 
	&\leq \frac{1}{2} \sup_{ \underset{t \in [0,\,T]}{y \in  [-M,\,M]}}|\partial^2_{xx}  F|(t,\,y)\, \cdot (M+1), 
\end{align*} 
which concludes the proof. 
\endproof 
\begin{proposition}\label{P_X_jump_comment2} 
	Let $X$ be a  c\`adl\`ag process on $[0,\,T]$ satisfying condition \eqref{StandAss}, 
	and let  $F$ be a function of class $C^{0,1}$. 
	Then 
	\begin{align} 
		&|(F(s,X_{s-} + x)-F(s,X_{s-})|^2\,\one_{\{\vert x \vert \leq 1\}}\ast \mu^X \in \mathcal{A}_{\textup{loc}}^+,\label{A_Aloc}\\ 
		&|x\,\partial_x F(s,X_{s-})|^2\,\one_{\{\vert x \vert \leq 1\}}\ast \mu^X \in \mathcal{A}_{\textup{loc}}^+.\label{K_Aloc} 
	\end{align} 
In particular, the integrands in \eqref{A_Aloc}-\eqref{K_Aloc} belong to $\mathcal G^2_{\rm loc}(\mu^X)$, so that the stochastic integrals 
	\begin{align} 
		&((F(s,X_{s-} + x)-F(s,X_{s-}))\,\one_{\{\vert x \vert \leq 1\}}\ast (\mu^X-\nu^X),\label{A_mart}\\ 
		&x\,\partial_x F(s,X_{s-})\,\one_{\{\vert x \vert \leq 1\}}\ast (\mu^X-\nu^X)\label{K_mart} 
	\end{align}
	are well-defined and are two square integrable purely discontinuous martingales.
\end{proposition} 
\proof 
Proceeding as  in the proof of Proposition \ref{P_X_jump_comment}, we 
consider the  sequence  of stopping times $(\tau_n)_{n \geq 1}$  
defined  in Remark \ref{R_localisation}-(ii) for the  
process $Y_t = (X_{t-},  \sum_{s < t}|\Delta X_{s}|^2)$.  Fix $\tau = \tau_n$, and  let $M$ such that $\sup_{t \in [0,\,T]}|Y_{t \wedge \tau}\,\one_{\{\tau >0 \}}|\leq M$. 
For any $t \in [0,T]$, we have 
\begin{align*}  
	&\sper{\int_{]0,\,t \wedge \tau]\times \R}|(F(s,X_{s-} + x)-F(s,X_{s-})|^2\, 
		\one_{\{\vert x \vert \leq 1\}}\,\mu^X(ds,\,dx)}\nonumber\\ 
	&\leq \sup_{ \underset{t \in [0,\,T]}{y \in  [-M,\,M]}}|\partial_{x} F|^2(t,\,y)\, \sper{\sum_{0<s < {t \wedge \tau}}|\Delta X_s|^2 \,\one_{\{|\Delta X_s| \leq 1\}} \one_{\{\tau>0\}}+|\Delta X_{\tau}|^2 \,\one_{\{|\Delta X_{\tau}| \leq 1\}}\,\one_{\{\tau>0\}}}\nonumber\\ 
	&\leq \sup_{ \underset{t \in [0,\,T]}{y \in  [-M,M]}}|\partial_{x}  F|^2(t,\,y)\, \cdot (M+1), 
\end{align*} 
and 
\begin{align*} 
	&\sper{\int_{]0,\,t \wedge \tau]\times \R} |x\,\partial_x F(s,X_{s-})|^2\, 
		\one_{\{\vert x\vert \leq 1\}}\,\mu^X(ds,\,dx)}\nonumber\\ 
	&=\sper{\int_{]0,\,t \wedge \tau]\times \R}|x|^2\,|\partial_{x} F|^2(t,\, X_{s-})\,\one_{\{|x| \leq 1\}}\,\mu^X(ds,\,dx)}\nonumber\\ 
	&\leq  \sup_{ \underset{t \in [0,\,T]}{y \in  [-M,\,M]}}|\partial_{x} F|^2(t,\,y)\, \sper{\sum_{0<s < {t \wedge \tau}}|\Delta X_s|^2 \,\one_{\{|\Delta X_s| \leq 1\}} \,\one_{\{\tau>0\}}+|\Delta X_{\tau}|^2 \,\one_{\{|\Delta X_{\tau}| \leq 1\}}\,\one_{\{\tau>0\}}}\nonumber\\ 
	&\leq \sup_{ \underset{t \in [0,\,T]}{y \in [-M,\,M]}}|\partial_{x} F|^2(t,\,y)\, 
	\cdot (M+1). 
\end{align*} 
By 
Lemma 2.4-2. in  \cite{BandiniRusso2}, the integrands in   \eqref{A_mart} and \eqref{K_mart} belong to $\mathcal G^2_{\rm loc}(\mu^X)$. 
Then 
the conclusion follows by Theorem 11.21, point 3), in \cite{chineseBook}.
\endproof 

Taking $F(t,x) = x$, we have the following.
 \begin{corollary}\label{R_x_smalljumps_G2(mu)} 
	Let $X$ be a c\`adl\`ag process satisfying condition \eqref{StandAss}. 
	Then 
	\begin{equation}\label{x_small_jumps_int} 
		x\,\one_{\{|x| \leq 1\}}\,\in \mathcal{G}^2_{\rm loc}(\mu^X), 
	\end{equation} 
	and the stochastic integral 
	\begin{equation}\label{x_small_jumps_mtg} 
		\int_{]0,\,t]\times \R}x\,\one_{\{|x| \leq 1\}}\,(\mu^X- \nu^X)(ds\,dx) 
	\end{equation} 
	is well-defined and defines a purely discontinuous square integrable local martingale. 
\end{corollary}

We state now the main result of the section. Even if the following decomposition  is not canonical, it will serve as a significant tool for the case when $F$ is of class $C^{0,1}$.
\begin{proposition}\label{C_reformulation_C12_Ito} 
	Let $X$ be a finite quadratic variation c\`adl\`ag process and $F: [0,\,T]\times \R \rightarrow \R $ a function of class $C^{1,2}$. Then we have 
	\begin{align}\label{Ito_formula_C2_jump_measure} 
	F(t,X_t)&=F(0,X_0)+\int_0^t\partial_s F(s,X_s)\,ds +\int_0^t\partial_x F(s,X_{s})\,d^{-}X_s+\frac{1}{2}\int_0^t\partial_{xx}^2 F(s,X_{s})\,d[X,X]_s^c \nonumber\\ 
	&+\int_{]0,\,t]\times \R}\,(F(s,X_{s-}+x)-F(s,X_{s-}))\,\one_{\{x \leq 1\}}\,(\mu^X- \nu^X)(ds\,dx)\nonumber\\ 
	&-\int_{]0,\,t]\times \R}\,x\,\partial_x F(s,X_{s-})\,\one_{\{x \leq 1\}}\,(\mu^X- \nu^X)(ds\,dx)\nonumber\\ 
	&+\int_{]0,\,t]\times \R}\,(F(s,X_{s-}+x)-F(s,X_{s-})- x\,\partial_x F(s,X_{s-}))\,\one_{\{x > 1\}}\,\mu^X(ds\,dx)\nonumber\\ 
	&+\int_{]0,\,t]\times \R}\,(F(s,X_{s-}+x)-F(s,X_{s-})- x\,\partial_x F(s,X_{s-}))\,\one_{\{x \leq 1\}}\,\nu^X(ds\,dx). 
	\end{align} 
\end{proposition}
 

\proof 
This result is a consequence of Theorem \ref{P_Ito_C2_cadlag}. Indeed, by Remark \ref{R_Jumps_measure}, the term 
$$
\sum_{s\leq t}\,[F(s,X_s)-F(s,X_{s-})- \partial_x F(s,X_{s-})\,\Delta \,X_{s}]
$$
equals
\begin{align}\label{interm}
&\int_{]0,\,t]\times \R}\,(F(s,X_{s-}+x)-F(s,X_{s-})- x\,\partial_x F(s,X_{s-}))\,\mu^X(ds\,dx)\nonumber\\
&= \int_{]0,\,t]\times \R}\,(F(s,X_{s-}+x)-F(s,X_{s-})- x\,\partial_x F(s,X_{s-}))\,\one_{\{x > 1\}}\,\mu^X(ds\,dx)
\nonumber\\
&+ \int_{]0,\,t]\times \R}\,(F(s,X_{s-}+x)-F(s,X_{s-})- x\,\partial_x F(s,X_{s-}))\,\one_{\{x \leq 1\}}\,\mu^X(ds\,dx).
\end{align}
	We set 
	\begin{align*} 
	W_s(x)&=  \left(F(s,X_{s-} + x)-F(s,X_{s-})- x\,\partial_x F(s,X_{s-}) 
	\right) \one_{\{\vert x \vert \leq 1\}}. 
	\end{align*} 
 By 
	Propositions \ref{P_X_jump_comment}, $\vert W \vert \ast \mu^X$ belongs to $\mathcal{A}^+_{\rm loc}$,  and consequently $\vert W \vert \ast\nu^X$ belongs to $\mathcal{A}^+_{\rm loc}$
	so that
the last integral in \eqref{Ito_formula_C2_jump_measure}  is well-defined.
		Setting
	\begin{align*}  
	K_s(x)&=\left(F(s,X_{s-} + x)-F(s,X_{s-})\right)\one_{\{\vert x\vert \leq 1\}},\\ 
 	Y_s(x)&=x\,\partial_x F(s,X_{s-})\one_{\{\vert x \vert \leq 1\}},
	\end{align*}
Proposition \ref{P_X_jump_comment2} insures that 
   $K \ast (\mu^X-\nu^X)$, $Y \ast (\mu^X-\nu^X)$ are well-defined and are purely discontinuous square integrable 
 local martingales. 
For those reasons, the second integral of the right-hand side of \eqref{interm} gives us the contributions of second,  third and fifth lines of \eqref{Ito_formula_C2_jump_measure}.
\endproof

\section{About weak Dirichlet processes}\label{Section_WDprocesses} 
 
\subsection{Basic definitions and properties} 
 We consider again the filtration 
 $(\mathcal{F}_t)_{t \geq 0}$ introduced at Section 
\ref{Sec_tecnicalities_JumpCalculus}. 
Without further mention, 
the underlying  filtration  will be indeed $(\mathcal{F}_t)_{t \geq 0}$.

\begin{definition} 
Let $X$ be an $({\mathcal F}_t)$-adapted process. 
	We say that  $X$  is  $({\mathcal F}_t)$-orthogonal if 
	$[X, N] = 0$  for every $N$ continuous local $({\mathcal F}_t)$-martingale. 
\end{definition} 
 
 
\begin{remark}\label{R_purelydisc_mart} 
Basic examples of $(\mathcal F_t)$-orthogonal processes are purely discontinuous $(\mathcal F_t)$-local martingales. 
Indeed, according to Theorem 7.34 in \cite{chineseBook} and the comments above, 
any  $(\mathcal F_t)$-local martingale, null at zero, is a purely discontinuous local martingale if and only if it is  $(\mathcal F_t)$-orthogonal. 
\end{remark} 
 
\begin{proposition}\label{P1} 
	If $M$ is a purely discontinuous $(\mathcal F_t)$-local martingale, then  
	$$ 
	[M,M]_t= \sum_{s \leq t} (\Delta M_s)^2. 
	$$ 
\end{proposition} 
\proof 
The result  follows from Theorem 5.2, Chapter I, in \cite{JacodBook}, and Proposition \ref{P_Ito_cov}-(i). 
\endproof 
 
\begin{definition}\label{D_dirichlet} 
	We say that an $(\mathcal{F}_t)$-adapted process $X$ 
 is a \textbf{Dirichlet process} 
	if it admits a decomposition $X=M+A$, where $M$ is a local martingale and $A$ is a finite quadratic variation process with $[A,A]=0$. 
\end{definition} 
\begin{definition}\label{D_weak_dirichlet} 
	We say that $X$ is an $(\mathcal{F}_t)$-adapted 
 \textbf{weak Dirichlet process} 
	if it admits a decomposition $X=M+A$, where $M$ is a local martingale and the process $A$ is $(\mathcal F_t)$-orthogonal. 
\end{definition} 
\begin{definition}\label{D_special_weak_dirichlet} 
	We say that   an $(\mathcal{F}_t)$-adapted process $X$ is  a 
	\textbf{special weak Dirichlet process} 
	if it admits a decomposition of the type above 
such that, in addition, 
$A$ is  predictable. 
\end{definition} 
 
\begin{proposition}\label{P_uniq_spec_weak_Dir_decomp} 
A  Dirichlet process $X$ given in Definition \ref{D_dirichlet} is 
a special weak Dirichlet process, and $[X,X]=[M,M]$. 
\end{proposition} 
\proof
Let $N$ be a continuous local martingale. By Proposition \ref{P_2.15}, $[N,A]=0$. This shows in particular that $X$ is a weak Dirichlet process. It is even special since $A$ is continuous, see Lemma \ref{L_bracket_jumps}.
 By the bilinearity of the covariation, $[X,X]=[M,M]$.  \endproof

\begin{corollary} 
	Let $X$ be an $(\mathcal F_t)$-Dirichlet process. Then	\begin{itemize} 
		\item[(i)] $[X,X]= [M^c,M^c] + \sum_{s \leq \cdot} (\Delta X_s)^2$; 
		\item[(ii)]$[X,X]^c= [M^c,M^c]$. 
	\end{itemize} 
\end{corollary} 
\proof 
(ii) follows by (i) and  Lemma \ref{L_bracket_jumps}-(ii). Concerning (i), by the bilinearity of the covariation, and by the definitions of purely discontinuous local martingale (see Remark \ref{R_purelydisc_mart}) and of  Dirichlet process, we have  
\begin{align*} 
[X,X]_t &= [M^c,M^c]_t + [M^d,M^d]_t\\ 
& = [M^c,M^c]_t + \sum_{s \leq t} (\Delta M_s^d)^2\\ 
&=[M^c,M^c]_t + \sum_{s \leq t} (\Delta X_s)^2, 
\end{align*} 
where the second equality holds because of Proposition \ref{P1}. 
\endproof

\begin{proposition}\label{P_unique_decomp} 
Let $X$ be a special weak Dirichlet process of the type 
\begin{equation}\label{Mc+Md+A} 
X=M^c+M^d+A, 
\end{equation} 
where $M^c$ is a continuous local martingale, and $M^d$ is a purely discontinuous local martingale. 
Supposing that $A_0=M^d_0=0$, the decomposition \eqref{Mc+Md+A}  is unique. 
In that case the decomposition $X = M^c + M^d + A$ will be called the {\bf canonical decomposition} of $X$. 
\end{proposition} 
%
\proof 
Assume that we have 
two  decompositions $X=M^c + M^d +A=(M')^{c}+(M')^{d}\normalcolor + A'$, with $A$ and $A'$ predictable, 
 verifying $[A,N]=[A',N]=0$  for every continuous local martingale $N$. 
 We set $\tilde{A}=A - A'$, $\tilde{M^c}=M^c - (M')^{c}$ and $\tilde{M^d}=M^d - (M')^{d}$. 
By linearity,  
\begin{equation} \label{TildeMart}
\tilde M^c + \tilde M^d + \tilde A=0.
\end{equation} 
We have 
 \begin{align*} 
 0&=[ \tilde M^c + \tilde M^d +\tilde A, \tilde M^c ]  
 = [\tilde M^c, \tilde M^c ] + [\tilde  M^d, \tilde M^c ]+ [ \tilde A, \tilde M^c]\\ 
 &= [\tilde M^c, \tilde M^c], 
 \end{align*} 
 therefore $\tilde M^c=0$ since $\tilde M^c$ is a continuous martingale 
vanishing at zero. 
 It follows in particular that $\tilde{A}$ is a predictable local martingale, hence a  continuous local martingale, see e.g., Corollary 2.24 and Corollary 2.31 in \cite{JacodBook}. 
In particular, calculating the bracket of both sides of
\eqref{TildeMart} against $\tilde M^d$,
	\begin{displaymath} 
	0=[\tilde{M}^d,\tilde{M}^d] + [\tilde{A},\tilde{M}^d]= [\tilde{M}^d,\tilde{M}^d] 
	\end{displaymath} 
	and, since $\tilde{M}^d_0 = 0$, we deduce that $\tilde{M}^d=0$ and therefore $\tilde{A}=0$. 
\endproof 
 
\begin{remark} 
Every $(\mathcal F_t)$-special weak Dirichlet process is of the type \eqref{Mc+Md+A}. 
Indeed, every local martingale $M$ can be decomposed as the sum of a continuous local martingale $M^c$ and a purely discontinuous local martingale $M^d$, see Theorem  4.18, Chapter I, in \cite{JacodBook}. 
\end{remark}

 A very simple example of (special) $(\mathcal F_t)$-weak Dirichlet process is given below. This result extends Remark 2.4-(ii) in \cite{cjms}, valid if $f$ is continuous. 
 \begin{proposition}\label{P_new}
 	Let $f: [0,\,T] \rightarrow \R$ be a c\`adl\`ag (deterministic) function. Then it is $(\mathcal F_t)$-orthogonal; consequently it is a (special) $(\mathcal F_t)$-weak Dirichlet process.
 \end{proposition}
 \proof 
 Let $N$ be an $(\mathcal F_t)$-continuous local martingale.
 We need to show that 
 $$
 [f, N]^{ucp}_{\varepsilon}(t) = 
 \int_0^t \frac{ds}{\varepsilon} \,(f(
 (s+ \varepsilon) \wedge t)- f(s))\,(N(
 (s+ \varepsilon) \wedge t)- N(s))
 $$
 converges u.c.p. to zero. Previous expression is the difference 
 $$
 I^{+ucp}(\varepsilon, t, f, dN)- I^{-ucp}(\varepsilon, t, f, dN),
 $$
 where  $I^{-ucp}(\varepsilon, t, f, dN)$ is given by \eqref{Appr-ucp-int}, and 
 $$
 I^{+ucp}(\varepsilon, t, f, dN):= \int_0^t \frac{ds}{\varepsilon} \,f(
 (s+ \varepsilon) \wedge t)\,(N(
 (s+ \varepsilon) \wedge t)- N(s)).
 $$
By Proposition \ref{P_Ito_cov}-(ii), $I^{-ucp}(\varepsilon, t, f, dN)$ converges u.c.p. to the It\^o-Wiener integral $\int_{]0,\,t]}f(s-)\,d N_s$.
It remains to prove that $I^{+ucp}(\varepsilon, t, f, dN)$ converges u.c.p. to the same quantity. We have 
\begin{align*}
I^{+ucp}(\varepsilon, t, f, dN)	&=  \int_0^t \frac{ds}{\varepsilon} \,f(
 (s+ \varepsilon) \wedge t)\,\int_s^{(s +\varepsilon)\wedge t} d N(u)\\
 &=  \int_0^t \frac{ds}{\varepsilon} \,f(s+ \varepsilon)\,\int_s^{(s +\varepsilon)\wedge t} d N(u) + J(\varepsilon, t), 
\end{align*}
where 
\begin{align*}
|J(\varepsilon, t)| &:=\left| \int_0^t \frac{ds}{\varepsilon} \,(f((s+ \varepsilon)\wedge t)-f(s+ \varepsilon))\,\int_s^{(s +\varepsilon)\wedge t} d N(u)\right|\\
&\leq  \int_{t- \varepsilon}^t \frac{ds}{\varepsilon} \,|f(t)-f(s+ \varepsilon)|\,\delta(N, \varepsilon),
 \end{align*}
and $\delta(N, \varepsilon)$ is the modulus of continuity of $N$. So 
$$
|J(\varepsilon, t)| \leq 2 \, ||f||_{\infty}\,\delta(N, \varepsilon) \underset{\varepsilon \rightarrow 0}{\longrightarrow} 0, \quad \textup{u.c.p.}
$$
On the other hand, by stochastic Fubini's theorem, 
\begin{equation}\label{S_F_T}
\int_0^t \frac{ds}{\varepsilon} \,f(s+ \varepsilon)\,\int_s^{(s +\varepsilon)\wedge t} d N(u) = 
\int_0^t d N(u)\,\int_{(u-\varepsilon)-}^u \frac{ds}{\varepsilon} \,f(s+ \varepsilon).
\end{equation}
Since 
$$
\int_0^T d [N, N]_u\,\left(\int_{(u-\varepsilon)-}^u \frac{ds}{\varepsilon} \,(f(s+ \varepsilon)-f(u))\right)^2\underset{\varepsilon \rightarrow 0}{\longrightarrow} 0 \quad \textup{in probability},
$$
by  \cite{Karatzas1991Brownian}, Problem 2.27, Chapter 3, \eqref{S_F_T} goes u.c.p. to $\int_0^t d N(u) \, f(u) =\int_0^t d N(u) \, f(u-)$, which concludes the proof.
 \endproof

 \begin{remark}
A weak Dirichlet processes may be of finite quadratic variation or not.
 \begin{enumerate}
 \item	
 Assume that $f$ in  Proposition \ref{P_new} has no quadratic variation. This  constitutes a simple example of weak Dirichlet process without quadratic variation. For instance, when $f$ is a single path of a fractional Brownian motion of Hurst index $H < 1/2$, $f$ is a deterministic continuous function, whose quadratic variation is infinite.
\item
Let  $F$ be a function of class $C^{0,1}$, and $X$   an $(\mathcal F_t)$-semimartingale. By the subsequent Theorem \ref{T_C1_dec_weak_Dir}, and since a semimartingale is a weak Dirichlet process with finite quadratic variation (see Remark \ref{R_partWD_WD}-2.), $F(t, X_t)$ is an $(\mathcal F_t)$-weak Dirichlet process.   When  $X$ is continuous, typical examples arise in the framework of SDEs with time-dependent distributional drift, see  for instance [3].
If $F$ depends roughly on time, then $F(t,X_t)$ is expected not to have finite quadratic variation.
\item Let us consider two independent $(\mathcal F_t)$-Brownian motions  $W$ and $B$. 
By Proposition 2.10 of   [2], $\int_0^t B_{t-s}\,d W_s$ is  an  $(\mathcal F_t)$-martingale orthogonal process, therefore an $(\mathcal F_t)$-weak Dirichlet process. By Remark 2.16-2) in [2], it is also a finite quadratic variation process, and $[X,X] = t^2/2$. Similar examples appear in Section 4 of \cite{cjms}.
 \item
At Section 5.3 we have introduced the notion of $(\mathcal F_t)$-particular weak Dirichlet, see Definition \ref{R_PWD_cov}. 
If $A'$ is of zero quadratic variation, then $X$ is of finite quadratic variation.
By Proposition \ref{Sec_paticular},  $[X,X]=[M,M] + \sum_{s \leq \cdot} (\Delta V_s)^2+ 2\,\sum_{s \leq \cdot} \Delta V_s\,\Delta M_s$.
 For instance, this happens when $A'$ is a fractional Brownian motion with Hurst index $H > 1/2$,   independent of $M$.
\end{enumerate}
 \end{remark}

 \begin{proposition}\label{P_TODO_stopped_SWD}
Let $X$ be an $(\mathcal F_t)$-special weak Dirichlet, and $\tau$ an  $(\mathcal F_t)$-stopping time. Then $X^{\tau}=M^{\tau} + A^{\tau}$ is still an $(\mathcal F_t)$-special weak Dirichlet process.
 \end{proposition}
 \proof
Obviously $M^\tau$ is an $(\mathcal F_t)$-local martingale. On the other hand, $A^\tau$ is an $(\mathcal F_t)$-predictable process, see Proposition 2.4-b), Chapter I, in \cite{JacodBook}. Let now $N$ be a continuous $(\mathcal F_t)$-local martingale. It remains to prove that $[A^\tau, N]=0$. By Proposition \ref{P_conv_as}, it will be enough to show that 
$$
C_{\varepsilon}(A^\tau, N)\underset{\varepsilon \rightarrow 0}{\longrightarrow} 0 \quad \textup{u.c.p.}
$$
We have
$$
C_{\varepsilon}(A^\tau, N) := \frac{
1}{\varepsilon} \int_0^t (N_{s+\varepsilon}-N_s)\,(A_{(s+\varepsilon)\wedge \tau}-A_s)\,ds = C_\varepsilon(A, N) + J(\varepsilon, t), 
$$
where 
$$
|J(\varepsilon, t)| \leq 2\,\delta(N, \varepsilon)\,||A||_{\infty}\,T \underset{\varepsilon \rightarrow 0}{\longrightarrow} 0 \quad \textup{u.c.p.}
$$
Since $[A,N]=0$,  again by Proposition \ref{P_conv_as},  
$C_{\varepsilon}(A, N)\underset{\varepsilon \rightarrow 0}{\longrightarrow} 0$ u.c.p.,
which concludes the proof.  
 \endproof
 
\begin{proposition}\label{P_SemiMtgSpecialWD} 
	Let $S$ be an $(\mathcal F_t)$-semimartingale which is a special weak Dirichlet process. Then it is a special semimartingale. 
\end{proposition} 
\proof 
Let $S = M^1+ V$ such that $M^1$ is a local martingale and $V$ is a bounded variation process. Let moreover $S= M^2 + A$, where a predictable 
$(\mathcal F_t)$-orthogonal process. 
Then $0=V-A+M$, where $M=M^2-M^1$. So $A$ is a predictable semimartingale. By Corollary 8.7  in \cite{chineseBook}, 
 $A$ is a special semimartingale, and so by additivity $S$ is a special semimartingale as well. 
\endproof

\subsection{Stability of weak Dirichlet processes under $C^{0,1}$ transformation}

\begin{theorem}\label{T_C1_dec_weak_Dir} 
Let $X = M + A$ be a c\`adl\`ag weak Dirichlet process of finite quadratic variation, 
 and $F:[0,T] \times \R\rightarrow \R$.
\begin{itemize}
\item[(i)]If $F$ is  of class $C^{0,1}$, we have  
\begin{align}\label{CO1_Ito_formula_weak_D} 
& F(t,X_t)=F(0,X_0)+\int_0^t\partial_x F(s,X_{s-})\,d M_s\nonumber\\ 
&+ \int_{]0,\,t]\times \R} (F(s,X_{s-}+x)-F(s,X_{s-}))\,\one_{\{|x| \leq 1 \}}\,(\mu^X-\nu^X)(ds\,dx)\nonumber\\ 
&- \int_{]0,\,t]\times \R} x\,\partial_x F(s,X_{s-})\,\one_{\{|x| \leq 1 \}}\,(\mu^X-\nu^X)(ds\,dx)\nonumber\\ 
&+ \int_{]0,\,t]\times \R} (F(s,X_{s-}+x)-F(s,X_{s-})-x\,\partial_x F(s,X_{s-}))\,\one_{\{|x| > 1 \}}\,\mu^X(ds\,dx)+\Gamma^F(t), 
\end{align} 
where $\Gamma^F: C^{0,1}\rightarrow \mathbb{D}^{ucp}$ is a continuous linear map
such that, for every $F \in C^{0,1}$, it fulfills the following  properties. 
\begin{itemize} 
\item[(a)] $[\Gamma^F, N] = 0$  for every $N$ continuous local martingale. 
\item[(b)] If $A$ is predictable, then $\Gamma^F$ is  predictable. 
\end{itemize} 
\item[(ii)]
If $F$ is  of class $C^{1,2}$, \eqref{CO1_Ito_formula_weak_D} holds with
\begin{align} 
\Gamma^F(t)&:=\int_0^t\partial_s F(s,X_s)\,ds + \int_0^t\partial_x F(s,X_s)\,d^{-}A_s +\int_0^t\partial_{xx}^2 F(s,X_s)\,d[X,X]_s^c\nonumber\\ 
&  + \int_{]0,\,t]\times \R} (F(s,X_{s-}+x)-F(s,X_{s-})-x\,\partial_x F(s,X_{s-}))\one_{\{|x| \leq 1 \}}\,\nu^X(ds\,dx).\label{Gamma_C12_weakD_BIS} 
\end{align}
\end{itemize}
\end{theorem} 

 \begin{remark}
 	Point (a) in Theorem \ref{T_C1_dec_weak_Dir}-(i) implies in  particular that $F(s,X_s)$ is a weak Dirichlet process when $X$ is a weak Dirichlet  process. 
 \end{remark}

\begin{remark}\label{R_gamma_pred} 
Proposition \ref{ucp_predic}  implies that, when $A$ is predictable,  
	$\Gamma^F$ in \eqref{Gamma_C12_weakD_BIS} is a predictable process for any $F \in C^{1,2}$. 
\end{remark}

\proof
We start by proving  item (ii).
Expressions 
\eqref{CO1_Ito_formula_weak_D}-\eqref{Gamma_C12_weakD_BIS} follow 
by Proposition \ref{C_reformulation_C12_Ito}, in particular 
 by \eqref{Ito_formula_C2_jump_measure}. 
We remark  that, since $M$ is a local martingale and 
 $\partial_x F(s, X_s)$ is a c\`adl\`ag process,  by Proposition \ref{P_Ito_cov}-(ii) 
 we have 
\begin{eqnarray*} 
\int_0^t\partial_x F(s,X_s)\,d^{-} X_s  &=& \int_0^t\partial_x F(s,X_s)\,d^{-} M_s  + \int_0^t\partial_x F(s,X_s)\,d^{-} A_s\\ 
 &=& \int_0^t\partial_x F(s,X_{s-})\,d M_s  + \int_0^t\partial_x F(s,X_s)\,d^{-} A_s. 
\end{eqnarray*} 
Let us now prove item (i). Let $F$ be a function of class $C^{0,1}$. Setting
	\begin{align*}  
	K_s(x)&=\left(F(s,X_{s-} + x)-F(s,X_{s-})\right)\one_{\{\vert x\vert \leq 1\}},\\ 
 	Y_s(x)&=x\,\partial_x F(s,X_{s-})\one_{\{\vert x \vert \leq 1\}},
	\end{align*}
Proposition \ref{P_X_jump_comment2} insures that 
   $K \ast (\mu^X-\nu^X)$, $Y \ast (\mu^X-\nu^X)$ are well-defined and are purely discontinuous square integrable 
 local martingales. 
On the other hand, the fourth integral of the right hand-side of \eqref{CO1_Ito_formula_weak_D} is well-defined since, by Remark \ref{R_Jumps_measure}, it equals 
$$
\sum_{s\leq t}\,[F(s,X_s)-F(s,X_{s-})- \partial_x F(s,X_{s-})\,\Delta \,X_{s}]\,\one_{\{|\Delta \,X_{s}| >1\}}.
$$
The sum is finite since there are at most a finite number of  jumps whose size is larger than one.
Previous considerations, in agreement with \eqref{CO1_Ito_formula_weak_D}, allow us to set 
\begin{align} \label{GammaF} 
	\Gamma^F(t) :=& F(t,X_t)- F(0,X_0)-\int_0^t\partial_x F(s,X_{s-})\,d M_s \\ 
	& -\int_{]0,\,t]\times \R} \left\{F(s,X_{s-}+x)-F(s,X_{s-})-x\,\partial_x F(s,X_{s-})\right\}\,\one_{\{|x| > 1 \}}\,\mu^X(ds\,dx)\nonumber\\ 
	& -\int_{]0,\,t]\times \R} \left\{F(s,X_{s-}+x)-F(s,X_{s-})-x\,\partial_x F(s,X_{s-})\right\}\,\one_{\{|x| \leq 1 \}}\,(\mu^X-\nu^X)(ds\,dx).\nonumber 
\end{align} 
The second step consists in 
proving that $C^{0,1} \ni F \mapsto \Gamma^F(t)$ is continuous with respect to the u.c.p. topology. 
For this we first observe  that the map $F \mapsto F(t,X_t)- F(0,X_0)$ fulfills the mentioned continuity. Moreover, if $F^n \rightarrow F$ in $C^{0,1}$, then $\int_0^t(\partial_x F^n-\partial_x F)(s,X_{s-})\,d M_s$ converges to zero u.c.p. since 
$\partial_x F^n(s,X_{s-})$ converges to $\partial_x F(s,X_{s-})$ in 
 $\mathbb{L}^{ucp}$, see Chapter II Section 4 in \cite{protter}. 
 
Let us consider the second line of \eqref{GammaF}. For almost all 
 fixed $\omega$, 
the process $X$ has  a finite number of jumps, 
$s_i= s_i(\omega), 1 \le i \le N(\omega)$, 
 larger than one. 
Let $F^n\rightarrow F$  in $C^{0,1}$. Since the map is linear we can suppose that $F=0$. 
\begin{align*} 
&\sup_{0 <t \leq T}\bigg |\int_{]0,\,t]\times \R} \left\{F^n(s,X_{s-}(\omega)+x)-F^n(s,X_{s-}(\omega))-x\,\partial_x F^n(s,X_{s-}(\omega))\right\}\,\one_{\{|x| > 1 \}}\,\mu^X(\omega,ds\,dx)\bigg|\\ 
&\leq  \int_{]0,\,T]\times \R} \left\vert 
F^n(s,X_{s-}(\omega)+x)-F^n(s,X_{s-}(\omega))-x\,\partial_x F^n(s,X_{s-}(\omega))\right\vert \,\one_{\{|x| > 1 \}}\,\mu^X(\omega,ds\,dx) \\ 
& = \sum_{i=1}^{N(\omega)} \left\vert F^n(s_i,X_{s_i}(\omega))-F^n({s_i},X_{{{s_i}-}}(\omega))-\Delta X_{{{s_i}}}(\omega)\,\partial_x F^n(s_i,X_{{s_i}-}(\omega)) 
\right\vert\,\one_{\{|\Delta X_{{{s_i}}}(\omega)| > 1 \}}  
\underset{n \rightarrow \infty}{\rightarrow} 0. 
\end{align*} 
This  shows in particular that 
\[\int_{]0,\,\cdot]\times \R} \{F^n(s,X_{s-}(\omega)+x)-F^n(s,X_{s-}(\omega))-x\,\partial_x F^n(s,X_{s-}(\omega))\}\,\one_{\{|x| > 1 \}}\,\mu^X(\omega,ds\,dx) \underset{n \rightarrow \infty}{\rightarrow} 0\quad \textup{u.c.p.} 
\] 
and so the map defined by the second line in \eqref{GammaF} is continuous.

The following proposition shows the continuity properties of the  last term in \eqref{GammaF} with respect to $F$, which finally  allows to conclude the continuity of the map $ \Gamma^F: C^{0,1} \rightarrow \mathbb D^{ucp}$.
\begin{proposition}\label{P_cont_ucp_int_stoc} 
	The map 
	\begin{eqnarray*} 
	I: C^{0,1} &\rightarrow& \mathbb{D}^{ucp} \\ 
		g &\mapsto& \int_{]0,\cdot]\times \R}G^g\,(s,\,X_{s-},\,x)\,\one_{\{|x| \leq 1\}}\,(\mu^X-\nu^X)(ds\,dx), 
	\end{eqnarray*} 
	where 
	\begin{equation}\label{Ggn} 
	G^g\,(s,\,\xi,\,x)=g(s, \,\xi + x)-g(s,\xi)-x\,\partial_\xi g(s,\,\xi), 
	\end{equation} 
	is continuous. 
\end{proposition} 
\proof 
We consider  the  sequence  $(\tau_{l})_{l \geq 1}$ of increasing stopping times introduced in Remark \ref{R_localisation}-(ii) for the  process $Y_t = (X_{t-},  \sum_{s < t}|\Delta X_{s}|^2)$. 
Since $\Omega = \cup_l\,\{\omega: \tau_l(\omega)>T\}$ a.s.,  
the result is proved if we show that, for every fixed 
$\tau = \tau_l$, 
\begin{displaymath} 
g \mapsto 
\one_{\{ \tau>T\}} (\omega) 
\int_{]0,\,\cdot] \times \R}G^g(s,\,X_{s-},\,x)\,\one_{\{|x| \leq 1\}}\,(\mu^X-\nu^X)(ds\,dx) 
\end{displaymath} 
is  continuous. 
Let $g^n\rightarrow g$  in $C^{0,1}$. Then $G^{g^n}\rightarrow G^g$ in $C^{0}([0,\,T]\times \R^2)$. Since the map is linear we can suppose that $g=0$. 
Let $\varepsilon_0 >0$.  
We aim at showing that 
\begin{equation}\label{prob_conv} 
\P\bigg(\sup_{t \in [0,\,T]} 
\bigg| 
\one_{\{ \tau>T\}} (\omega) 
\int_{]0,\,t]\times \R}G^{g^n}(s,\,X_{s-},\,x)\,\one_{\{|x| \leq 1\}}\,(\mu^X-\nu^X)(ds\,dx) 
\bigg| > \varepsilon_0 
\bigg) \underset{n \rightarrow \infty}{\longrightarrow} 0. 
\end{equation} 
Let $W^{n}_s(x)$ (resp. by $\hat{W}^{n}_s$) denote the random field 
 $G^{g^n}(s,\,X_{s-},\,x)\,\one_{\{|x| \leq 1\}}$ (resp. the process 
$\int_{\R}G^{g^n}(s,\,X_{s-},\,x)\,\one_{\{|x| \leq 1\}}\,\hat{\nu}^X(dx)$), and define 
\begin{align*} 
& I^n_t := 
\int_{]0,\,t]\times \R}W^{n}_s(x)\,(\mu^X-\nu^X)(ds\,dx). 
\end{align*} 
The convergence \eqref{prob_conv} will follow if we show that 
\begin{equation}\label{prob_conv2} 
\P\big(\sup_{t \in [0,\,T]} 
|I^n_{t \wedge \tau}| > \varepsilon_0 
\big) \underset{n \rightarrow \infty}{\longrightarrow} 0. 
\end{equation} 
For every process $\phi=(\phi_t)$, we indicate the stopped process at $\tau$ by 
$ 
\phi_{t}^{\tau} (\omega) := \phi_{ t \wedge \tau (\omega)} (\omega). 
$ 
We prove below that 
\begin{equation}\label{W2_A_loc} 
(|W^n|^2 \ast \mu^X)^{\tau} \in \mathcal{A}^+. 
\end{equation} 
As a matter of fact, let $M$ such that $\sup_{t \in [0,\,T]}|Y_{t \wedge \tau}\,\one_{\{\tau >0 \}}|\leq M$. 
Recalling Remark \ref{R_Jumps_measure}, an obvious Taylor expansion yields 
\begin{align}\label{estimate_G2} 
&\sper{\int_{]0,\,t \wedge \tau]\times \R}|W^{n}_s(x)|^2\,\mu^X(ds,\,dx)}\nonumber\\ 
&\leq 2 \sup_{ \underset{t \in [0,\,T]}{y \in  [-M,\,M]}}|\partial_x g^n|^2(t,\,y)\, \sper{\sum_{0 <s < {\tau}}|\Delta X_s|^2 \,\one_{\{|\Delta X_s| \leq 1\}} \,\one_{\{\tau >0\}}+|\Delta X_{\tau}|^2 \,\one_{\{|\Delta X_{\tau}| \leq 1\}}\,\one_{\{\tau >0\}}}\nonumber\\ 
&\leq 2\sup_{ \underset{t \in [0,\,T]}{y \in  [-M,\,M]}}|\partial_x g^n|^2(t,\,y)\, \cdot (M+1),  
\end{align} 
which shows \eqref{W2_A_loc}.
By Lemma 2.4-1. in \cite{BandiniRusso2},
it follows 
that $W^n\,\one_{[0,\,\tau]} \in \mathcal{G}^2(\mu^X)$. 
Consequently,  by Proposition 3.66 of \cite{jacod_book},  
\begin{equation}\label{EA2} 
I_{t \wedge \tau}^n \textup{\,\,is a purely discontinuous square integrable martingale}. 
\end{equation} 
On the other hand, $W^n \in \mathcal{G}^2_{\rm loc}(\mu^X)$, 
and by Theorem 11.12, point 3), in \cite{chineseBook} 
	,   it follows that 
\begin{equation}\label{predictable_bracket} 
\langle I^n,\,I^n \rangle_{t} 
= \int_{]0,\,t]\times \R}|W^{n}_s(x)|^2\,\nu^X(ds\,dx)- 
\sum_{0 <s \leq t}\, |\hat{W}^{n}_s|^2\leq \int_{]0,\,t]\times \R}|W^{n}_s(x)|^2\,\nu^X(ds\,dx). 
\end{equation} 
Taking into account \eqref{EA2}, we can apply  Doob's inequality.   
Using estimates \eqref{estimate_G2}, \eqref{predictable_bracket} and \eqref{EA2}, we get 
\begin{align*} 
	\P\left[\sup_{t \in [0,\,T]}|I_{t \wedge \tau}^n|> \varepsilon_0\right] 
	&\leq \frac{1}{\varepsilon_0^2}\,\sper{|I_{T \wedge \tau}^n|^2}\\
	&=\frac{1}{\varepsilon_0^2}\,\sper{\langle I^n,\,I^n\rangle_{T \wedge \tau}}\\ 
	&\leq \frac{2\,(M+1)}{\varepsilon_0^2} \sup_{ \underset{t \in [0,\,T]}{y \in  [-M,\,M]}}|\partial_x g^n|^2(t,\,y). 
\end{align*} 
Therefore, since $\partial_x g^n \rightarrow 0$ in $C^0$ as $n$ goes to infinity, 
\[ 
\lim_{n \rightarrow \infty}	\P\left[\sup_{t \in [0,\,T]}|I_{t \wedge \tau}^n|> \varepsilon_0\right] =0. 
\] 
\qed 
 
We continue the proof of item (i) of Theorem \ref{T_C1_dec_weak_Dir}. 
It remains to prove items (a) and (b). 
 
(a) 
We have to prove that, for any continuous local martingale $N$, we have 
\begin{align*} 
&\bigg[F(\cdot, X)-\int_0^{\cdot} \partial_x F(s,X_{s-})\,dM_s \\ 
& -\int_{]0,\,\cdot]\times \R} \left\{F(s,X_{s-}+x)-F(s,X_{s-})-x\,\partial_x F(s,X_{s-})\right\}\,\one_{\{|x| > 1 \}}\,\mu^X(ds\,dx)\nonumber\\ 
& -\int_{]0,\,\cdot]\times \R} \left\{F(s,X_{s-}+x)-F(s,X_{s-})-x\,\partial_x F(s,X_{s-})\right\}\,\one_{\{|x| \leq 1 \}}\,(\mu^X-\nu^X)(ds\,dx), N\bigg]=0. 
\end{align*} 
We set 
\begin{align*} 
Y_t&=\int_{]0,\,t]\times \R}  W_s(x)\,\one_{\{|x| \leq 1 \}}\,(\mu^X-\nu^X)(ds\,dx),\\ 
Z_t &=\int_{]0,\,t]\times \R} W_s(x)\,\one_{\{|x| > 1 \}}\,\mu^X(ds\,dx), 
\end{align*} 
with 
$
W_s(x) = F(s,X_{s-}+x)-F(s,X_{s-})-x\,\partial_x F(s,X_{s-})$. 

Since $X$ has almost surely a finite number of jumps larger than one, by Remark \ref{R_Jumps_measure}  $Z$ is a bounded variation process.
Being $N$  continuous, Proposition \ref{P_quad_var} insures that 
$\left[Z, N\right] = 0$.
By Proposition  \ref{P_X_jump_comment2} $W^2 \one_{\{|x| \leq 1 \}} \ast \mu^X \in {\mathcal A}^+_{\rm loc}$, 
therefore $W \one_{\{|x| \leq 1 \}}$ belongs to 
${\mathcal G}^2_{\rm loc}(\mu^X)$ as well, see  Lemma 2.4-2. in \cite{BandiniRusso2}. 
In particular, by Theorem 11.21, point 3), in \cite{chineseBook}
, $Y$  is a purely discontinuous (square integrable) 
  local martingale. 
Recalling that a local  $(\mathcal F_t)$-martingale, null at zero, is a purely discontinuous martingale if and only if it is  $(\mathcal F_t)$-orthogonal (see Remark \ref{R_purelydisc_mart}), 
 from Proposition \ref{P_Ito_cov}-(i) we have 
$
\left[Y, N\right] = 0
$.
From Proposition \ref{P_Ito_cov}-(iii), and the fact that $[M,N]$ is continuous, it follows that 
\begin{eqnarray*} 
\left[\int_0^{\cdot} \partial_x F(s,X_{s-})\,dM_s, N\right] = \int_0^{\cdot} \partial_x F(s,X_{s-})\,d \left[M, N\right]_s. 
\end{eqnarray*} 
Therefore it remains to check that 
\begin{equation}\label{identity} 
[F(\cdot, X),N]_t = \int_0^{\cdot} \partial_x F(s,X_{s-})\,d \left[M, N\right]_s. 
\end{equation} 
To this end, we  evaluate the limit  of 
\begin{align*} 
&\frac{1}{\varepsilon}\int_0^t (F((s+\varepsilon)\wedge t, X_{(s+\varepsilon)\wedge t})-F(s, X_s))\, (N_{(s+\varepsilon)\wedge t}- N_s)\,ds\\ 
&=\frac{1}{\varepsilon}\int_0^t (F((s+\varepsilon)\wedge t, X_{(s+\varepsilon)\wedge t})-F((s+\varepsilon)\wedge t, X_s))\, (N_{(s+\varepsilon)\wedge t}- N_s)\,ds\\ 
&\quad + \frac{1}{\varepsilon}\int_0^t (F((s+\varepsilon)\wedge t, X_s)-F(s, X_s))\, (N_{(s+\varepsilon)\wedge t}- N_s)\,ds\\ 
&=: I_1(\varepsilon, \,t)+  I_2(\varepsilon, \,t). 
\end{align*} 
Concerning the term $I_1(\varepsilon, \,t)$, it can be decomposed as 
$$ 
I_1(\varepsilon, \,t)= I_{11}(\varepsilon,\,t) +I_{12}(\varepsilon,\,t) +I_{13}(\varepsilon,\,t), 
$$ 
where 
\begin{align*} 
I_{11}(\varepsilon,\,t) &=\frac{1}{\varepsilon}\int_0^t \partial_x F(s, X_{s})\, (N_{(s+\varepsilon)\wedge t}- N_s)(X_{(s+\varepsilon)\wedge t}-X_{s})\,ds,\\ 
I_{12}(\varepsilon,\,t)&=\frac{1}{\varepsilon}\int_0^t (\partial_x F((s+\varepsilon)\wedge t, X_{s})- \partial_x F(s, X_{s}))\, (N_{(s+\varepsilon)\wedge t}- N_s)(X_{(s+\varepsilon)\wedge t}-X_{s})\,\,ds,\\ 
I_{13}(\varepsilon,\,t)&= \frac{1}{\varepsilon}\int_0^t \left(\int_0^1 (\partial_x F((s+\varepsilon)\wedge t, X_{s} + a (X_{(s+\varepsilon)\wedge t}-X_{s}))- \partial_x F((s+\varepsilon) \wedge t, X_{s}))\,da\right)\cdot\nonumber\\ 
&\cdot (N_{(s+\varepsilon)\wedge t}- N_s)(X_{(s+\varepsilon)\wedge t}-X_{s})\,ds. 
\end{align*} 
Notice that the brackets $[X,X]$, $[X,N]$ and $[N,N]$  exist. 
Indeed, $[X,X]$  exists by definition, $[N,N]$  exists by Proposition \ref{P_Ito_cov}-(i). Concerning 
  $[X,N]$, it 
can be decomposed as 
$
[X,N]=[M,N]+ [A,N]$,
where $[M,N]$ exists by Proposition \ref{P_Ito_cov}-(i) and  $[A,N]=0$   by hypothesis, since $A$ comes 
from the weak Dirichlet decomposition of $X$. 
 
Then,  from  Corollary \ref{C_id_RV_fw}-2) 
and Proposition \ref{P_equiv_bracket} we have 
\begin{equation}\label{I11_conv} 
I_{11}(\varepsilon,\,t) \underset{\varepsilon \rightarrow 0}{\longrightarrow}\, \int_0^t \partial_x F(s, X_{s-})\, d[M,N]_s\quad \textup{u.c.p.} 
\end{equation} 
At this point, we have  to prove the u.c.p. convergence to zero of  the remaining terms  $I_{12}(\varepsilon,t)$, $I_{13}(\varepsilon,t)$, $I_{2}(\varepsilon,t)$. 
First,  since $\partial_x F$ is uniformly continuous on each compact, we have 
\begin{align}\label{I_12_ineq} 
\vert I_{12}(\varepsilon,\,t)\vert \leq \delta\Big( \partial_x F\bigg|_{[0,\,T]\times \mathbb{K}^X}  ;\,\varepsilon\Big)\, 
 \sqrt{   [X,X]_{\varepsilon}^{ucp} [N,N]_{\varepsilon}^{ucp}}, 
\end{align} 
where $\mathbb{K}^X$ is the (compact) set $\{X_t(\omega) : t \in [0,\,T]\}$. 
When $\varepsilon$ goes to zero, the 
modulus of continuity component in \eqref{I_12_ineq}  converges  to zero a.s., 
while the remaining term u.c.p. converges to $ \sqrt{[X,X]_t [N,N]_t}$ by definition. 
Therefore, 
\begin{equation}\label{I12_conv} 
I_{12}(\varepsilon,\,t) \underset{\varepsilon \rightarrow 0}{\longrightarrow}\, 0\quad \textup{u.c.p.} 
\end{equation}

Let us then evaluate $I_{13}(t,\,\varepsilon)$. 
Since  $[X,X]^{ucp}_\varepsilon$, $[N,N]^{ucp}_\varepsilon$ u.c.p. converge, there exists of a sequence $(\varepsilon_n)$ such that 
$[X,X]^{ucp}_{\varepsilon_n}$, $[N,N]^{ucp}_{\varepsilon_n}$ converges uniformly a.s. 
respectively to $[X,X]$, $[N,N]$. 
We fix a realization $\omega$ outside a null set. Let $\gamma >0$. 
We enumerate the jumps of $X(\omega)$  on $[0,\,T]$ by $(t_i)_{i \geq 0}$. 
Let $M = M(\omega)$ such that 
\begin{displaymath} 
\sum_{i = M+1}^{\infty} |\Delta X_{t_i}|^2 \leq \gamma^2. 
\end{displaymath} 
We define $A(\varepsilon_n,M)$ and $B(\varepsilon_n,M)$ as in \eqref{Aepsilon}-\eqref{Bepsilon}.
The term $I_{13}(\varepsilon_n,\,t)$ can be decomposed as the sum of two  terms: 
\begin{align*} 
I^A_{13}(\varepsilon_n,\,t)&=\sum_{i = 1}^M \int_{t_i - \varepsilon_n}^{t_i}\frac{ds} 
{\varepsilon_n}\, 
\one_{]0,\,t]}(s)\,(X_{(s+ \varepsilon_n)\wedge t}-X_s)(N_{(s+ \varepsilon_n)\wedge t}-N_s)\cdot \\ 
&\quad \cdot\int_0^1 (\partial_x F((s+ \varepsilon_n)\wedge t,\,X_s + a (X_{(s+ \varepsilon_n)\wedge t}-X_s))-\partial_x F((s+ \varepsilon_n)\wedge t,\,X_s))\,da,\nonumber\\ 
I^B_{13}(\varepsilon_n,\,t)&= \frac{1} 
{\varepsilon_n}\int_{]0,\,t]} 
(X_{(s+ \varepsilon_n)\wedge t}-X_s)(N_{(s+ \varepsilon_n)\wedge t}-N_s)\,R^B(\varepsilon_n,s,t,M)\,ds, 
\end{align*} 
with 
$$ 
R^B(\varepsilon_n,s,t, M)= \one_{B(\varepsilon_n,M)}(s)\,\int_0^1 [\partial_x F((s+ \varepsilon_n)\wedge t,\,X_s + a (X_{(s+ \varepsilon_n)\wedge t}-X_s))-\partial_x F((s+ \varepsilon_n)\wedge t,\,X_s)]\,da. 
$$ 
By Remark \ref{R_Billingsley}, we have for every $s$, $t$, 
\begin{align*} 
R^B(\varepsilon_n,s,t,M) &\leq \delta\bigg( \partial_x F \bigg|_{[0,\,T]\times \mathbb{K}^X} 
,\,\sup_{l} \sup_{\underset{|r-a| \leq \varepsilon_n}{r, a \in [t_{l-1},\,t_{l}]}} |X_{a}-X_{r}|\bigg), 
\end{align*} 
so that Lemma \ref{Lem_Billingsley} applied successively to the intervals $[t_{l-1},\,t_{l}]$ implies 
\begin{align*} 
R^B(\varepsilon_n,s,t,M) &\leq \delta\big( \partial_x F \big|_{[0,\,T]\times \mathbb{K}^X} 
,3\gamma\big). 
\end{align*} 
Then 
$$ 
|I_{13}^B(\varepsilon_n,\,t)| \leq \delta\big( \partial_x F \big|_{[0,\,T]\times \mathbb{K}^X} 
,3\gamma\big) \sqrt{ [N,N]_{\varepsilon_n}^{ucp}(T)\,[X,X]_{\varepsilon_n}^{ucp}(T)}, 
$$ 
and we get 
\begin{equation}\label{I13B_limsup} 
\limsup_{n \rightarrow \infty} \sup_{t \in [0,T]}|I_{13}^B(\varepsilon_n,t)|\leq \delta\big( \partial_x F \big|_{[0,\,T]\times \mathbb{K}^X} 
,3\gamma\big)\sqrt{ [N,N]_T\,[X,X]_T}. 
\end{equation} 
 
Concerning 
$I^A_{13}(\varepsilon_n,\,t)$, 
 we apply Lemma \ref{L_ucp_big_jumps} to $Y = (Y^1,Y^2,Y^3)= (t,X,N)$ and 
\begin{align*} 
\phi(y_1, y_2) 
=(y^2_1-y^2_2)\,(y^3_1-y^3_2)\int_0^1 [\partial_x F(y_1^1,\,y^2_2 + a (y^2_1-y^2_2))-\partial_x F(y_1^1,\,y^2_2)]\,da. 
\end{align*} 
Then $I^A_{13}(\varepsilon_n,t)$ converges uniformly in $t \in [0,\,T]$, as $n$ goes to infinity, to 
\begin{equation}\label{limi_IA13} 
\sum_{i = 1}^M 
\one_{]0,\,t]}(t_i)\,(X_{t_i}-X_{t_{i}-})(N_{t_i}-N_{t_{i}-})\int_0^1 [\partial_x F(t_i,\,X_{t_{i}-} + a (X_{t_i}-X_{t_{i}-}))-\partial_x F(t_i,\,X_{t_{i}-})]\,da. 
\end{equation} 
In particular, \eqref{limi_IA13} equals zero since $N$ is a continuous process. 
Then, recalling \eqref{I13B_limsup}, we have 
$$ 
\limsup_{n \rightarrow \infty}\sup_{t \in [0,\,T]} |I_{13}(\varepsilon_n,\,t)|\leq \delta( \partial_x F 
,3\,\gamma)\sqrt{[N,N]_T\,[X,X]_T}.
$$ 
Letting $\gamma$ go to zero, we conclude that 
\begin{equation}\label{I13_conv} 
\limsup_{n \rightarrow \infty}\sup_{t \in [0,\,T]} |I_{13}(\varepsilon_n,\,t)|=0. 
\end{equation}

It remains to show 
the u.c.p. convergence to zero of 
$I_2(\varepsilon,\,t)$, 
 as $\varepsilon \rightarrow 0$. 
To this end, let us  write it  as the sum of the two terms 
\begin{align*} 
I_{21}(\varepsilon,\,t) &=\frac{1}{\varepsilon}\int_0^t (F(s+\varepsilon, X_s)-F(s, X_s))\, (N_{(s+\varepsilon)\wedge t}- N_s)\,ds,\\ 
I_{22}(\varepsilon,\,t) &=\frac{1}{\varepsilon}\int_0^t (F((s+\varepsilon)\wedge t, X_s)-F(s+ \varepsilon, X_s))\, (N_{(s+\varepsilon)\wedge t}- N_s)\,ds. 
\end{align*} 
Concerning $I_{21}(\varepsilon,\,t)$, it can be written as 
\begin{align}\label{I2} 
I_{21}(\varepsilon,\,t) & 
= \int_{]0,\,t]}J_{\varepsilon}(r)\,d N_r 
\end{align} 
with 
\begin{displaymath} 
J_{\varepsilon}(r) = 
\int_{[(r-\varepsilon)_+, \, r[} \frac{F(s+\varepsilon, X_s)-F(s, X_s)}{\varepsilon}\,\,ds. 
\end{displaymath} 
Since 
$J_{\varepsilon}(r)\rightarrow 0$ pointwise, it follows from the  Lebesgue dominated convergence theorem that 
\begin{equation}\label{karatzas} 
\int_0^T J^2_{\varepsilon}(r) \, d\langle N,\,N \rangle_r \limit^{\P} 0 \qquad \textup{as}\,\, \varepsilon \rightarrow 0. 
\end{equation} 
Therefore, according to \cite{Karatzas1991Brownian}, Problem 2.27 in Chapter 3, 
\begin{equation}\label{I21_conv} 
\lim_{\varepsilon \rightarrow 0} \sup_{t \in [0,T]}|I_{21}(\varepsilon,t)|=0. 
\end{equation} 
 
As far as $I_{22}(\varepsilon,t)$ is concerned, we have 
\begin{align*} 
|I_{22}(\varepsilon,\,t)| &\leq \frac{1}{\varepsilon}\int_{t-\varepsilon}^{t} |F(t, X_s)-F(s+ \varepsilon, X_s)|\, |N_{t}- N_s|\leq 2\,\delta\big(F \big|_{[0,\,T]\times \mathbb{K}^X} 
,\varepsilon\big)\,||N||_{\infty} 
\end{align*} 
and we get 
\begin{equation}\label{I22_conv} 
\limsup_{\varepsilon \rightarrow 0} \sup_{t \in [0,T]}|I_{22}(\varepsilon,t)|=0. 
\end{equation} 
This concludes the proof of item (a). 
 
(b)  Let $F^n$ be a sequence   $C^{1,2}$ functions such that  $F^n \rightarrow F$ 
and $\partial_x F^n \rightarrow \partial_x F$, uniformly on every compact subset. 
From item (ii),
the process $\Gamma^{F^n}(t)$ in \eqref{Gamma_C12_weakD_BIS} 
equals 
\begin{eqnarray*} 
	&&\int_0^t\partial_s F^n(s,X_s)\,ds + \int_0^t\partial_x F^n(s,X_s)\,d^{-}A_s +\int_0^t\partial_{xx}^2 F^n(s,X_s)\,d[X,X]_s^c  \\ 
	&& + \int_{[0,\,t]\times \R} (F^n(s,X_{s-}+x)-F^n(s,X_{s-})-x\,\partial_x F^n(s,X_{s-}))\,\one_{\{|x| \leq 1\}}\nu^X(ds\,dx), 
\end{eqnarray*} 
which is predictable, see Remark \ref{R_gamma_pred}. 
Since, by  point (a),    the map $\Gamma^F: C^{0,1}\rightarrow \mathbb{D}^{ucp}$ is continuous, 
$\Gamma^{F^n}$ converges to $\Gamma^{F}$ u.c.p. 
Then $\Gamma^{F}$ is predictable because it is the u.c.p. limit of predictable processes. 
\endproof 
 
\subsection{A class of particular weak Dirichlet processes}\label{Sec_paticular}
The notion of Dirichlet process is a natural extension of the one of semimartingale only in the continuous case. 
Indeed, if $X$ is a c\`adl\`ag process, which is also Dirichlet, then $X=M+A'$ with $[A',A']=0$, and therefore $A'$ is continuous because of Lemma \ref{L_bracket_jumps}. 
This class does not include all  the c\`adl\`ag semimartingale $S=M+V$, perturbed by a zero quadratic variation process $A'$. 
Indeed, if $V$ is not continuous, $S +A'$  is not necessarily a Dirichlet process,
even though $X$ is a weak Dirichlet process. 
Notice that, in general, it is 
even not a special weak Dirichlet process, since $V$ is generally not predictable.\\ 
 
We propose then the following natural extension of the semimartingale notion in the weak Dirichlet framework. 
\begin{definition}\label{D_PWD} 
	We say that $X$ is an $({\mathcal F}_t)$-\textbf{particular weak 
 Dirichlet process} 
	if it admits a decomposition $X=M+A$, where $M$ is an  $({\mathcal F}_t)$-local martingale,  $A=V+A'$ with 
 $V$ being a bounded variation adapted  process and  $A'$  a continuous 
adapted process  $({\mathcal F}_t)$-orthogonal process  
 such that $A'_0 = 0$. 
\end{definition} 
\begin{remark}\label{R_partWD_WD} 
	\begin{enumerate} 
	\item A particular weak Dirichlet process is a weak Dirichlet process. Indeed 
	by Proposition \ref{P_quad_var} we have   $[V,N] =0$, so 
	$$ 
	[A' + V, N]= [A',N]+[V,N] =0. 
	$$ 
	\item By construction an $(\mathcal F_t)$-semimartingale is a particular $(\mathcal F_t)$-weak Dirichlet process, so also an $(\mathcal F_t)$-weak Dirichlet process.
	\item There exist processes that are (even special) weak Dirichlet and not  particular weak Dirichlet. 
	As a matter of fact, let for instance consider the deterministic 
process $A_t = \one_{\Q \cap [0,\,T]}(t)$. Then $A$ is predictable and $[A,N]=0$ for any $N$ continuous local 
martingale, since, the fact that $A \equiv 0$ Lebesgue a.e. implies that $[A,N]^{ucp}_\varepsilon \equiv 0$. 
Moreover, since $A$ is totally discontinuous, it can not have bounded variation, so that $A$ is special 
weak Dirichlet but not a  particular weak Dirichlet process. 
	\end{enumerate} 
\end{remark}

The result below is an extension of Proposition \ref{P_uniq_spec_weak_Dir_decomp}.
\begin{proposition}\label{R_PWD_cov}
Let $X=M+V+A'$ be an-$(\mathcal F_t)$-particular weak Dirichlet process, with $A'$ of zero quadratic variation. Then $X$ is of finite quadratic variation, and 
$$
[X,X]=[M,M] + \sum_{s \leq \cdot} (\Delta V_s)^2+ 2\,\sum_{s \leq \cdot} \Delta V_s\,\Delta M_s.
$$
\end{proposition}
\proof
By Proposition \ref{P_quad_var}, $[V,V]_t=\sum_{s \leq t} (\Delta V_s)^2$ and $[M,V]_t=\sum_{s \leq t} \Delta V_s\,\Delta M_s$.
By Proposition \ref{P_2.15}, $[M, A']+ [V, A']\equiv 0$. The result follows by the bilinearity of the covariation.
\endproof

In Propositions \ref{P_int_dec_par_WD} and \ref{P_char_part_weak_D}  
we extend some properties valid for semimartingales to the case of particular weak Dirichlet processes. 

\begin{proposition}\label{P_int_dec_par_WD} 
Let $X$ be an  $(\mathcal F_t)$-adapted c\`adl\`ag process satisfying  assumption \eqref{StandAss}. $X$ is a particular weak Dirichlet process if and only if there exist a continuous local martingale $M^c$, a predictable process $\alpha$ of the type $\alpha^S+ A'$, where $\alpha^S$ is predictable with bounded variation, $A'$ is a  $({\mathcal F}_t)$-adapted continuous orthogonal process, $\alpha_0^S= A'_0=0$, 
and 
\begin{equation}\label{integral_repr} 
	X =  M^c + \alpha + (x\,\one_{\{|x| \leq 1\}})\ast (\mu^X - \nu^X)+  (x\,\one_{\{|x| > 1\}})\ast \mu^X. 
\end{equation} 
In this case, 
\begin{equation}\label{alpha_jumps} 
	\Delta \alpha_t = \left ( \int_{|x| \leq 1}x\,\hat{\nu}_t^X(dx)\right),\quad t \in [0,\,T], 
\end{equation} 
where $\hat{\nu}^X$ has been defined in \eqref{hatnu}. 
\end{proposition} 
\begin{remark}\label{R_alpha_orth}
Assume that \eqref{integral_repr} holds.
Then the process $\alpha$ is $(\mathcal{F}_t)$-orthogonal. Indeed,  for every  $(\mathcal{F}_t)$-local  martingale $N$,  
  $[A',N]= 0$ 
and $[\alpha^S,N] =  0$ by 
Proposition \ref{P_quad_var}. 	
\end{remark}

\proof 
If we suppose that  decomposition \eqref{integral_repr} holds, then $X$ is a particular weak Dirichlet process satisfying 
\[ 
X = M + V+ A', \quad   M =   M^c +  (x\,\one_{\{|x| \leq 1\}})\,\ast (\mu^X - \nu^X), \quad V =  \alpha^S +(x\,\one_{\{|x| > 1\}})\ast \mu^X. 
\] 
Conversely,  suppose that $X=M+V+A'$ is a particular weak Dirichlet process. 
Since $S=M+V$ is a  semimartingale, by  Theorem 11.25  in \cite{chineseBook}, 
 it can be decomposed as 
\[ S =   S^c+ \alpha^S + (x\,\one_{\{|x|\leq 1\}})\ast (\mu^S-\nu^S) +(x\,\one_{\{|x|> 1\}})\ast \mu^S, 
\] 
where $\mu^S$ is the jump measure of $S$ and $\nu^S$ is  the associated L\'evy system, 
 $S^c$ a  continuous local martingale, $\alpha^S$ a predictable process with finite variation such 
that $\alpha^S_0 = 0$ and 
\[ 
\Delta \alpha^S_s =\left ( \int_{|x| \leq 1}x\,\hat{\nu}^S_s(dx)\right). 
\] 
Consequently,  since $A'$ is adapted and continuous, with $A'_0 = 0$, we have 
$$ 
X= S + A'=   S^c+ (\alpha^S+ A') + (x\,\one_{\{|x|\leq 1\}})\ast (\mu^X-\nu^X) +(x\,\one_{\{|x|> 1\}})\ast \mu^X 
$$ 
and \eqref{integral_repr} holds  with $\alpha = \alpha^S + A'$  and  $M^c = S^c$. 
On the other hand, since $\Delta \alpha = \Delta \alpha^S$, 
 \eqref{alpha_jumps} follows. 
\endproof 
 
The following condition on  $X$ will play a fundamental role in the sequel: 
	\begin{equation}\label{CNS} 
	|x|\,\one_{\{|x| > 1\}} \ast \mu^X \,\,\in \mathcal A_{\rm loc}^+. 
	\end{equation} 
\begin{proposition}\label{P_char_part_weak_D} 
Let $X$ be  a particular $({\mathcal F}_t)$-weak Dirichlet process 
verifying the jump assumption \eqref{StandAss}. 
 $X$ is a special weak Dirichlet process if and only if \eqref{CNS} holds. 
\end{proposition} 
\proof 
Suppose the validity of \eqref{CNS}. We can decompose 
\[ 
(x\,\one_{\{|x| > 1\}})\ast \mu^X = (x\,\one_{\{|x| > 1\}})\ast (\mu^X-\nu^X) + (x\,\one_{\{|x| > 1\}})\ast \nu^X. 
\] 
By Proposition \ref{P_int_dec_par_WD}, using the notation of \eqref{integral_repr}, we obtain
	\begin{equation}\label{canonical_decomp} 
 X = M + A, \quad   M =   M^c + M^d, \quad A =  \alpha +(x\,\one_{\{|x| > 1\}})\ast \nu^X, 
	\end{equation} 
where $M^d= x\,\ast (\mu^X - \nu^X)$. 
First,  the process $\alpha + (x\,\one_{\{|x| > 1\}})\ast \nu^X$ is predictable.
Let  $N$ be a continuous local martingale. Then,  $[A,N]=0$ by Remark \ref{R_alpha_orth} 
and Proposition \ref{P_quad_var}, because $(x\,\one_{\{|x| > 1\}})\ast \nu^X$ has bounded variation. 
Consequently $X$ is a special Dirichlet process. 
 
 
Conversely, let  $X = M + V + A'$ be a particular weak Dirichlet process,  with $V$ bounded variation. 
We suppose that $X$ is a special weak Dirichlet process. Since $[A',N]=0$ for every continuous local martingale, then by additivity $X-A'$ is still a special weak Dirichlet process, $A'$ being continuous adapted. 
But $S:=X-A'= M+V$ is a semimartingale, and by Proposition 
\ref{P_SemiMtgSpecialWD} it is a special semimartingale. 
By Corollary 11.26 in \cite{chineseBook}, 
$|x|\,\one_{\{|x| > 1\}} \ast \mu^S \,\,\in \mathcal A_{\rm loc}^+$,
where $\mu^S$ is the jump measure of $S$. 
On the other hand, since $A'$ is continuous, $\mu^S $ coincides with $\mu^X$ and  \eqref{CNS} holds. 
\endproof 
 
 
 
We give the following result on the stochastic integration theory. 
\begin{proposition}\label{P_prop_int_stoc_munu} 
	Let $W\in \mathcal{G}^1_{\rm loc}(\mu^X)$, and define 
	$ 
	M^d_t=\int_{[0,t]\times \R} W_s(x)\,(\mu^X- \nu^X)(ds\,dx). 
	$ 
	Let moreover $(Z_t)$ be a predictable process such that 
	\begin{equation}\label{prop_Z} 
	\sqrt{\sum_{s \leq \cdot} Z_s^2 |\Delta M^d_s|^2} \in \mathcal{A}^+_{\rm loc}. 
	\end{equation} 
	Then \,\,$\int_{0}^\cdot Z_s\,d M^d_s$ \,\,is a local martingale and equals 
	\begin{equation}\label{g_munu_mtg} 
	\int_{]0,\cdot]\times \R} Z_s\,W_s(x)\,(\mu^X- \nu^X)(ds\,dx). 
	\end{equation} 
\end{proposition} 
\proof 
The conclusion  follows by the definition of the  stochastic integral \eqref{g_munu_mtg}, see Definition 3.63 in \cite{jacod_book}, 
provided we check the following three conditions. 
\begin{itemize} 
	\item[(i)]$\int_{0}^\cdot Z_s\,d M_s$ is a local martingale. 
	\item[(ii)] $\int_0^{\cdot} Z_s\,d M_s$ is a purely discontinuous local martingale; 
	in agreement with 
	Theorem 7.34,  in \cite{chineseBook}, we will show 
	$[\int_0^{\cdot} Z_s\,d M_s, N]=0$ for every $N$ continuous local martingale vanishing at zero. 
	\item[(iii)] 
	$\Delta \left(\int_0^{\cdot} Z_s\,d M_s\right)_t= \int_{\R} Z_t\,W_t(e)\,(\mu(\{t\}, de)-\nu(\{t\},de)), \quad t \in [0,\,T]$. 
\end{itemize} 
We prove now the validity of (i), (ii) and (iii). 
Condition \eqref{prop_Z} is equivalent to 
$ 
\sqrt{\int_0^t Z_s^2 \, d[M,M]_s} \in \mathcal{A}^+_{\rm loc}. 
$ 
According to Definition 
2.46 in \cite{jacod_book}, 
$\int_0^t Z_s\,dM_s$ is the unique local martingale satisfying 
\begin{equation}\label{JumpMart} 
\Delta \left(\int_0^{\cdot} Z_s\,d M_s\right)_t= Z_t\,\Delta M_t, \ t \in [0,T]. 
\end{equation} 
This implies in particular item (i). 
By Theorem 29, Chapter II, in \cite{protter}, it follows that 
$$ 
\left[\int_0^{\cdot}Z_s\,d M_s,N\right] = \int_0^{\cdot}Z_s\,d[M,N]_s, 
$$ 
and item (ii) follows because  $M$ is orthogonal to $N$, see Theorem 7.34,  in \cite{chineseBook} and the comments above, together with 
obvious localization arguments.
Finally, by   Definition 
3.63 in \cite{jacod_book}, 
taking into account \eqref{JumpMart}, 
$\Delta \left(\int_0^{\cdot} Z_s\,d M_s\right)_t$ equals 
$$ 
Z_t\,\Delta M_t = 
\int_{\R}Z_t \,W_t(e)\,(\mu(\{t\}, de)-\nu(\{t\}, de)), 
$$ 
for every $t \in [0,\,T]$, and this shows item (iii). 
\endproof

\begin{remark}\label{R_on_prp_Z} 
	Recalling that 
	$\sqrt{[M,M]_t} \in \mathcal{A}^+_{\rm loc}$ for any local martingale $M$ (see, e.g. Theorem 2.34 and Proposition 2.38 in \cite{jacod_book}), condition \eqref{prop_Z} is verified if for  instance  $Z$ is locally bounded. 
\end{remark}

\subsection{Stability of special weak Dirichlet processes under $C^{0,1}$ transformation} 


At this point, we   investigate  the stability properties of the class of  special weak Dirichlet processes. 
We start with an important property. 
\begin{proposition}\label{P_Md} 
	Let $X$ be a special weak Dirichlet process  with its canonical decomposition  $X=M^c+M^d+A$. We suppose that  assumptions \eqref{StandAss}, \eqref{CNS} 
	are verified. 
	Then 
	\begin{equation}\label{disc_mtg_munu} 
	M^d_s=\int_{]0,s]\times \R} x\,\,(\mu^X- \nu^X)(dt\,dx). 
	\end{equation} 
\end{proposition} 
\proof 
Taking into account assumption \eqref{StandAss},  Corollary \ref{R_x_smalljumps_G2(mu)} together with condition \eqref{CNS} insures  the fact that the 
right-hand side of \eqref{disc_mtg_munu} is well-defined. 
By definition, it is the unique  purely discontinuous local martingale whose jumps are indistinguishable from 
\[ 
\int_{\R} x\,\mu^X(\{t\},dx)-\int_{\R} x\,\nu^X(\{t\},dx). 
\] 
It remains to prove that 
\begin{equation} \label{E431} 
\Delta M^d_t  = \int_{\R} x\,\mu^X(\{t\},dx)-\int_{\R} x\,\nu^X(\{t\},dx),\,\, 
\textup{up to indistinguishability}. 
\end{equation} 
We have   
$\Delta M^d_t =  \Delta X_t - \Delta A_t$, for all $t \geq 0$.
We recall that 
$ 
\Delta X_t = \int_{\R}\,x\,\mu^X(\{t\},dx).
$ 
Thus \eqref{E431} can be established by additivity if we prove that, for any predictable time $\tau$, 
\begin{equation}\label{A_toprove}
\Delta A_{\tau}\,\one_{\{\tau < \infty\}} = \int_{\R}\,x\,\one_{\{\tau < \infty\}}\,\nu^X(\{\tau\},dx). 
\end{equation} 
Indeed in this case the Predictable Section Theorem (see e.g. Proposition 2.18, Chapter I, in \cite{JacodBook}) would insure that $\Delta A_t$ and $\int_{\R}\,x\,\nu^X(\{t\},dx)$ are indistinguishable. 

 Let us prove \eqref{A_toprove}. Let $\tau$ be a predictable time.
Let $(\tau_l)$  be a sequence of localizing stopping times, $\tau_l \rightarrow +\infty$ as $l$ goes to $+\infty$, such that $(M^d)^{\tau_l}$ is a martingale. In particular $(\Delta M^d)^{\tau_l}_\tau$ is integrable.
Since $\tau$ is a stopping time and  $A$ is a predictable process, $A_\tau$ is $\mathcal F_{\tau-}$-measurable by Proposition 2.4-a), Chapter I, in \cite{JacodBook}.
We set $\Omega_n :=\{A_\tau \leq n\}$: 
$\Omega_n \in \mathcal{F}_{\tau-}$, $\cup_n \Omega_n = \Omega$
a.s.,  so that 
 $\Delta (M^d)^{\tau_l}_\tau\,\one_{\{\tau < \infty\}}\,\one_{\Omega_n}$, 
$\Delta A^{\tau_l}_\tau\,\one_{\{\tau < \infty\}}\,\one_{\Omega_n}$ are integrable. 
By additivity, also 
\begin{equation}\label{Integr_deltaXtau}
\Delta X^{\tau_l}_\tau\,\one_{\{\tau < \infty\}}\,\one_{\Omega_n} \textup{ is integrable}. 
\end{equation}

Being $A$ predictable, $A^{\tau_l}$ (and therefore $\Delta A^{\tau_l}$) is also predictable,  Proposition 2.4-b), Chapter I, in \cite{JacodBook}. 
Since   $\tau$ is a  stopping time, again by item a) of the same proposition, $\Delta A_{\tau}^{\tau_l}\,\one_{\{\tau < \infty\}}\,\one_{\Omega_n}$ is $\mathcal F_{\tau-}$-measurable, so that
\begin{align} \label{A_eq}
\Delta A_{\tau}^{\tau_l}\,\one_{\{\tau < \infty\}}\,\one_{\Omega_n}&= \sper{\Delta X_{\tau}^{\tau_l}\,\one_{\{\tau < \infty\}}\,\one_{\Omega_n}-\Delta (M^d_{\tau})^{\tau_l}\,\one_{\{\tau < \infty\}}\,\one_{\Omega_n}|\mathcal{F}_{\tau-}}.
\end{align}
Now, by Corollary 2.31, Chapter I,  in \cite{JacodBook},  recalling the definition of predictable projection, for any  martingale $L$ and  any predictable time $\tau$, we have 
$\sper{\Delta L_\tau \,\one_{\{\tau < \infty\}}|\mathcal{F}_{\tau-}}=0$.
Therefore, by Remark \ref{R_Jumps_measure}, identity \eqref{A_eq}  gives 
\begin{align}\label{A_eq2}
\Delta A_{\tau}^{\tau_l}\,\one_{\{\tau < \infty\}}\,\one_{\Omega_n}& = \sper{\int_{\R}\,x\,\one_{\{\tau < \infty\}}\,\one_{\Omega_n}\,\mu^{X^{\tau_l}}(\{\tau\},dx)|\mathcal{F}_{\tau-}} \nonumber\\ 
& = \int_{\R}\,x \,\one_{\{\tau < \infty\}}\,\one_{\Omega_n}\,\nu^{X^{\tau_l}}(\{\tau\},dx)\quad \textup{a.s.}, 
\end{align} 
where for the latter equality we have used 
Proposition 
	1.17, point b), Chapter II, in \cite{JacodBook}.
By the subsequent Remark \ref{R_int_modulo},  the right-hand side of \eqref{A_eq2} is  well-defined.

Now we remark that, for any  nonnegative random field $U \in \tilde{\mathcal P}$, 
\begin{align*}
	\sper{\int_{[0,\,T] \times\R}\,U_s(x)\,\one_{]]0,\,\tau_l]]}(s)\,\nu^{X^{\tau_l}}(ds,dx) }&=\sper{\int_{[0,\,T] \times\R}\,U_s(x)\,\one_{]]0,\,\tau_l]]}(s)\,\mu^{X^{\tau_l}}(ds,dx) }\\
	&=\sper{\sum_{s \leq T}\,U_s(\Delta X_s^{\tau_l})\,\one_{]]0,\,\tau_l]]}(s)\,\one_{\{\Delta X^{\tau_l}_s\neq 0\}}}\\
	&=\sper{\sum_{s \leq T}\,U_s(\Delta X_s)\,\one_{]]0,\,\tau_l]]}(s)\,\one_{\{\Delta X_s \neq 0\}}} \\
	&=\sper{\int_{[0,\,T] \times\R}\,U_s(x)\,\one_{]]0,\,\tau_l]]}(s)\,\mu^{X}(ds,dx) }\\
&=\sper{\int_{[0,\,T] \times\R}\,U_s(x)\,\one_{]]0,\,\tau_l]]}(s)\,\nu^{X}(ds,dx) }.
\end{align*}
Therefore, 
\begin{equation}\label{id_nu}
\one_{]]0,\,\tau_l]]}(s)\,\nu^{X^{\tau_l}}(ds,dx) =\one_{]]0,\,\tau_l]]}(s)\,\nu^{X}(ds,dx).
\end{equation}
So, identity \eqref{A_eq2} becomes
\begin{align}\label{5.36bis}
\Delta A_{\tau\wedge \tau_l}\,\one_{]]0,\,\tau_l]]}(s)\,\one_{\{\tau < \infty\}}\,\one_{\Omega_n}& = \int_{\R}\,x\,\nu^{X}(\{\tau\},dx) \,\one_{]]0,\,\tau_l]]}(s)\,\one_{\{\tau < \infty\}}\,\one_{\Omega_n}\quad \textup{a.s.}
\end{align}
We  prove that 
$
\int_{\R}\,|x|\,\nu^{X}(\{\tau\},dx) 
$
is finite a.s.
For every $n\in\N$, we set $\tilde{\Omega}^{n}:= \{\tau_n > T\} \cap \Omega_n$. Obviously, $\Omega = \cup_{n} \tilde{\Omega}^{n}$ a.s. We have
\begin{align*}
	\sper{\int_{\R}\,|x|\,\nu^{X}(\{\tau\},dx)\,\one_{\tilde{\Omega}^{n}}} &= \sper{\int_{\R}\,|x|\,\nu^{X}(\{\tau\},dx)\,\one_{]]0,\,\tau_n]]}(s)\,\one_{\tilde{\Omega}_n}}\\
	&= \sper{\int_{\R}\,|x|\,\nu^{X^{\tau_n}}(\{\tau\},dx)\,\one_{\tilde{\Omega}_n}},
\end{align*}
which is finite by Remark \ref{R_int_modulo}.
The latter equality follows by \eqref{id_nu}.
This implies that  $
\int_{\R}\,|x|\,\nu^{X}(\{\tau\},dx) 
$
is finite a.s. on $\tilde{\Omega}^{n}$ and therefore on  $\Omega$.
The conclusion follows letting first $n$, then $l$, go to infinity in \eqref{5.36bis}.

\endproof

\begin{remark}\label{R_int_modulo}
For every $l, n \in \N$, we have 
	$$
	\sper{\int_{\R}|x| \,\one_{\{\tau < \infty\}}\,\one_{\Omega_n}\,\nu^{X^{\tau_l}}(\{\tau\},dx)} < \infty.
	$$
	This follows from  Proposition 
	1.17, point b), Chapter II, in \cite{JacodBook}, Remark \ref{R_Jumps_measure}, and \eqref{Integr_deltaXtau}. 
In particular this ensures that the right-hand side of \eqref{A_eq2} is
well-defined.
\end{remark}

 

	\begin{lemma}\label{L_int_cond_special_wD} 
		Let $X$ be a c\`adl\`ag process  satisfying 
condition \eqref{CNS}. Let also   $F:[0,T] \times \R \rightarrow \R$  be  a function of class $C^{0,1}$ such that 
		\begin{equation}\label{E_C} 
			\int_{]0,s]\times \R} |F(t,X_{t-} +x )-F(t,X_{t-})-x\,\partial_x F(t,X_{t-})|\,\one_{\{|x| >1\}}\,\mu^X(dt\,dx) \in \mathcal{A}^{+}_{\rm loc}. 
		\end{equation} 
		Then	 
		\begin{align} 
		&\int_{]0,s]\times \R} x\,\partial_x F(t,X_{t-})\one_{\{|x| >1\}}\,\mu^X(dt\,dx)\,\,\in \mathcal{A}^{+}_{\rm loc},\label{x_derivx_Aloc}\\ 
		&\int_{]0,s]\times \R} |F(t,X_{t-} +x )-F(t,X_{t-})|\,\one_{\{|x| >1\}}\,\mu^X(dt\,dx) \in \mathcal{A}^{+}_{\rm loc}.\label{E_C_old} 
		\end{align} 
	\end{lemma} 
	 \begin{remark}\label{R_bounded_derivative} 
	 	Condition \eqref{E_C} is automatically verified if 
	 	$X$ is a c\`adl\`ag 
 process  satisfying  \eqref{CNS}  and  $F:[0,T] \times \R \rightarrow \R$   is a  function of $C^{0,1}$ class with  $\partial_x F$ bounded. 
	 \end{remark} 
	\proof 
	 Condition \eqref{CNS} together the fact that the process $(\partial_x F(t,X_{t-}))$ is locally bounded implies \eqref{x_derivx_Aloc}; then condition \eqref{E_C_old} follows from   \eqref{E_C} and \eqref{x_derivx_Aloc}. 
 \endproof

\begin{theorem}\label{T_C1_dec_specialweak_Dir} 
		Let $X$ be a special weak Dirichlet process  of finite quadratic variation  with its canonical decomposition  $X=M^c+M^d+A$. 
	We denote by $\tilde{C}^{0,1}$ the space of functions $F:[0,T] \times \R \rightarrow \R$ 
	satisfying condition \eqref{E_C}. Then we have the following. 
	\begin{itemize} 
			\item[(1)] 
	For every  $F$
	of class  $\tilde{C}^{0,1}$,  
			$Y = F(\cdot,X)$ is a special weak Dirichlet process, with decomposition $Y= M^F + A^F$, where 
	\begin{align*} 
	M^F_t&= F(0,X_0)+ \int_{]0,t]}\partial_x F(s,X_s)\,d (M^c+M^d)_s\\ 
	&+ \int_{]0,\,t]\times \R} (F(s,X_{s-}+x)-F(s,X_{s-})-x\,\partial_x F(s,X_{s-}))\,(\mu^X-\nu^X)(ds\,dx), 
	\end{align*} 
	and $F \mapsto A^F$,  
	$\tilde{C}^{0,1}\rightarrow \mathbb{D}^{ucp}$,   is a linear map and, for every $F \in \tilde{C}^{0,1}$, $A^F$ is a predictable $(\mathcal F_t)$-orthogonal  process. 
	\item[(2)] If moreover condition \eqref{CNS} holds, 
	$M^F$ reduces to   
			\begin{align*} 
			M^F_t= F(0,X_0)+\int_0^t\partial_x F(s,X_s)\,d M^c_s+ \int_{]0,\,t]\times \R} (F(s,X_{s-}+x)-F(s,X_{s-}))\,(\mu^X-\nu^X)(ds\,dx). 
			\end{align*} 
		\end{itemize} 
\end{theorem} 

\proof 
(1)  For every $F$ of class $\tilde{C}^{0,1}$, we set  
\begin{equation} 
\label{E_2Gamma} 
A^F = \Gamma^F  + \bar V^F, 
\end{equation} 
where $\Gamma^F$ has been defined in Theorem \ref{T_C1_dec_weak_Dir}, and  
$$ 
\bar V^F_t := \int_{]0,\,t]\times \R} (F(s,X_{s-}+x)-F(s,X_{s-})-x \,\partial_x F(s,X_{s-}))\,\one_{\{|x| > 1\}}\,\nu^X(ds\,dx), 
$$ 
which is well-defined by assumption \eqref{E_C}.  
 
The map $F \mapsto A^F$ is linear since $F \mapsto \Gamma^F$ and $F \mapsto \bar V^F$ are linear. 
Given $F \in \tilde{C}^{0,1}$, $A^F$ is a $(\mathcal F_t) $-orthogonal process 
by Theorem \ref{T_C1_dec_weak_Dir}-(a),  
taking into account that  $[\bar V^F,N]=0$ by Proposition \ref{P_quad_var}. 
Using  decomposition \eqref{E_2Gamma}, Theorem \ref{T_C1_dec_weak_Dir}-(b) and the fact that $\bar V$ is predictable, it follows that   $A^F$ is  predictable. 
 
(2) follows from (1), if we  show that  
\[ 
\int_{0}^t \partial_x F(s,X_{s-})\,d M^d_s=\int_{]0,t]\times \R} x\,\partial_x F(s,X_{s-})\,(\mu^X- \nu^X)(ds\,dx), 
\] 
This follows from  Proposition \ref{P_prop_int_stoc_munu} and  Proposition \ref{P_Md}, taking into account Remark 
\ref{R_on_prp_Z}.	 
 

\endproof 
\begin{remark}\label{C_C1_dec_specialweak_Dir+BV} 
In Theorem \ref{T_C1_dec_specialweak_Dir}-(2),  condition  \eqref{CNS} is verified for instance if $X$ 
is a particular weak Dirichlet process, see  Proposition \ref{P_char_part_weak_D}. 
\end{remark}

 When $X$ is a special weak Dirichlet process  with finite quadratic  and finite energy,  and $F: \R \rightarrow \R$ is a function of class $C^1_b$,  Corollary 3.2 in \cite{cjms}  states that $F(X)$ is again a special weak Dirichlet process. 
 Below we extend that result.
\begin{corollary}\label{C_Coquet}
Let $X$ be a special weak Dirichlet process  
of finite quadratic variation  
with its canonical decomposition  $X=M^c+M^d+A$, and such that $\sper{\sum_{s \leq T} |\Delta X_s|^2 } < \infty$. 
Let $F:[0,T] \times \R \rightarrow \R$ of class $C^{0,1}$, with $\partial_x F$ bounded. 
Then $Y_t = F(t,X_t)$ is a special weak Dirichlet process, with decomposition $Y= M^F + A^F$, where  
$A^F$ is a predictable $(\mathcal F_t)$-orthogonal  process, and  
\begin{align*} 
			M^F_t= F(0,X_0)+\int_0^t\partial_x F(s,X_s)\,d M^c_s+ \int_{]0,\,t]\times \R} (F(s,X_{s-}+x)-F(s,X_{s-}))\,(\mu^X-\nu^X)(ds\,dx). 
\end{align*} 
\end{corollary}
\proof
The result  will be a consequence of Theorem \ref{T_C1_dec_specialweak_Dir}, provided we verify conditions \eqref{CNS} and \eqref{E_C}. Condition \eqref{CNS}  follows from 
\begin{align*}
	\sum_{s \leq T} |\Delta X_s|\,\one_{\{|\Delta X_s| > 1\}} &\leq \sum_{s \leq T} |\Delta X_s|^2\,\one_{\{|\Delta X_s| > 1\}} \leq \sum_{s \leq T} |\Delta X_s|^2 \in L^1(\Omega). 
\end{align*}
Condition \eqref{E_C} follows by Remark \ref{R_bounded_derivative}. 
\endproof

\subsection{The case  of special weak Dirichlet processes 
 without continuous local martingale.} 
We end this section by considering the  case of 
 special weak Dirichlet processes with canonical decomposition $X=M+A$ 
where $M=M^d$ is a purely discontinuous local martingale. In particular 
there is   no continuous martingale part. 
In this framework, under the assumptions of   
Theorem  \ref{T_C1_dec_specialweak_Dir}, if assumption \eqref{CNS} in  verified, then item (2) says that  
\begin{equation} \label{C01_special_WD_formulaMD} 
 	F(t,X_t)=F(0,X_0) 
	+ \int_{]0,\,t]\times \R} (F(s,X_{s-}+x)-F(s,X_{s-}))\,(\mu^X-\nu^X)(ds\,dx)+A^F(t). 
\end{equation} 
Since in the above formula no derivative  appears, a natural question appears: 
 is it possible to state  a chain rule  \eqref{C01_special_WD_formulaMD} when 
$F$ is not of class  $C^{0,1}$? Indeed we have the following result, which does not require any weak Dirichlet structure on $X$.

We first introduce some notations. Let $E$ be a closed subset of $\R$ on which $X$ takes values. Given a  c\`adl\`ag function $\varphi: [0,\,T] \rightarrow \R$, we denote by $\mathcal C_{\varphi}$ the set of times $t \in [0,\,T]$ for which there is a left (resp. right) neighborhood $I_{t-} = ]t-\varepsilon,\,t[$ (resp. $I_{t+} = [t,\,t+\varepsilon[$) such that $\varphi$ is  constant on $I_{t-}$ and $I_{t+}$.
We introduce the following assumption.
\begin{hypothesis}\label{H_new}
There exists $\mathcal C \in [0,\,T]$ such that for $\omega$ a.s. $\mathcal C \supset \mathcal C_{X(\omega)}$, and 
\begin{itemize}
		\item $\forall t \in \mathcal C$, $t \mapsto F(t,x)$ is continuous $\forall x \in E$;
		\item $\forall t \in \mathcal C^c$, $x \in E$,  $(t, x)$ is a continuity point of $F$.
\end{itemize}	
\end{hypothesis} 

\begin{remark}\label{R_1}
	Hypothesis \ref{H_new}   is fulfilled in two typical situations.
	\begin{enumerate}
	\item $\mathcal C=[0,\,T]$. Almost surely $X$ admits a finite number of jumps and $t \mapsto F(t,x)$ is continuous $\forall x \in E$.
	\item $\mathcal C=\emptyset$ and  $F|_{[0,\,T] \times E}$ is continuous.
	\end{enumerate}
\end{remark}
From now on, we denote by $\Delta F(t,X_t)$ the jump of the time $t$ the process $(F(t,X_t))$.

\begin{remark}\label{R_2}
Assume that Hypothesis \ref{H_new} holds. Then
\begin{itemize}
\item[(i)]	
$F(t, X_t)$ is necessarily a  c\`adl\`ag process.
\item[(ii)] $\forall t \in [0,\,T]$, $\Delta F(t, X_t)= F(t,X_t)- F(t, X_{t-})$.
\end{itemize}
\end{remark}

\begin{proposition}\label{P_C00_chain_rule} 
	Let $X$ be an adapted c\`adl\`ag process. 
	Let   $F:[0,T] \times \R \rightarrow \R $ be a 
function satisfying Hypothesis \ref{H_new}.
Assume that the following holds. 
	\begin{itemize} 
		\item[(i)] $F(t,X_t)$ is an $(\mathcal F_t)$-orthogonal process such that $\sum_{s \leq T} |\Delta F(s,X_s)| < \infty$, a.s.
		\item[(ii)] $\int_{]0,\,\cdot]\times \R}|F(s,X_{s-}+x)-F(s,X_{s-})|\,\mu^X(ds\,dx)\,\,\in \mathcal{A}^+_{\rm loc}$. 
	\end{itemize} 
	Then 
	\begin{equation}\label{C00_decomp} 
		F(t,X_t)=F(0,X_0) 
		+ \int_{]0,\,t]\times \R} (F(s,X_{s-}+x)-F(s,X_{s-}))\,(\mu^X-\nu^X)(ds\,dx)+A^F(t), 
	\end{equation} 
	where $A^F$ is a predictable $(\mathcal F_t)$-orthogonal  process. In particular, $F(t,X_t)$ is a special weak Dirichlet process. 
\end{proposition} 
\proof
Since by 
condition (i),  $\sum_{s \leq t}|\Delta F(s,X_s)|$ is finite a.s.,  the process 
$ 
B_t =
\sum_{s \leq t}\Delta F(s,X_s) $ is well-defined and 
has bounded variation.
We set $A'_t := F(t,X_t)
-B_t$.
$ A'$ 
is a continuous process, and is 
$(\mathcal F_t)$-orthogonal by additivity,  since  by assumption $F(t,X_t)$ is $(\mathcal F_t)$-orthogonal and  $B$ is $(\mathcal F_t)$-orthogonal by  Proposition \ref{P_quad_var}.
Recalling the definition of the jump measure $\mu^X$, 
by Remark \ref{R_2}
and  condition (ii), we get 
\begin{align*} 
B_t&=\sum_{0 <s \leq t}(F(s,X_{s-}+\Delta X_s)-F(s,X_{s-}))\\ 
&=\int_{]0,\,t]\times \R}(F(s,X_{s-}+x)-F(s,X_{s-}))\,\mu^X(ds\,dx)\\ 
&=\int_{]0,\,t]\times \R}(F(s,X_{s-}+x)-F(s,X_{s-}))\,(\mu^X-\nu^X)(ds\,dx)\\ 
&+\int_{]0,\,t]\times \R}(F(s,X_{s-}+x)-F(s,X_{s-}))\, 
\nu^X(ds\,dx). 
\end{align*} 
Finally, decomposition \eqref{C00_decomp} holds with 
\begin{align}\label{Gamma_C0} 
A^F(t):=A'_t +\int_{]0,\,t]\times \R}(F(s,X_{s-}+x)-F(s,X_{s-}))\, 
\nu^X(ds\,dx). 
\end{align} 
The  process $A^F$ in \eqref{Gamma_C0} is clearly predictable. The $(\mathcal F_t)$-orthogonality property of  $A^F$ follows from the orthogonality of $A'$ and by Proposition \ref{P_quad_var},  noticing that the integral term in  \eqref{Gamma_C0} is a bounded variation process. 
\qed

We provide below  an example of $(\mathcal F_t)$-weak Dirichlet process, intervening in the framework of Piecewise Deterministic Processes, see e.g., Section 4.3  in \cite{BandiniRusso2}.
\begin{proposition}\label{P_P}
	Let $(T_n)_{n \geq 1}$ be a sequence of strictly increasing stopping times, such that $T_n \uparrow \infty$ a.s. Let $X$ be a (c\`adl\`ag) process of the type
	$$
	X_t = \sum_{n = 1}^\infty \one_{[T_{n-1}, \,T_n[}(t)\, R_{n-1}(t), 
	$$
	where $R_n$ is an $(\mathcal F_{T_n})$-measurable c\`adl\`ag process. Let $F: [0,\,T] \times \R \rightarrow \R$. We suppose one of the two following assumptions. 
	\begin{itemize}
	\item[(i)] $R_n(t)$ are constant $(\mathcal F_{T_n})$-random variables, $t \mapsto F(t,x), \,\, \forall x \in E$, is continuous.
	\item[(ii)]	$F$ is continuous. 
	\end{itemize}
	Then $F(\cdot, X)$ is $(\mathcal F_t)$-orthogonal.
\end{proposition}
\proof
Let $N$ be an $(\mathcal F_t)$-continuous local  martingale.  
By Proposition \ref{P_Ito_cov}-(ii), 
\begin{equation}\label{first}
\int_0^t F(s, X_{s-})\,(N_{(s+ \varepsilon) \wedge t}-N_s)\,\frac{ds}{\varepsilon}\underset{\varepsilon \rightarrow 0}{\longrightarrow} \,\, \int_0^t F(s, X_{s-})\,d N_s \quad \textup{u.c.p.} 
\end{equation}
 By the definition of covariation
it will be enough to prove 
\begin{equation}\label{(2)}
	\int_0^t F((s+ \varepsilon)\wedge t, X_{(s + \varepsilon) \wedge t})\,(N_{(s+ \varepsilon) \wedge t}-N_s)\,\frac{ds}{\varepsilon} \underset{\varepsilon \rightarrow 0}{\longrightarrow} \,\, \int_0^t F(s, X_{s-})\,d N_s \quad \textup{u.c.p.} 
\end{equation}
Indeed, by subtracting  \eqref{first} from \eqref{(2)}, we  would obtain $[F(\cdot, X), \,N]=0$.\\
\eqref{(2)} will be the consequence of 
\begin{equation}\label{(1)}
	\int_0^t F(s+ \varepsilon, \tilde{X}_{s + \varepsilon})\,(N_{(s+ \varepsilon) \wedge t}-N_s)\,\frac{ds}{\varepsilon} \underset{\varepsilon \rightarrow 0}{\longrightarrow} \,\, \int_0^t F(s, X_{s-})\,d N_s \quad \textup{u.c.p.},
\end{equation}
where $\tilde X$ is an extension of $X$ to $]T, \infty[$ by $X_T$, and $F$ is extended to $\R_+ \times \R$, setting 
$F(s,x)= F(T,x)$, if $s >T$, $x \in \R$. 
This happens because   $N$ is uniformly continuous on $[0,\,T]$ and $F$ is locally bounded.\\ 
Let us now  concentrate on proving \eqref{(1)}. We set $\tau_n := T_n \wedge T$. By convention, we set $T_0=0$.  The left-hand side of \eqref{(1)} gives $J_1(t, \varepsilon) + J_2(t,\varepsilon)$,
where 
\begin{align*}
	J_1(t, \varepsilon) &= \sum_{n=1}^\infty\,\int_0^t \one_{[\tau_{n-1}, \,\tau_{n}-\varepsilon[}(s)\, F(s+ \varepsilon, \tilde{X}_{s + \varepsilon})\,(N_{(s+ \varepsilon) \wedge t}-N_s)\,\frac{ds}{\varepsilon}\\
	J_2(t, \varepsilon) &= \sum_{n=1}^\infty\,\int_0^t \one_{[\tau_{n}-\varepsilon,\,\tau_{n}[}(s)\, F(s+ \varepsilon, \tilde{X}_{s + \varepsilon})\,(N_{(s+ \varepsilon) \wedge t}-N_s)\,\frac{ds}{\varepsilon},
\end{align*}
with the convention that $\one_{[a,\,b[}(s)\equiv 0$ if $b \leq a$. Now, since there is a.s. only a finite number of $(\tau_n)$ in $[0,\,T]$, we have a.s.
$$
\sup_{t \in [0,\,T]} |J_2(t, \varepsilon)| \leq ||F||_\infty\,\delta(N, \varepsilon) \underset{\varepsilon \rightarrow 0}{\longrightarrow} \,\,0,
$$
where $\delta(N, \cdot)$ denotes the modulus of continuity of $N$. Concerning $J_1(\varepsilon, t)$, it will be enough to show that 
\begin{equation}\label{(3)}
	\int_0^t \mathcal Y^\varepsilon_s\,(N_{(s+ \varepsilon) \wedge t}-N_s)\,\frac{ds}{\varepsilon} \underset{\varepsilon \rightarrow 0}{\longrightarrow} \,\, \int_0^t \sum_{n =1}^\infty \one_{[\tau_{n-1}, \,\tau_{n}[}(s)\, F(s, X_{s-})\,d N_s \quad \textup{u.c.p.},
\end{equation}
and 
\begin{eqnarray*} 
		\mathcal Y^\varepsilon_s= 
		\left\{ 
		\begin{array}{ll} 
			F(s+ \varepsilon, R_{n-1}(s + \varepsilon)) &\quad s \in [\tau_{n-1}, \,\tau_{n}-\varepsilon[,\\ 
			0 &\quad  s \in [\tau_{n}-\varepsilon, \,\tau_{n}[.
		\end{array} 
		\right. 
	\end{eqnarray*}
	Clearly, $\one_{[\tau_{n-1}, \,\tau_{n}-\varepsilon[}(s)\,Y^\varepsilon_s$ is $\mathcal F_s$-measurable. So $Y^\varepsilon$ is $(\mathcal F_s)$-adapted. Since, by assumption, $Y^\varepsilon$ is c\`adl\`ag, it is then progressively measurable. The left-hand side of \eqref{(3)}, using stochastic Fubini's theorem, gives 
	\begin{equation}\label{(4)}
		\int_0^t d N_u \int_{(u- \varepsilon)-}^u \mathcal Y^\varepsilon_s \,\frac{
		ds}{\varepsilon}.
	\end{equation}
	To show that \eqref{(4)} goes u.c.p. to the right-hand side of \eqref{(3)}, taking into account Problem 2.26, Chapter 3, in \cite{Karatzas1991Brownian}, it is enough to show that 
$$
\int_0^T d [N, N]_u \left(\int_{(u-\varepsilon)-}^u  \mathcal Y^\varepsilon_s \, \frac{
ds}{\varepsilon} - \sum_{n =1}^\infty \one_{[\tau_{n-1}- \varepsilon, \,\tau_{n}[}(u)\, F(u, X_{u}) \right)^2 \underset{\varepsilon \rightarrow 0}{\longrightarrow} 0 \quad \textup{in probability.}
$$
That convergence will be proved to be also a.s. Indeed 
$|\mathcal Y^\varepsilon_s| \leq ||F||_{\infty}$ and $|F(s, X_{s-})| \leq ||F||_{\infty}$, 
where $||F||_{\infty}:= \sup_{t \in [0,\,T]} |F(t,x)|$, with $x \in \textup{Im}(X_t, \,t \in [0,\,T])$. Moreover, for fixed $\omega$, for $s \in [\tau_{n-1},\,\tau_n[$, for $\varepsilon$ small enough, 
\begin{align*}
	\mathcal Y^\varepsilon_s = F(s+ \varepsilon, \, R_{n-1}(s + \varepsilon)) &\longrightarrow F(s, \,R_{n-1}(s)) =  F(s, X_{s-}),
\end{align*}
in both cases (i) and (ii) of the assumption. 
This shows the validity of \eqref{(1)} and completes the proof. 
\endproof
\begin{example}
Previous proposition provides further examples of weak-Dirichlet processes.
In particular, let us consider the two following cases, developed respectively  in Sections 4.2 and 4.3 of \cite{BandiniRusso2}, see also \cite{BandiniFuhrman}, \cite{BandiniConfortola}, \cite{BandiniPDMPs}, \cite{BandiniBSDE}.
\begin{enumerate}
\item $X$ pure jump process.
\item $X$ Piecewise Deterministic Markov Processes.
\end{enumerate}
Then Proposition \ref{P_P} applies to the case 1. with assumption (i), and to the case 2. with assumption (ii).
\end{example}

\appendix 
\renewcommand\thesection{Appendix} 
\section{} 
\renewcommand\thesection{\Alph{subsection}} 
\renewcommand\thesubsection{\Alph{subsection}} 
 

\subsection{Additional results on calculus via regularization}\label{Section_Appendix} 
 
In what follows, we are given a filtered  probability space $(\Omega,\mathcal{F}, (\mathcal{F}_t),\P)$, and  an integer-valued random measure $\mu$.

For every functions $f, g$ defined on $\R$, let now set 
\begin{align} 
\tilde I^{-}(\varepsilon, t,f, dg)=&\int_{]0,\,t]} f(s)\,\frac{g(s+\varepsilon)-g(s)}{\varepsilon}\,ds,\label{Int_cont}\\ 
C_{\varepsilon}(f,g)(t)=& \frac{1}{\varepsilon}\,\int_{]0,\,t]}\,(f(s+\varepsilon)-f(s)) (g(s+\varepsilon)-g(s))\,ds.\label{Cov_cont} 
\end{align} 
\begin{definition}\label{D_forw_int} 
	Assume  that $X, Y$ are   two c\`adl\`ag processes. We say that \textbf{the forward integral of $Y$ with respect to $X$ exists in the pathwise sense}, 
	if there exists some process $(I(t),t\geq 0)$ such that, 
	for all subsequences $(\varepsilon_n)$, there is a subsequence $(\varepsilon_{n_k})$ and a null set $\mathcal{N}$:
	\begin{displaymath} 
	\forall \omega \notin \mathcal{N}, \qquad \lim_{k \rightarrow \infty}\,|\tilde I^{-}(\varepsilon_{n_k}, t,Y, dX)(\omega)- I(t)(\omega)|=0 
	\qquad \forall t \geq 0.
	\end{displaymath} 
\end{definition} 
\begin{definition}\label{D_cov} 
	Let  $X, Y$ be two c\`adl\`ag processes. 
	\textbf{the covariation between $X$ and $Y$ (the quadratic variation of $X$) 
		exists in the pathwise sense}, 
	if there exists a c\`adl\`ag process $(\Gamma(t),t\geq 0)$ such that, for all subsequences $(\varepsilon_n)$ 
	there is a subsequence $(\varepsilon_{n_k})$ and a null set $\mathcal N$: 
	\begin{displaymath} 
	\forall \omega \notin \mathcal N, \qquad \lim_{k \rightarrow \infty}\,|C_{\varepsilon_{n_k}}(X,Y)(t)(\omega)- \Gamma(t)(\omega)|=0 
	\qquad \forall t \geq 0.
	\end{displaymath} 
\end{definition} 
\begin{proposition}\label{P_conv_as} 
	Let $X,Y$ be two c\`adl\`ag processes. Then 
	\begin{eqnarray} 
	I^{-ucp}(\varepsilon, t, Y, dX) &=& \tilde I^{-}(\varepsilon, t, Y, dX) + R_1(\varepsilon, t)\label{rel_int}\\ 
	\left[X,Y\right]^{ucp}_{\varepsilon}(t)  &=&  C_{\varepsilon}(X,Y)(t) + R_2(\varepsilon, t)\label{rel_cov}, 
	\end{eqnarray} 
	where 
	\begin{equation} \label{E_conv_unif} 
	R_i(\varepsilon, t)(\omega) \underset{\varepsilon \rightarrow 0}{\longrightarrow} 0\quad i=1,2,\,\,\forall t \in [0,\,T], \,\, \forall \omega \in \Omega. 
	\end{equation} 
	Moreover, if  $X$  is continuous, then the convergence in \eqref{E_conv_unif} holds u.c.p. 
\end{proposition} 
\proof 
We fix $t \in [0,\,T]$. Let $\gamma > 0$. 
The definition of  right continuity in $t$  insures that  there exists  $\delta >0$ small enough such that 
\begin{align*} 
|X(t) - X(a)| \le \gamma \quad \textup{ if }\,\, a-t < \delta, \,a>t,\\ 
|Y(t) - Y(a)| \le \gamma \quad \textup{ if }\,\, a-t < \delta, \,a>t. 
\end{align*} 
We start  proving \eqref{rel_int}. 
From  decomposition \eqref{int_ucp_dec} and the definition of $\tilde I^{-}(\varepsilon, t, Y, dX)$  we get 
\begin{eqnarray*} 
	I^{-ucp}(\varepsilon, t, Y, dX)-\tilde I^{-}(\varepsilon, t, Y, dX) 
	&=&\frac{1}{\varepsilon}\int_{(t-\varepsilon)_+}^{t}Y(s)\,[X(t)-X(s)]\,ds\\ 
	&&-\frac{1}{\varepsilon}\int_{(t-\varepsilon)_+}^{t}Y(s)\,[X(s+ \varepsilon)- X(s)]\,ds\\ 
	&=& 
	\frac{1}{\varepsilon}\int_{(t-\varepsilon)_+}^{t}Y(s)\,[X(t)-X(s+ \varepsilon)]\,ds 
	=: R_1(\varepsilon, t). 
\end{eqnarray*} 
Choosing $\varepsilon < \delta$ we get 
$|R_1(\varepsilon, t)| \leq \gamma\, ||Y||_{\infty}$,
and since $\gamma$ is arbitrary, we conclude that $R_1(\varepsilon, t)\rightarrow  0$ as $\varepsilon$ goes to zero, for every $t \in [0,\,T]$. 
 
It remains to show \eqref{rel_cov}. To this end we evaluate 
\begin{align*} 
[X,Y]_{\varepsilon}^{ucp}(t)-C_{\varepsilon}(X,Y)(t) 
&= \frac{1}{\varepsilon}\int_{(t-\varepsilon)_+}^{t}[X(t)-X(s)]\,[Y(t)-Y(s)]\,ds\\ 
&-\frac{1}{\varepsilon}\int_{(t-\varepsilon)_+}^{t}[X(s+ \varepsilon)- X(s)]\,[Y(s+ \varepsilon)- Y(s)]\,ds\\ 
&=: R_2(\varepsilon, t). 
\end{align*} 
We have 
\begin{align*} 
R_2(\varepsilon,t) 
&= \frac{1}{\varepsilon}\int_{(t-\varepsilon)_+}^{t}[X(t)-X(s)]\,[Y(t)-Y(s)]\,ds\\ 
&-\frac{1}{\varepsilon}\int_{(t-\varepsilon)_+}^{t}[X(s+ \varepsilon)- X(s)]\,[Y(t)- Y(s)]\,ds\\ 
&+ \frac{1}{\varepsilon}\int_{(t-\varepsilon)_+}^{t}[X(s+ \varepsilon)- X(s)]\,[Y(t)- Y(s)]\,ds\\ 
&-\frac{1}{\varepsilon}\int_{(t-\varepsilon)_+}^{t}[X(s+ \varepsilon)- X(s)]\,[Y(s+ \varepsilon)- Y(s)]\,ds\\ 
&= \frac{1}{\varepsilon}\int_{(t-\varepsilon)_+}^{t}[X(t)-X(s+ \varepsilon)]\,[Y(t)-Y(s)]\,ds\\ 
&+\frac{1}{\varepsilon}\int_{(t-\varepsilon)_+}^{t}[X(s+ \varepsilon)- X(s)]\,[Y(t)- Y(s+ \varepsilon)]\,ds. 
\end{align*} 
Choosing $\varepsilon < \delta$, the absolute value of previous expression is smaller than 
$2 \,\gamma \,(||Y||_{\infty}+||X||_{\infty})$.
Since $\gamma$ is arbitrary, $R_2(\varepsilon,t) \rightarrow  0$ as $\varepsilon$ goes to zero, 
for every $t \in [0,\,T]$.
Suppose now that $X$ is continuous. The expression of $R_2(\varepsilon,t)$ 
can be uniformly (in $t$)  bounded by $2 \delta(X, \varepsilon)\, \Vert Y \Vert_\infty$, where $\delta(X,\cdot)$ denotes the modulus of continuity of $X$; 
on the other hand $R_1(\varepsilon,t) \le 2 \delta(X, \varepsilon)\, \Vert Y \Vert_\infty, \forall t \in [0,T]$. 
This concludes the proof of Proposition \ref{P_conv_as}. 
\endproof 
\begin{corollary}\label{C_id_RV_fw} 
	Let $X,Y$ be two c\`adl\`ag processes. 
	\begin{itemize} 
		\item[1)] If the forward integral of $Y$ with respect to $X$ exists, then it exists in the pathwise sense. In particular,  there is a null set $\mathcal N$ and, 
		for any sequence $(\varepsilon_n)\downarrow 0$, 
		a subsequence $(\varepsilon_{n_k})$ such that 
		\begin{eqnarray}\label{conv_ucp_int} 
		\tilde I^{-}(\varepsilon_{n_k}, t, Y, dX)(\omega) \underset{k \rightarrow \infty}{\longrightarrow} \left(\int_{]0,\,t]} Y_s \,d^{-}X_s\right)(\omega) \quad \forall t \in [0,\,T],\,\,\forall \omega \notin \mathcal N. 
		\end{eqnarray} 
		\item[2)] If the covariation between $X$ and $Y$  exists, then it exists in the pathwise sense. 
		In particular, there is a null set $\mathcal N$ and, for any sequence $(\varepsilon_n)\downarrow 0$,  a subsequence 
		$(\varepsilon_{n_k})$  such that 
		\begin{eqnarray}\label{conv_ucp_cov} 
		C_{\varepsilon_{n_k}}(X,Y)(t)(\omega) \underset{k \rightarrow \infty}{\longrightarrow} 
		\left[X,Y\right]_t(\omega) \quad \forall t \in [0,\,T],\,\,\forall \omega \notin \mathcal N. 
		\end{eqnarray} 
	\end{itemize} 
\end{corollary} 
\proof 
The result is a direct application of Proposition \ref{P_conv_as}. 
\endproof 
\begin{lemma}\label{L_ucp_conv} 
	Let $g: [0,\,T] \rightarrow \R$ be a c\`agl\`ad process, $X$ be a  c\`adl\`ag process such that 
	the quadratic variation of $X$ 
	exists in the pathwise sense, see Definition \ref{D_cov}.   Setting (improperly) $[X,X] = \Gamma$, we have 
	\begin{equation}\label{EAS} 
	\int_{0}^{s}g_t\,(X_{(t + \varepsilon)\wedge s}-X_t)^2\,\frac{dt}{\varepsilon} \overset{\varepsilon \rightarrow 0}{\longrightarrow}\int_{0}^{s}g_{t}\,d[X,X]_t\quad \textup{u.c.p.} 
	\end{equation} 
\end{lemma} 
\proof 
We have to prove that 
\begin{equation}\label{P_conv} 
\sup_{s \in [0,\,T]}\bigg|\int_{0}^{s}g_t\,(X_{(t + \varepsilon)\wedge s}-X_t)^2\,\frac{dt}{\varepsilon} -\int_{0}^{s}g_{t}\,d[X,X]_t\bigg|\overset{P}{\longrightarrow}0 \quad \text{as $\varepsilon$ goes to zero}. 
\end{equation} 
Let $\varepsilon_n$ be a sequence converging to zero. Since  $[X,X]$ 
exists in the pathwise sense, there is a subsequence $\varepsilon_{n_k}$, that we still symbolize by $\varepsilon_n$, such that 
\begin{equation}\label{DDD} 
C_{\varepsilon_{n}}(X,X)(t) \, \limit^{n \rightarrow \infty} \,[X,X]_t,\quad \forall t \in [0,\,T]\,\,\textup{a.s.} 
\end{equation} 
Let $\mathcal N$ be a null set such that 
\begin{equation}\label{mathcalN} 
C_{\varepsilon_{n}}(X,X)(\omega,t)  \overset{n \rightarrow \infty}{\longrightarrow} [X,X]_t(\omega)\quad \forall t \in [0,\,T],\,\,\forall \omega \notin \mathcal N. 
\end{equation} 
From here on we fix  $\omega \notin  \mathcal N$. 
We have to prove that 
\begin{equation}\label{P_conv_n} 
\sup_{s \in [0,\,T]}\bigg|\int_{0}^{s}g_t\,(X_{(t + \varepsilon_n)\wedge s}-X_t)^2\,\frac{dt}{\varepsilon_n} -\int_{0}^{s}g_{t}\,d[X,X]_t\bigg|\, \limit^{n \rightarrow \infty} \, 0. 
\end{equation} 
We will do it in two steps. 
 
\emph{Step 1.} We consider first the case of a c\`agl\`ad process  $(g_t)$ with  a finite number of jumps. 
 
Let us fix $\gamma >0$, $\varepsilon >0$.   
We enumerate by $(t_i)_{i \geq 0}$ 
the set of jumps of $X(\omega)$  on $[0,\,T]$, 
 union $\{T\}$. 
Without restriction of generality, we will assume that the jumps of $(g_t)$ are included in $\{t_i\}_{i \geq 0}$. 
Let $N = N(\omega)$ such that 
\begin{equation} 
\sum_{i = N+1}^{\infty} |\Delta X_{t_i}|^2 \leq \gamma^2, \quad \sum_{i = N+1}^{\infty} |\Delta g_{t_i}| =0.\label{inequalityGamma} 
\end{equation} 
We define $A(\varepsilon,N)$ and $B(\varepsilon,N)$ as in \eqref{Aepsilon}-\eqref{Bepsilon}.
The term inside the supremum in \eqref{P_conv} can be written as 
\begin{displaymath} 
\frac{1}{\varepsilon} \int_{]0,\,s]} g_t\,(X_{(t + \varepsilon)\wedge s}-X_t)^2\,dt-\int_{]0,\,s]} g_t\,d [X,X]_t = J_{1}(s,\,\varepsilon)+J_{2}(s,\,\varepsilon)+J_{3}(s,\,\varepsilon), 
\end{displaymath} 
where 
\begin{align*} 
J_{1}(\varepsilon,\,N,\,s)&=\frac{1}{\varepsilon} \int_{]0,\,s]\cap A(\varepsilon,N)} g_t\,(X_{(t + \varepsilon)\wedge s}-X_t)^2\,dt-\sum_{i=1}^N \one_{]0,\,s]}(t_i)\,(\Delta X_{t_i})^2\,g_{t_{i}},\\ 
J_{2}(\varepsilon,\,N,\,s)&=\frac{1}{\varepsilon} \int_{]0,\,s]\cap B(\varepsilon,N)} g_t\,(X_{t + \varepsilon}-X_t)^2\,dt 
-\int_{]0,\,s]} g_t\,d [X,X]^c_t-\sum_{i=N+1}^\infty \one_{]0,\,s]}(t_i)\,(\Delta X_{t_i})^2\,g_{t_{i}},\\ 
J_{3}(\varepsilon,\,N,\,s)&=\frac{1}{\varepsilon} \int_{]0,\,s]\cap B(\varepsilon,N)} g_t\,\left[(X_{(t + \varepsilon)\wedge s}-X_t)^2-(X_{t + \varepsilon}-X_t)^2\right]\,dt. 
\end{align*} 
Applying Lemma \ref{L_ucp_big_jumps}  to $J_{1}(\varepsilon,\,N,\,s)$, with $Y=(Y^1,Y^2)=(t,X)$ and $\phi(y_1,y_2)=g_{y_2^1}(y_1^2-y_2^2)^2$, we get 
\begin{equation}\label{C1} 
\lim_{\varepsilon \rightarrow 0}\sup_{s \in [0,\,T]} |J_{1}(\varepsilon,\,N,\,s)| = 0. 
\end{equation} 
Concerning $J_{3}(\varepsilon,N,s)$, we have 
\begin{align*} 
|J_{3}(\varepsilon,N,s)| 
&\leq \frac{||g||_{\infty}}{\varepsilon}\left( 
\int_{s-\varepsilon}^s \one_{B(\varepsilon,N)}(t)\,(|X_{t + \varepsilon}-X_{t}|^2 + |X_{s}-X_{t}|^2) \,\frac{dt}{\varepsilon} 
\right). 
\end{align*} 
We recall that 
$
B(\varepsilon,N) = \bigcup_{i = 1}^N ]t_{i-1},\,t_i - \varepsilon]$.
From Remark \ref{R_Billingsley} it follows that, for every $t \in ]t_{i-1},\,t_i - \varepsilon]$ and $s >t$,  $(t+\varepsilon)\wedge s \in [t_{i-1},\,t_i]$. Therefore Lemma \ref{Lem_Billingsley} applied successively to the intervals $[t_{i-1},\,t_i]$ implies that 
\begin{equation}\label{C3} 
\limsup_{\varepsilon \rightarrow 0}\sup_{s \in [0,\,T]}|J_{3}(\varepsilon,\,N,\,s)|\leq   18 \gamma^2\,||g||_{\infty}. 
\end{equation} 
%
It remains to evaluate the uniform limit of 	$J_{2}(\varepsilon_n,N,s)$. We show, in a first moment, that, 
for fixed $s \in [0,\,T]$, we have the pointwise convergence 
\begin{align} 
\label{DDE} 
J_{2}(\varepsilon_n,N,s)&=\frac{1}{\varepsilon_n} \int_{]0,\,s]\cap B(\varepsilon_n,N)} g_t\,(X_{t + \varepsilon_n}-X_t)^2\,dt 
-\int_{]0,\,s]} g_t\,d [X,X]^c_t-\sum_{i=N+1}^\infty \one_{]0,\,s]}(t_i)\,(\Delta X_{t_i})^2\,g_{t_{i}}\nonumber\\ 
&\underset{n \rightarrow \infty}{\rightarrow}\,0, \quad \forall s \in [0,\,T]. 
\end{align} 
For this, it will be useful to show that
\begin{equation}\label{wek_convB} 
\frac{dt}{\varepsilon_{n}} \one_{B(\varepsilon_n, N)}(t)\,(X_{t + \varepsilon_{n}}-X_{t})^2 
\Rightarrow d\bigg(\sum_{\underset{i=N+1}{t_i \leq t}}^\infty (\Delta X_{t_i})^2 + [X,X]^c_t\bigg). 
\end{equation} 
It will be enough to show  that, $\forall s \in [0,\,T]$, 
\begin{equation}\label{B39bis} 
\int_{0}^{s}  \frac{dt}{\varepsilon_{n}} \one_{B(\varepsilon_n, N)}(t)\,(X_{t + \varepsilon_{n}}-X_{t})^2 
\rightarrow_{n\rightarrow\infty} \sum_{\underset{i=N+1}{t_i \leq s}}^\infty (\Delta X_{t_i})^2 + [X,X]^c_s. 
\end{equation} 
By \eqref{conv_ucp_cov} in Corollary \ref{C_id_RV_fw}-2)
and Lemma \ref{L_bracket_jumps},  we have 
\begin{equation}\label{to_substract} 
\int_{0}^{s}(X_{t + \varepsilon_n}-X_t)^2\,\frac{dt}{\varepsilon_n} \overset{n \rightarrow \infty}{\longrightarrow} [X,X]^c_s + \sum_{t_i\leq s} (\Delta X_{t_i})^2 \quad \forall s \in [0,\,T]. 
\end{equation} 
On the other hand, we can show that  
\begin{equation}\label{A18_BIS} 
\int_0^s \frac{dt}{\varepsilon_{n}}  \one_{A(\varepsilon_n, N)}(t)\,(X_{t + \varepsilon_{n}}-X_{t})^2 \overset{n \rightarrow \infty}{\longrightarrow} \sum_{\underset{i=1}{t_i \leq s}}^N (\Delta X_{t_i})^2 \quad \forall s \in [0,\,T]. 
\end{equation} 
Indeed 
\begin{align*} 
&\Big|\int_0^s \frac{dt}{\varepsilon_{n}}  \one_{A(\varepsilon_n, N)}(t)\,(X_{t + \varepsilon_{n}}-X_{t})^2 - \sum_{\underset{i=1}{t_i \leq s}}^N (\Delta X_{t_i})^2\Big| \\ 
&\leq \Big|\int_0^s \frac{dt}{\varepsilon_{n}}  \one_{A(\varepsilon_n, N)}(t)\,(X_{(t + \varepsilon_{n}) \wedge s}-X_{t})^2 - \sum_{\underset{i=1}{t_i \leq s}}^N (\Delta X_{t_i})^2\Big|\\ 
&+ \Big|\int_0^s \frac{dt}{\varepsilon_{n}}  \one_{A(\varepsilon_n, N)}(t)\,(X_{t + \varepsilon_{n}}-X_{t})^2-\int_0^s \frac{dt}{\varepsilon_{n}}  \one_{A(\varepsilon_n, N)}(t)\,(X_{(t + \varepsilon_{n}) \wedge s}-X_{t})^2\Big| \quad \forall s \in [0,\,T]. 
\end{align*} 
The first addend  converges to zero by Lemma 
\ref{L_ucp_big_jumps} applied to  $Y=X$ and $\phi(y)=(y_1-y_2)^2$. 
The second one converges to zero by similar arguments as those we have used to prove  
Proposition 
\ref{P_conv_as}. This establishes \eqref{A18_BIS}. 
Subtracting \eqref{to_substract} and \eqref{A18_BIS}, we get \eqref{B39bis}, 
and so 
\eqref{wek_convB}.

We remark that the left-hand side of \eqref{wek_convB} are positive measures. 
Moreover,  
we notice that, since the jumps of $g$ are included in $\{t_1,..., t_N\}$,  
$t \mapsto g_t(\omega)$ is $\mu$-continuous, where $\mu$ is the measure on 
the right-hand side of \eqref{wek_convB}. 
At this point,  
Portmanteau theorem  and \eqref{wek_convB} insure that $J_{2}(\varepsilon_n,N,s)$ 
converges to zero as $n$ goes to infinity, for every $s \in [0,\,T]$. 

Finally, we control the convergence  of $J_2(\varepsilon_n, N, s)$,  uniformly in $s$. 
We   make use of Lemma \ref{L_Fn_F_G}. 
We set 
\begin{eqnarray*} 
	G_n(s)&=& \frac{1}{\varepsilon_n}\int_{]0,\,s]}\one_{B(\varepsilon_n,N)}(t)\,(X_{t+ \varepsilon_n}-X_t)^2\,g_t\,dt,\\ 
	F(s)&=& \int_{]0,\,s]}g_{t}\,d[X,X]^c_t,\\ 
	G(s)&=& - \sum_{i=N+1}^\infty \one_{]0,\,s]}(t_i)\,(\Delta X_{t_i})^2\,g_{t_{i}}. 
\end{eqnarray*} 
By  \eqref{DDE}, $F_n := G_n +G$ converges pointwise to $F$ as $n$ goes to infinity. 
%
%
Since $G_n$ is continuous and increasing, $F$ is continuous and $G$ is c\`adl\`ag, Lemma \ref{L_Fn_F_G} implies that 
\begin{equation}\label{C2} 
\limsup_{n \rightarrow \infty}\sup_{s \in [0,\,T]}|J_{2}(\varepsilon_{n},\,N,\,s)|\leq 2 \gamma^2\,||g||_{\infty}. 
\end{equation} 
 
Collecting  \eqref{C1}, \eqref{C3} and \eqref{C2}, it follows that 
\begin{equation*} 
\limsup_{n \rightarrow \infty}\sup_{s \in [0,\,T]} 
\bigg| \int_0^s g_t\,(X_{(t + \varepsilon_{n})\wedge s}-X_{t})^2 \,\frac{dt}{\varepsilon_{n}}-\int_0^s g_t\,d[X,X]_t \bigg| 
\leq 
20\gamma^2 \,||g||_{\infty}. 
\end{equation*} 
Since $\gamma$ is arbitrarily small, \eqref{P_conv_n} follows. 
 
\emph{Step 2.} We treat now the case of a general   c\`agl\`ad process $(g_t)$. 
 
Let us fix $\gamma >0$, $\varepsilon >0$.  
Without restriction of generality, we can write 
$ 
g_{t}= g^{\gamma, BV}_{t} + g^{\gamma}_{t}, 
$ 
where $g^{\gamma,BV}_{t}$ is a process with a finite number of jumps and 
$g^{\gamma}_{t}$ is such that $|\Delta g^{\gamma}_t| \leq \gamma$ for every $t \in [0,\,T]$. 
From Step 1, we have 
\begin{align}\label{unif_est_gammaBV} 
I^{1,n}_s:=\int_0^s g^{\gamma,BV}_t\,(X_{(t + \varepsilon_{n})\wedge s}-X_{t})^2 \,\frac{dt}{\varepsilon_{n}}-\int_0^s g^{\gamma,BV}_t\,d[X,X]_t 
\end{align} 
converges to zero, uniformly in $s$, as $n$ goes to infinity.
 
Concerning $(g^{\gamma}_{t})$, by Lemma \ref{Lem_Billingsley}  we see that  there exists $\bar \varepsilon_0 = \bar\varepsilon_0(\gamma)$ such that 
\begin{equation}\label{Billg} 
\sup_{\underset{|a-t| \leq \bar \varepsilon_0}{a,\,t \in I}}|g^\gamma_a - g^\gamma_t|\leq 3 \gamma. 
\end{equation} 
At this point, we introduce the c\`agl\`ad  process 
\begin{equation}\label{def_gn} 
g^{k,\gamma}_t= \sum_{i = 0}^{2^k-1}\,g^{\gamma}_{i\, 2^{-k} T}\,\one_{]i 2^{-k}T,(i+1) 2^{-k}T]}(t), 
\end{equation} 
where  $k$ is   such that $2^{-k} < \bar \varepsilon_0$. From  \eqref{def_gn}, taking into account \eqref{Billg}, for every $t \in [0,\,T]$ there is $i \in \{0,..., 2 k -1\}$, such that  
\begin{align}\label{estimate_ggamma} 
|g^{\gamma}_t- g^{k,\gamma}_t|= |g^{\gamma}_t\,\one_{]i 2^{-k}\,T,(i+1) 2^{-k}\,T]}(t)-g^{\gamma}_{i\, 2^{-k}}|\leq 3 \gamma.
\end{align} 
We set 
\begin{align*} 
I^{2,n}_s&:=\int_0^s (g^{\gamma}_t- g^{k,\gamma}_t)\,(X_{(t + \varepsilon_{n})\wedge s}-X_{t})^2 \,\frac{dt}{\varepsilon_{n}}-\int_0^s (g^{\gamma}_t- g^{k,\gamma}_t)\,d[X,X]_t. 
\end{align*} 
From \eqref{estimate_ggamma} we have 
$\sup_{s \in [0,\,T]}|I^{2,n}_s|\leq 3 \gamma \, \Gamma$, 
%
with 
\begin{equation}\label{Gamma} 
\Gamma = \sup_{n \in \N, s \in [0,\,T]} \bigg|\int_0^s (X_{(t + \varepsilon_{n})\wedge s}-X_{t})^2 \, 
\frac{dt}{\varepsilon_{n}}\bigg|+[X,X]_T. 
\end{equation} 
Notice that $\Gamma$ is finite, since the term inside the absolute value in \eqref{Gamma} converges uniformly by Step 1 with $g =1$. 
On the other hand, by definition, $(g^{k,\gamma}_t)$ has a finite number of jumps, therefore from Step 1 we get that 
\begin{align}\label{unif_est_gamman} 
I^{3,n}_s&=\int_0^s g^{k,\gamma}_t\,(X_{(t + \varepsilon_{n})\wedge s}-X_{t})^2 \,\frac{dt}{\varepsilon_{n}}-\int_0^s g^{k,\gamma}_t\,d[X,X]_t 
\end{align} 
converges to zero, uniformly in $s$, as $n$ goes to infinity. 
Finally, collecting all the terms, we have 
\begin{align}\label{unif_est_final} 
&\limsup_{n \rightarrow \infty}\sup_{s \in [0,\,T]} 
\bigg| \int_0^s g_t\,(X_{(t + \varepsilon_{n})\wedge s}-X_{t})^2 \,\frac{dt}{\varepsilon_{n}}-\int_0^s g_t\,d[X,X]_t \bigg|\nonumber\\ 
&\leq \limsup_{n \rightarrow \infty}\sup_{s \in [0,\,T]}|I^{1,n}_s|+ \limsup_{n \rightarrow \infty} 
\sup_{s \in [0,\,T]}|I^{2,n}_s|+ \limsup_{n \rightarrow \infty}
\sup_{s \in [0,\,T]}|I^{3,n}_s| 
\nonumber\\ 
& \leq  3\,\gamma \Gamma.
\end{align} 
Since $\gamma$ is arbitrarily small,  \eqref{P_conv_n}  follows. 
\endproof

\begin{proposition}\label{P_equiv_bracket} 
	Let $X, Y$ be two c\`adl\`ag processes. The following properties are equivalent. 
	\begin{itemize} 
		\item[(1)] $[X,X]$, $[X,Y]$, $[Y,Y]$  exist in the pathwise sense.
		\item[(2)] Suppose the existence of processes $\Gamma_1$, $\Gamma_2$, $\Gamma_3$, fulfilling the following properties. For every c\`agl\`ad process $(g_t)$,\\ 
		\begin{align*} 
		\lim_{\varepsilon \rightarrow 0} \int_0^s g_t 
		\,\frac{(X((t+\varepsilon)\wedge s)-X(t))\,(Y((t+\varepsilon)\wedge s)-Y(t))}{\varepsilon}\,dt 
		&= \int_0^s g_t\,
		d \Gamma_1(t)\quad \textup{u.c.p.,}\\ 
		\lim_{\varepsilon \rightarrow 0} \int_0^s g_t \, 
		\,\frac{(X((t+\varepsilon)\wedge s)-X(t))^2}{\varepsilon}\,dt 
		&= \int_0^s g_t\,
		d \Gamma_2(t)
		\quad \textup{u.c.p.,}\\ 
		\lim_{\varepsilon \rightarrow 0} \int_0^s g_t \, 
		\frac{(Y((t+\varepsilon)\wedge s)-Y(t))^2}{\varepsilon}\,dt 
		&= \int_0^s g_t\,
		d \Gamma_3(t)\quad \textup{u.c.p.} 
		\end{align*} 
		In particular, setting $g=1$, we have $\Gamma_1 =[X,Y]$, $\Gamma_2 =[X,X]$, $\Gamma_3 =[Y,Y]$.
	\end{itemize} 
\end{proposition} 
\proof 
Without loss of generality, we first  reduce to the case $g \geq 0$.  Using polarity arguments of the type 
\begin{align*} 
[X+Y,X+Y]_t &= [X,X]_t+[Y,Y]_t + 2\,[X,Y]_t\\ 
[X+Y,X+Y]^{ucp}_{\varepsilon}(t) &= [X,X]^{ucp}_{\varepsilon}(t)+[Y,Y]^{ucp}_{\varepsilon}(t) + 2\,[X,Y]^{ucp}_{\varepsilon}(t), \\
C_{\varepsilon}(X+Y,X+Y)(t) &= C_{\varepsilon}(X,X)(t) +C_{\varepsilon}(Y,Y)(t) + 2\,C_{\varepsilon}(X,Y)(t), 
\end{align*} 
we can reduce to the case $X=Y$. 
 
(1) implies (2) by Lemma \ref{L_ucp_conv}. 
 
(1) follows from (2) choosing $g=1$ and  Corollary \ref{C_id_RV_fw}-2). 
 
\endproof 


\begin{remark}\label{R_equiv_mutual_brackets} 
	Let $X,Y$ be two c\`adl\`ag processes. Implication  (1) $\Rightarrow$ (2) in Proposition  \ref{P_equiv_bracket} with $g=1$, together with Corollary \ref{C_id_RV_fw}-2), implies that the following are equivalent: 
	\begin{itemize} 
		\item 
		$(X, Y)$ admits all its mutual brackets; 
		\item $[X,X]$, $[X,Y]$, $[Y,Y]$   exist in the pathwise sense. 
	\end{itemize} 
	In that case, the covariation processes  above equal the corresponding processes $\Gamma$  related to Definition \ref{D_cov}.
	In particular, the following properties are equivalent: 
	\begin{itemize} 
		\item 
		$X$ is a finite quadratic variation process; 
		\item $[X,X]$  exists in the pathwise sense. 
	\end{itemize} 
\end{remark} 

 
\begin{proposition}\label{P_WDP_finte_quad_var} 
	Let $X$ be a finite quadratic variation process. The following are equivalent. 
	\begin{itemize} 
		\item[(i)]  $X$ is a weak Dirichlet process; 
		\item[(ii)] $X=M+A$,\quad $[A,N]=0$ in the pathwise sense  for all $N$ continuous local martingale. 
	\end{itemize} 
\end{proposition} 
\proof 
(i) $\Rightarrow$ (ii) obviously. 
Assume now that (ii) holds. Taking into account Corollary \ref{C_id_RV_fw} 2), 
it is enough to prove that $[A,N]$ exists. 
Now, we recall that, whenever $M$ and $N$ are local martingale, $[M,N]$  exists by 
Proposition  \ref{P_Ito_cov}. 
Let $N$ be a continuous local martingale. 
By Remark 
\ref{R_equiv_mutual_brackets}, $[X,X]$ and $[N,N]$ exist in the pathwise sense. 
By additivity and item (ii), $[X,N]=[M,N]$  exists in the pathwise sense. 
Still by Remark \ref{R_equiv_mutual_brackets}, $(X,N)$ admits all its mutual brackets. 
Finally, by bilinearity 
$$ 
[A,N]=[X,N]-[M,N]=0. 
$$
\endproof 
 
\makeatletter 
\@addtoreset{equation}{subsection} 
\def\theequation{\thesection.\arabic{equation}} 
\makeatother 
 
\theoremstyle{plain} 
\newtheorem{Theorem}{Theorem}[subsection] 
\newtheorem{Lemma}[Theorem]{Lemma} 
\newtheorem{Proposition}[Theorem]{Proposition} 
\newtheorem{Corollary}[Theorem]{Corollary} 
 
\newtheorem{Definition}[Theorem]{Definition} 
\newtheorem{Hypothesis}[Theorem]{Hypothesis} 
\theoremstyle{remark} 
\newtheorem{Remark}[Theorem]{Remark} 
\newtheorem{Example}[Theorem]{Example}

 \small
\paragraph{Acknowledgements.} 
The authors are grateful to the Referee,  the Associated Editor
and  Editor in Chief for the comments on the first version of the paper
which have motivated them to significantly improve the paper.
The first named author benefited from the support of the Italian
 MIUR-PRIN 2010-11 ``Evolution differential problems: deterministic and stochastic approaches and their interactions''.
The paper was  partially written during the stay of the
second named author  at Bielefeld University, SFB 701 (Mathematik).
 
\addcontentsline{toc}{chapter}{Bibliography} 
\bibliographystyle{plain} 
\bibliography{BiblioLivreFRPV_TESI} 

\end{document}